\documentclass[11pt]{article}
\usepackage{graphicx}
\usepackage{amssymb}
\usepackage{epstopdf}
\begin{document}

\setcounter{secnumdepth}{5}
 
\newcommand{\R}{I\hspace{-1.5 mm}R}
 \newcommand{\N}{I\hspace{-1.5 mm}N}
\newcommand{\1}{I\hspace{-1.5 mm}I}
\newtheorem{proposition}{Proposition}[section]
\newtheorem{theorem}{Theorem}[section]
\newtheorem{lemma}[theorem]{Lemma}
\newtheorem{coro}[theorem]{Corollary}
\newtheorem{remark}[theorem]{Remark}
\newtheorem{ex}[theorem]{Example}
\newtheorem{claim}[theorem]{Claim}
\newtheorem{conj}[theorem]{Conjecture}
\newtheorem{definition}[theorem]{Definition}
\newtheorem{application}{Application}
 
\newcommand{\braket}[2]{\langle #1,#2 \rangle}

\newtheorem{corollary}[theorem]{Corollary}

\def\LX{{\cal L}(X)}
\def\LY{{\cal L}(Y)}
\def\LH{{\cal L}(H)}
 \def\ASD{{\cal L}_{\rm AD}(X)}
 \def\ASDY{{\cal L}_{\rm AD}(Y)}
\def\ASDH{{\cal L}_{\rm AD}(H)}
 \def\ASDP{{\cal L}^{+}_{\rm AD}(X)}
  \def\ASDYP{{\cal L}^{+}_{\rm AD}(Y)}
   \def\ASDHP{{\cal L}^{+}_{\rm AD}(H)}
 \def\CX{{\cal C}(X)}
\def\CY{{\cal C}(Y)}
\def\CH{{\cal C}(H)}
 \def\PX{{\cal A}(X)}
\def\PY{{\cal A}(Y)}
\def\PH{{\cal A}(H)}
\def\phi{{\varphi}}
\def\AH{A^{2}_{H}}
\newcommand{\la}{\lambda}

\title{Anti-selfdual Lagrangians: Variational resolutions of non self-adjoint equations and dissipative evolutions}
\author{ Nassif  Ghoussoub\thanks{Research partially supported by a grant
from the Natural Sciences and Engineering Research Council of Canada. The
author gratefully acknowledges the hospitality and support of the Centre
de Recherches Math\'ematiques in Montr\'eal where this work was initiated and the Universit\'e de Nice where it was completed.} 
\\
\small Department of Mathematics,
\small University of British Columbia, \\
\small Vancouver BC Canada V6T 1Z2 \\
\small {\tt nassif@math.ubc.ca} 
\\
%\today\\ 
\date{January 20, 2005}\\
}
\maketitle

\section*{Abstract} We develop the concept and the calculus of {\it anti-self dual (ASD) Lagrangians} which seems inherent to many questions in mathematical physics, geometry,  and differential equations. They are natural extensions of gradients of convex functions --hence of self-adjoint positive operators-- which usually drive dissipative systems, but also rich enough to provide representations for the superposition of such gradients with skew-symmetric operators which normally generate unitary flows. They yield variational formulations and resolutions for large classes of non-potential boundary value problems and initial-value parabolic equations. Solutions are minima of functionals of the form $I(u)=L(u, \Lambda u)$ (resp. $I(u)=\int_{0}^{T}L(t, u(t), \dot u (t)+\Lambda_{t}u(t))dt$) where $L$ is an anti-self dual Lagrangian and where $\Lambda_{t}$ are  essentially skew-adjoint operators. However, and just like the self (and antiself) dual equations of  quantum field theory (e.g. Yang-Mills) the equations associated to such minima are not derived from the fact they are critical points of the functional $I$, but because they are also zeroes of the Lagrangian $L$ itself. The approach has many advantages:  It solves variationally many equations and systems that cannot be obtained as Euler-Lagrange equations of action functionals, since they can involve non self-adjoint or other non-potential operators;  It also associates variational principles to variational inequalities, and to various dissipative initial-value first order parabolic problems.  These equations can therefore now be analyzed with the full range of methods --computational or not-- that are available for variational settings.  Most remarkable are the permanence properties that ASD Lagrangians possess making them more  pervasive than expected and quite easy to construct. 
 In this first of a series of  papers, we describe the basic theory of anti-self dual Lagrangians and some of its first applications involving mostly bounded linear operators.  In forthcoming papers, we extend the theory to deal with non bounded operators (\cite{GT2}) and with non-linear maps (\cite{G3}) . 
\newpage
\tableofcontents
\section{Introduction}
Non self-adjoint  problems such as the   {\it transport equation}:
\begin{equation}
\label{transport.eq.1}
 \left\{ \begin{array}{lcl}
    \hfill  -\Sigma_{i=1}^{n}a_{i}\frac{\partial u}{\partial x_{i}} +a_{0}u&=& \beta (u) +f  \hbox{\rm \, on \, $\Omega \subset \R^{n}$  
}
\\
 \hfill  u(x) &=& 0 \quad \quad \quad \quad \hbox{\rm on \quad $\Sigma_{-}$. } 
       \end{array}  \right.
   \end{equation}
 where ${\bf a}=(a_{i})_{i}:\Omega \to {\bf R}^{n}$ is a smooth vector field, $\beta$ is a convex function, $f\in L^{2}(\Omega)$, and where $\Sigma_{-}=\{x\in \partial \Omega;\, {\bf a}(x){\bf \cdot {n}}(x)<0\}$,  ${\bf n}$ being  the outer normal vector, are not of Euler-Lagrange type and their solutions are not normally obtained as critical points of functionals  of the form $\int_{\Omega}F(x, u(x), \nabla u(x)) dx$. 
 Similarly, dissipative initial value problems such as the heat equation or those describing   
 {\it porous media:}
\begin{equation}
   \left\{ \begin{array}{lcl}
   \hfill  Ê Ê\frac {\partial u}{\partial t} &=& \Delta u^m +f  \hbox{\rm \, on \, $\Omega\times
[0,T]$}
\\
   u(0,x) &=& u_0(x) \quad \quad \, \hbox{\rm on \quad $\Omega$, }\nonumber \\
       \end{array}  \right.
   \end{equation}
 are not normally solved by the methods of the calculus of variations since they do not correspond to Euler-Lagrange equations of action functionals of the form 
$\int_{0}^{T}L(t, x(t), \dot x (t) dt$.

However, physicists have managed to obtain variationally many of the basic first order equations of quantum field theory  by minimizing their associated action functionals. These are the celebrated self (antiself) dual equations of Yang-Mills, Seiberg-Witten and Ginzburg-Landau which are not derived from the fact they are critical points (i.e., from the corresponding Euler-Lagrange equations) but from the fact that they are zeros of the Lagrangian itself, which is the case as long as the action functional attains a natural and --a priori-- known minimum (See for example \cite{Jo}). 

From a totally different perspective,  Brezis and Ekeland formulated about 30 years ago in \cite{BE1} an intriguing minimization principle which can be associated to the heat equation and other gradient flows of convex energy functionals.  Again the applicability of their principle was conditional on identifying the minimum value of the functional. 
Later,  Auchmuty (\cite{Au1}, \cite{Au2}) proposed a framework in which he formalizes and generalizes the Brezis-Ekeland procedure in order to apply it to operator equations of non-potential type. However, the applicability of this variational  principle remained conditional on evaluating the minimum value and in most cases could not be used to establish existence and uniqueness of solutions. 

In this paper, we develop a general framework where such variational principles are applicable. It is based on the concept of {\it anti-selfdual (ASD) Lagrangians} which seems inherent to many important boundary value problems as well as parabolic evolution equations.  For such Lagrangians  $L$ and for  skew-adjoint operators $\Lambda_{t}$, solutions are obtained as minima of functionals of the form
\[
 I(u)=L(u, \Lambda u) \quad  \hbox{\rm or \quad $I(u)=\int_{0}^{T}L(t, u(t), \dot u(t)+\Lambda_{t}u(t))dt$}.
 \]
The minimal value will always be zero,  and  the equations associated to such minima are not derived from the fact they are critical points of the functional $I$, but because they are also zeroes of the Lagrangian $L$ itself. More specifically, the solutions will satisfy
  \[
  L(u, \Lambda u)+\langle u,\Lambda u\rangle=0 \quad {\rm and} \quad L(t, u(t), \dot u(t)+\Lambda_{t}u(t))+ \langle u (t),\dot u (t) \rangle=0, 
  \]
 for all time, which is reminiscent of the conservation laws enjoyed by Hamiltonians.   This provides variational formulations and complete proofs for the conditional results of Brezis-Ekeland, Auchmuty and others. 
 
 As importantly, we show that ASD Lagrangians possess remarkable  permanence properties making them more  prevalent than expected and quite easy to construct and/or identify. The variational game changes from  the  analytical proofs of existence of extremals for general Lagrangians, to a more  algebraic search of an appropriate ASD Lagrangian for which the minimization problem is remarkably simple. This makes  them efficient new tools for proving existence and  uniqueness results for a large array of differential equations.

The basic idea is simple and is an elaboration on our work in \cite{GT1} where  we gave  complete variational proofs of the existence and uniqueness of gradient flows of convex energy functionals, and the one in \cite{G1}, where we give a variational proof for the existence and uniqueness of solutions of certain non-linear transport equations. Starting with an equation of the form 
 \begin{equation}
 -Au\in \partial \varphi (u) 
 \end{equation}
it  is well known that it can be formulated --and sometimes solved--  variationally  whenever $A: X\to X^{*}$ is a selfadjoint bounded linear operator and $\varphi$ is an appropriate functional on $X$. Indeed, in this case it can be reduced to the equation $0\in \partial \psi (u)$, where $\psi$ is the functional 
  \begin{equation}
 \psi (u)=\varphi (u) +\frac{1}{2}\langle Au, u\rangle.
  \end{equation}
  A solution can then be obtained by minimization whenever $\varphi$ is convex and lower semi-continuous and whenever $A$ is positive (i.e., $\langle u, Au \rangle \geq 0$)  or better if $A$ is coercive (i.e., if for some $c>0$, $\langle u, Au \rangle \geq c\|u\|^{2}$ for all $u\in X$).  

   But this variational procedure fails when $A$ is not self-adjoint, or when $A$ is a non-potential operator (i.e., when $A$ is not a gradient vector field), and definitely when $A$ is not linear. In this case, the  Brezis-Ekeland procedure --as formalized by Auchmuty-- consists of simply minimizing  the functional 
   \begin{equation}
   I(u)=\varphi (u) +\varphi^{*}(-Au)+\langle u, Au \rangle
   \end{equation}
    where  $\varphi^{*}$ is the Legendre dual of $\varphi$ defined on $X^{*}$ by
  $
    \varphi^{*}(p)=\sup\{\langle x, p\rangle -\varphi (x); \, x\in X\}.
   $
  Legendre duality yields that $\alpha:=\inf_{u\in X}I(u)\geq 0$,  and the key observation made by several authors is the following simple\\
 {\bf Fact:} If the infimum $\alpha=0$ and if it is attained at $\bar u \in X$ 
   then we are in the limiting case of the Fenchel-Legendre duality, 
$
\varphi (\bar u) +\varphi^{*}(-A\bar u)=\langle \bar u, -A\bar u\rangle
$
 and therefore $-A\bar u\in \partial \varphi(\bar u)$.
 
 Note that the procedure does not require any assumption on $A$, and very general coercivity assumptions on $\varphi$ often ensure the existence of a minimum. However, the difficulty here is different from standard minimization problems in that besides the problem of existence of a minimum, one has to insure that the infimum is actually zero. This is obviously not the case for general operators $A$, though one can always write (and many authors did) the variational principle (5) for the operator equation (3). 

In this paper, we tackle the real difficulty of when the infimum $\alpha$ is actually zero and we try to identify a class of nonpotential operators $F(u)$ for which the equation and the initial-value problem 
\begin{equation}
0\in F(u) \quad \quad \quad {\rm and } \quad \quad \quad  \left\{ \begin{array}{lcl}
\label{eqn:laxmilgram999}
 \hfill 
-\dot u(t) &\in& F(u(t))\\
\hfill u(0)&=&u_{0}
 \\
\end{array}\right.
\end{equation}
can be solved by the above variational procedure. We show here that this is essentially the case whenever $F(u)=B u +\partial \varphi(u)$, where $\phi$ is a convex lower semicontinuous function and when $B$ is a skew-adjoint operator.  We note that --when $\Lambda$ is linear-- such operators form a very important subset of the class of maximal monotone operators for which there is already an extensive theory (\cite{Br}, \cite{Bar}). The interest here is in the new variational approach based on the concept of anti-selfdual Lagrangians which possesses remarkable  permanence properties that maximal monotone operators either do not satisfy or do so via substantially more elaborate methods.  In a forthcoming paper (\cite{G3}) we establish  similar results for operators of the form $
F(u)=\Lambda u +Bu +\partial \varphi(u)
$
where $\Lambda$ is an appropriate non-linear conservative operator, $B$ is linear and positive, and $\phi$ is convex, the  superposition of which is not normally covered by the theory of maximal monotone operators.

In this paper, we establish the algebraic structure of ASD Lagrangians, emphasizing issues on how to build and identify complex ASD Lagrangians from the more basic ones. To keep the key ideas transparent, we chose to deal with the case  when the operators are bounded and linear, leaving the more analytically involved cases of unbounded and nonlinear operators to forthcoming papers.  This --bounded linear-- case already has many interesting features, especially in boundary value problems of the form:
 \begin{equation}
 \left\{ \begin{array}{lcl}
\label{eqn:laxmilgram1000}
 \hfill  -Ax+f&\in &\partial \varphi (x)\\
\hfill  B(x)&=&a \\
\end{array}\right.
\end{equation}
  where $B$ is a boundary operator on $X$ (related to the positive operator $A$), as well as  parabolic evolution equations of the form:
 \begin{equation}
 \left\{ \begin{array}{lcl}
\label{eqn:1001}
 \hfill   -A_{t} x(t)-\dot{x}(t)  &\in &\partial \varphi (t, x(t)) \quad \hbox {\rm a.e. $t\in [0, T]$}\\
\hfill B_{t}(x(t))&=&a(t) \quad \quad \quad \quad \hbox{\rm a.e $t\in [0, T]$}\\
\hfill x(0)&=&x_{0}
\end{array}\right.
\end{equation}
where $x_{0}$ is a given initial value and where $a(t)$ is a prescribed boundary value. 

We  start by presenting -- in  section 2-- the special  variational properties of  the class of  $R$-Antiselfdual   Lagrangians, where $R$ is any automorphism of the state space. This should already give an idea of their relevance in the existence theory of certain PDEs, and will hopefully motivate the in-depth study of their permanence properties. Beyond this first section, we will 
only deal with the anti-symmetric case, i.e.,  when $R (x)=-x$, in which case $R$-antiselfdual Lagrangians will be called {\it anti-selfdual Lagrangians (ASD)}. We shall see that this class of Lagrangians already covers a great deal of applications which warranted that this paper as well as (\cite{GT2}, \cite{G3}) be solely devoted to this case. However, the theory involving other automorphisms $R$ will also be very useful, especially in applications to Hamiltonian systems and this will be developed in \cite{GT3}.

In section 3, we  establish  the basic permanence properties of {\it anti-self dual Lagrangians} as well as their special variational features while focussing on stationary equations and systems. This restrictive looking class turns out to be quite rich. In section (3) we deal with boundary value problems where appropriate selfdual boundary Lagrangians are appropriately added to the ``interior Lagrangian'' to make it anti-selfdual allowing us to solve problems with prescribed boundary terms. In section (4), we show how ASD Lagrangians ``lift'' to path spaces allowing us to solve with the same variational approach several parabolic equations --including gradient flows.  In section (5), we associate to each autonomous Lagrangian, a semi-group of contractions which emphasizes again that such Lagrangians are natural extensions of gradients of convex functions, of positive operators as well as of the ''superpositions'' of the two actions. In section 6, we give a glimpse on how the theory can help in solving variationally certain implicit PDEs, a project for future investigation. 

As mentioned above, in this paper we describe the basics of the ASD theory emphasizing its stability under various operations and its rich structure. So we stuck with the simplest of examples leaving more complicated PDE settings to forthcoming papers. In  (\cite{GT2}), we extend the theory to deal with  linear but unbounded operators, and in (\cite{G3}) we tackle various non-linear but appropriately defined ``skew-adjoint'' operators such as those appearing in the Navier-Stokes and other equations of hydrodynamics.    
Finally, I would like to thank Yann Brenier, Eric S\'er\'e, Leo Tzou and Abbas Moameni for the many extremely fruitful discussions and their  valuable input into this project.
 
\section{Basic variational properties of R-antiselfdual Lagrangians}

 We consider the class $\LX$ of  {\it convex Lagrangians} $L$ on a reflexive Banach space $X$: these are  all functions  $L:X\times X^{*} \to \R \cup \{+\infty\}$ which are convex and lower semi-continuous (in both variables) and which are not identically $+ \infty$. The Legendre-Fenchel dual (in both variables) of $L$ is defined at any pair $(q,y)\in X^{*}\times X$ by: 
 \[
    L^*( q,y)= \sup \{ \braket{q}{x} + \braket{y}{p} - L(x,p);\, x \in X, p \in X^{*}  \}
\]
  \begin{definition}\rm Given a bounded linear operator $R:X\to X$, say that:\\
  (1) $L$ is an {\it $R$-antiselfdual  Lagrangian} on $X\times X^{*}$, if  
\begin{equation}
\label{seldual}
L^*( p, x) =L(-Rx, -R^*p ) \quad \hbox{\rm for all $(p,x)\in X^{*}\times X$}. 
\end{equation}
(2) $L$ is  {\it partially $R$-antiselfdual}, if 
\begin{equation}
\label{seldual.atzero}
L^*( 0, x) =L(-Rx, 0 ) \quad \hbox{\rm for all $x\in X$}. 
\end{equation}
(3) $L$ is {\it $R$-antiselfdual on the graph of $\Lambda$}, the latter being a map from $X$ into $X^{*}$, if 
\begin{equation}
\label{seldual.atgraph}
L^*(\Lambda x, x) =L(-Rx, -R^*\circ \Lambda x ) \quad \hbox{\rm for all $x\in X$}. 
\end{equation}
(4) More generally, if $Y\times Z$ is any subset of $X\times X^{*}$, we shall say that $L$ is {\it $R$-antiself dual on the elements of $Y\times Z$} if $L^*( p, x) =L(-Rx, -R^*p )$  for all $(p,x)\in Y\times Z$.  
\end{definition}
A typical example of an $R$-antiselfdual   Lagrangian is $L(x,p)=\phi (R^{-1}x) +\phi^*(-p)$ and $M(x,p)=\phi (-x) +
\phi^*((R^*)^{-1}p)$ where $\phi$ is a convex lower semi-continuous function and $R$ is an invertible operator on $X$. More generally,  $L(x,p)=\phi (Rx) +\phi^*(-S^*p)$ is an $(S\circ R)^{-1}$-antiselfdual Lagrangian. Moreover, if $\Lambda: X\to X^*$ is such that $\Lambda \circ (S\circ R)^{-1}$ is skew-adjoint, then 
\[
L(x,p)=\phi (Rx) +\phi^*(-S^*\Lambda x-S^*p)
\]
is also an $(S\circ R)^{-1}$-antiselfdual Lagrangian.

Our basic premise in this paper is that many boundary value problems  can be solved by minimizing functionals of the form $I(x)=L(x,\Lambda x)$ where $L$ is  a $R$-Antiselfdual   Lagrangian and provided $\Lambda \circ R$ is a skew-adjoint operator. However, their main relevance to our study stems from the fact that --generically-- the infimum is actually equal to $0$.  It is this latter property that allows for novel variational formulations and resolutions of several basic PDEs and evolution equations, which --often because of lack of self-adjointness-- do not normally fit the Euler-Lagrange framework. 
\\
As mentioned above,  if $L$ is a $R$-Antiselfdual   Lagrangian and if  $\Lambda:X\to X^*$ is an operator such that $\Lambda  \circ R$ is skew adjoint, then the Lagrangian $L_{\Lambda}(x,p)=L(x, \Lambda x+p)$ is again $R$-Antiselfdual. In other words,  Minimizing $L(x,\Lambda x)$ amounts to minimizing $L_{\Lambda} (x,0)$ which is covered by the following very simple --yet far reaching-- proposition.  Again, its relevance comes from the evaluation of the minimum and not from the --more standard-- question about its attainability. 
 
 We  start by noticing that for a $R$-Antiselfdual   Lagrangian, we readily have:
 \begin{equation}
 \label{ineq.1}
L(Rx, R^*p)\geq -\langle Rx, p\rangle \quad \hbox{\rm  for every $(x, p) \in X\times X^{*}$}, 
\end{equation}
and if $L$ is partially anti-selfdual, then   
\begin{equation}
 \label{ineq.2}
I(x)=L(Rx,0)\geq 0 \quad \hbox{\rm  for every $x \in X$}, 
\end{equation}
So, we are looking into an interesting variational situation, where the minima can also be zeros of the functionals. Here are some necessary conditions for the existence of such minima.   
 
\begin{proposition} Let $L$ be a convex lower-semi continuous functional on a reflexive Banach space $X\times X^{*}$. Assume that $L$ is  a partially $R$-Antiselfdual   Lagrangian and that for some $x_{0}\in X$, the function  $p\to L(x_{0},p)$ is bounded above on  a neighborhood of the origin in $X^{*}$.
 Then there exists $\bar x\in X$, such that: 
 \begin{equation}
 \left\{ \begin{array}{lcl}
\label{eqn:existence}
  L( -R\bar x, 0)&=&\inf\limits_{x\in X} L(x,0)=0.\\
 \hfill (0, \bar x) &\in & \partial L (-R\bar x,0).
\end{array}\right.
 \end{equation}
 \end{proposition}
\noindent{\bf Proof:} This follows from the basic duality theory in convex optimization. Indeed, if $({\cal P}_{p})$ is the primal  minimization problem
 $h(p)=\inf\limits_{x\in X}L(x,p)$
 in such a way that $({\cal P}_{0})$ is the initial problem $h(0)=\inf\limits_{x\in X}L(x,0)$, 
then the dual problem $({\cal P}^{*})$ is 
$\sup_{y\in X}-L^{*}(0,y)$, 
and we have the weak duality formula
\[
\inf {\cal P}_{0}:=\inf_{x\in X}L(x,0) \geq \sup_{y\in X}-L^{*}(0,y):=\sup{\cal P}^{*}. 
\]
The ``partial $R$-Antiselfdual  ity'' of $L$ gives that 
\begin{equation}
\label{weak.duality}
\inf_{x\in X}L(x,0) \geq \sup_{y\in X}-L^{*}(0,y)=\sup_{y\in X}-L(-Ry,0).
 \end{equation}
Note now that  $h$ is convex on $X^{*}$ and that its Legendre conjugate satisfies  $h^{*}(y)=L^{*}(0,y)=L(-Ry,0)$ on $X$.
 If now  $h$ is subdifferentiable  at $0$ (i.e., if the problem $({\cal P}_{0})$ is stable), then for any $\bar x \in \partial h(0)$, we have $h(0) +h^{*}(\bar x)=0$, which means that
\[
-\inf_{x\in X}L(x,0)=-h(0)=h^{*}(\bar x)=L^{*}(0,\bar x)=L(-R\bar x,0)\geq  \inf_{x\in X}L(x,0).
\]
It follows that $\inf_{x\in X}L(x,0)=L(-R\bar x,0)\leq 0$ and in view of (\ref{ineq.2}), we get that  the infimum of $({\cal P})$  is zero and  attained at $-R\bar x$, while the supremum of $({\cal P}^{*})$ is attained at $\bar x$.
In this case we can write
\[
L(-R\bar x,0)+L^*(0,\bar x)=0
\]
which yields that $(0, \bar x) \in  \partial L (-R\bar x,0)$.

 If now for some $x_{0}\in X$, the function  $p\to L(x_{0},p)$ is bounded above on  a neighborhood of the origin in $X^{*}$,  then  $h(p) \leq \inf\limits_{x\in X}L(x,p) \leq L(x_{0},p)$ and therefore $h$ is subdifferentiable at $0$ and we are done.

\begin{remark}\rm  The above holds under the condition that $x\to L(Rx,0)$ is coercive in the following sense: 
\begin{equation}
\label{coercive.0}
\lim\limits_{\|x\|\to \infty}\frac{L(Rx,0)}{\|x\|}=+\infty.
\end{equation}
Indeed since $h^{*}(y)=L^{*}(0,y)=L(Ry,0)$ on $X$, we get that  that $h^{*}$ is coercive on $X$, which means that $h$ is bounded above on neighborhoods of zero in $X^{*}$.
  \end{remark}

\begin{remark} \rm The proof above requires only that $L$ is a Lagrangian satisfying 
\begin{equation}
L^{*}(0,x)\geq L(-Rx,0) \geq 0 \quad {\rm for\, all}\,  x\in X.
\end{equation}
 Now we can deduce the following
\begin{theorem} Let $R:X\to X$ be a bounded linear operator on a reflexive Banach space $X$ and let $\Lambda :X\to X^{*}$ be another operator such that  $\Lambda \circ R$ is skew adjoint. Let $L$ be a Lagrangian on $X$ that is $R$-antiselfdual  on the graph of  $-\Lambda^*$, and assume that  $\lim\limits_{\|x\|\to \infty}\frac{L(Rx,\Lambda Rx)}{\|x\|}=+\infty$.
 Then there exists $\bar x\in X$, such that: 
 \begin{equation}
 \left\{ \begin{array}{lcl}
\label{eqn:existence.2}
  L( -R\bar x, -\Lambda R\bar x)&=&\inf\limits_{x\in X} L(x,\Lambda x)=0.\\
 \hfill (-\Lambda^* \bar x, \bar x) &\in & \partial L (-R\bar x,-\Lambda R\bar x).
\end{array}\right.
 \end{equation}
 \end{theorem}
 \noindent{\bf Proof:} We first prove that the Lagrangian defined as $M(x,p)=L(x, \Lambda x+p)$ is partially $R$-Antiselfdual. 
Indeed fix $(q, y)\in X^{*}\times X$, set $r=\Lambda x+p$ and write:
   \begin{eqnarray*}
  M^{*} (q,y)&=&\sup\{\langle q, x\rangle +  \langle y, p\rangle- L(x, \Lambda x+p); (x,p)\in X\times X^{*}\}\\
  &=& \sup\{\langle q, x\rangle +  \langle y, r-\Lambda x\rangle- L(x, r); (x,r)\in X\times X^{*}\}\\
  &=& \sup\{\langle q-\Lambda^* y, x\rangle +  \langle y, r\rangle- L(x, r); (x,r)\in X\times X^{*}\}\\
  &=&  L^{*}(q-\Lambda^* y, y).
 \end{eqnarray*}  
 If $q=0$, then $
 M^{*} (0,y)=L^{*}(-\Lambda^* y, y)=L(-Ry,R^*\Lambda^*y)=L(-Ry, -\Lambda Ry)=M(-Ry,0)$, and $M$ is therefore partially $R$-Antiselfdual.
 
It follows from the previous proposition applied to $M$, that there exists $\bar x\in X$ such that:
\[
 L( -R\bar x, -\Lambda R\bar x)=M(-R\bar x,0)=\inf\limits_{x\in X} M(x,0)=\inf\limits_{x\in X} L(x,\Lambda x)=0.
 \]
Now note that
\[
 L( -R\bar x, -\Lambda R\bar x)=L(-R \bar x,R^*\Lambda^*\bar x)=L^*(-\Lambda^*\bar x, \bar x),
\] 
hence
\[
 L( -R\bar x, -\Lambda R\bar x)+L^*(-\Lambda^*\bar x, \bar x)=0=\langle (-R\bar x,-\Lambda R\bar x), (-\Lambda^* \bar x, \bar x)\rangle.
\]
It follows from the limiting case of Legendre duality  that $ (-\Lambda^* \bar x, \bar x) \in  \partial L (-R\bar x,-\Lambda R\bar x)$.
 
\end{remark}

\section{Permanence properties of Anti-selfdual Lagrangians}

The concept of $R$-Antiselfduality for a general automorphism $R$  is relevant for dealing with certain Hamiltonian systems  \cite{GT3} and will be pursued in full generality  in a forthcoming paper \cite{G4}. We shall however concentrate in the sequel  on the class of  {\it anti-selfdual Lagrangians (ASD)}, meaning those $R$-Antiselfdual   Lagrangians corresponding to the identity operator $R(x)=x$. In other words, \\
 (1)  $L$ is said to be an {\it anti-selfdual Lagrangian} on $X\times X^{*}$, if  
\begin{equation}
\label{seldual}
L^*( p, x) =L(-x, -p ) \quad \hbox{\rm for all $(p,x)\in X^{*}\times X$}. 
\end{equation}
(2) $L$ is  {\it partially anti-self dual}, if 
\begin{equation}
\label{seldual.atzero}
L^*( 0, x) =L(-x, 0 ) \quad \hbox{\rm for all $x\in X$}. 
\end{equation}
(3) $L$ is {\it anti-self dual on the graph of $\Lambda$}, the latter being a map from $X$ into $X^{*}$, if 
\begin{equation}
\label{seldual.atgraph}
L^*(\Lambda x, x) =L(-x, -\Lambda x ) \quad \hbox{\rm for all $x\in X$}. 
\end{equation}
(4) More generally, if $Y\times Z$ is any subset of $X\times X^{*}$, we shall say that $L$ is {\it anti-self dual on the elements of $Y\times Z$} if $L^*( p, x) =L(-x, -p )$  for all $(p,x)\in Y\times Z$. \\
 
 Denote by  $\ASD$ the class of anti-selfdual (ASD) Lagrangians on a given Banach space $X$. We shall see that this is already a very interesting and natural class of Lagrangians as they appear in several basic PDEs and evolution equations. The basic example of an anti-selfdual Lagrangian is given by a function $L$ on $X\times X^{*}$, of the form
\begin{equation}
L(x,p)=\varphi (x) +\varphi^{*}(-p)
\end{equation}
where $\varphi$ is a convex and lower semi-continuous function on $X$ and $\varphi^{*}$ is its Legendre conjugate on $X^{*}$.  We shall call them the {\it Basic ASD-Lagrangians}. A key element of this theory is that the family of ASD Lagrangians is much richer and goes well beyond convex functions and their conjugates, since they are naturally compatible with skew-symmetric operators.  Indeed if $\Lambda: X\to X^{*}$ is skew-symmetric  (i.e., $\Lambda^{*}=-\Lambda$),  the Lagrangian 
\begin{equation}
M(x,p)=\varphi (x) +\varphi^{*}(-\Lambda x-p)
\end{equation}
is also anti-self dual, and if in addition $\Lambda$ is invertible then the same holds true for 
\begin{equation}
N(x,p)=\varphi (x+\Lambda^{-1}p) +\varphi^{*}(\Lambda x).
\end{equation}

\subsection*{Basic properties of ASD Lagrangians}
The class $\ASD$ enjoys a remarkable number of  permanence properties. Indeed,  we define on the class of Lagrangians $\LX$ the following operations: 
\begin{itemize}
 \item {\bf Scalar multiplication:} If $\lambda>0$ and $L\in \LX$, define the Lagrangian $\lambda {\bf \cdot} L$ on $X\times X^{*}$ by:
 \[
( \lambda {\bf \cdot} L)(x, p)=\lambda^{2}L(\frac{x}{\lambda},  \frac{p}{\lambda}).
 \]
 \item {\bf Addition:} If $L, M\in \LX$, define the Lagrangian $L+M$ on $X\times X^{*}$ by:
 \[
 (L\oplus M)(x,p)=\inf\{L(x, r) + M(x,p-r); r\in X^{*}\}
 \]
\item {\bf Convolution:} If $L, M\in \LX$, define the Lagrangian $L\star M$ on $X\times X^{*}$ by:
\[
(L\star M) (x,p)=\inf\{L(z, p) + M(x-z,p); z\in X\}
\]
\item {\bf Right operator shift:}   If $L\in \LX$ and $\Lambda :X\to X^{*}$ is a bounded linear operator, define the Lagrangian $L_{\Lambda}$ on $X\times X^{*}$ by
\[
L_{\Lambda} (x,p):=L(x, \Lambda x+p).
\] 
 \item {\bf Left operator shift:} If $L\in \LX$ and if $\Lambda :X\to X^{*}$ is an invertible  operator, define the Lagrangian $ {}_{\Lambda}L$ on $X\times X^{*}$ by:
  \[
  {}_{\Lambda}L(x, p):=L(x +\Lambda^{-1}p, \Lambda x).
  \] 
\item {\bf Free product:}\,  If $\{L_{i}; i\in I\}$ is a finite family of Lagrangians on reflexive Banach spaces $\{X_{i}; i\in I\}$,  define the Lagrangian $L:=\Sigma_{i\in I} L_{i}$ on  $(\Pi_{i\in I} X_{i}) \times (\Pi_{i\in I} X^{*}_{i})$ by 
\[
L((x_{i})_{i}, (p_{i})_{i})=\Sigma_{i\in I} L_{i}(x_{i}, p_{i}).
\]
\item {\bf Twisted A-product:} If $L \in \LX$ and $M\in \LY$ where $X$ and $Y$ are two reflexive spaces, then for any bounded linear operator $A: X\to Y^{*}$, define the Lagrangian $L\oplus_{A}M$  on $(X\times Y)\times (X^{*}\times Y^{*})$ by 
\[
(L\oplus_{A}M) ((x,y), (p.q)):=L(x, A^{*}y+p)+M(y, -Ax+q). 
\]
\item {\bf A-antidualisation:} If $\phi$ is any convex function on $X\times Y$ and $A$ is any bounded linear operator $A: X\to Y^{*}$, define the Lagrangian $L\oplus_{\rm as} {A}$  on $(X\times Y)\times (X^{*}\times Y^{*})$ by 
 \[
\phi\oplus_{\rm as} A((x,y), (p.q))=\phi (x, y)+\phi^{*}(-A^{*}y-p, Ax-q).
\]
 \end{itemize}
 The above defined convolution operation should not be confused with the standard convolution for $L$ and $M$ as convex functions in both variables. It is easy to see that in the case where $L(x,p)=\phi (x) +\phi^{*}(-p)$ and $M(x,p)=\psi (x) +\psi^{*}(-p)$, addition corresponds to taking 
\[
(L\oplus M) (x,p)=(\phi +\psi) (x) +\phi^{*}\star \psi^{*}(-p)
\]
while convolution reduces to:
\[ 
(L\star M) (x,p)=(\phi \star\psi) (x) +(\phi^{*}+ \psi^{*})(-p). 
\]
which also means that they are dual operations. We do not know whether this is true in general, but for the sequel we shall only need the following:
 
 \begin{lemma}  Let $X$ be a reflexive Banach space and consider two Lagrangians $L$ and $M$ in $ \LX$. Then the following hold:
 \begin{enumerate}
 \item  If  $\lambda >0$,  then  $(\lambda{\bf \cdot} L) \star M= \lambda {\bf \cdot} (L\star M)$.
 \item $(L\oplus M)^{*} \leq L^{*}\star M^{*}$ and $(L\star M)^{*} \leq L^{*}\oplus M^{*}$.
 \item If $L$ or $M$ is a basic ASD Lagrangian, then $(L\oplus M)^{*} = L^{*}\star M^{*}$ and $(L\star M)^{*} = L^{*}\oplus M^{*}$.
 \item If $L$ and $M$ are in $ \ASD$, then $L^{*}\oplus M^{*}(q,y)=L\star M (-y,-q)$ for every $(y,q)\in X\times X^{*}$. 
 \end{enumerate}
 \end{lemma}
 \noindent{\bf Proof:} (1)  is straightforward. To prove (2), fix $(q, y)\in X^{*}\times X$ and write:
  \begin{eqnarray*}
  &&(L\star M)^{*} (q,y)\\
  &=&\sup\{\langle q, x\rangle +  \langle y, p\rangle- L(z, p)- M(x-z,p); (z, x, p)\in X\times X\times X^{*}\}\\
&=&  \sup\{\langle q, v+z\rangle +  \langle y, p\rangle- L(z, p)- M(v,p); (z, v, p)\in X\times X\times X^{*}\}\\
  &=&  \sup\{\langle q, v+z\rangle +\sup\{\langle y, p\rangle- L(z, p)- M(v,p); p\in X^{*} \}; (z, v)\in X\times X\}\\
   &=&  \sup\limits_{ (z, v)\in X\times X}\left\{\langle q, v+z\rangle +\inf\limits_{w\in X}\{\sup\limits_{p_{1}\in X^{*}}(\langle w, p_{1}\rangle - L(z,p_{1}))+\sup\limits_{p_{2}\in X^{*}}(\langle y-w, p_{2}\rangle- M(v,p_{2}))\}\right\}\\
   &\leq&\inf\limits_{w\in X}\left\{\sup\limits_{(z, p_{1})\in X\times X^{*}}\{ \langle q, z\rangle +\langle w, p_{1}\rangle -L(z,p_{1}))\}+\sup\limits_{(v, p_{2})\in X\times X^{*}}\{\langle q, v\rangle+\langle y-w, p_{2}\rangle - M(v,p_{2})\right\}\\
     &=&  \inf\limits_{w\in X}\left\{L^{*}(q,w)+M^{*}(q, y-w)\right\}\\
 &=& (L^{*}\oplus M^{*})(q, y).
 \end{eqnarray*} 
For (3) assume that $M(x,p)=\phi (x) +\phi^{*}(-p)$ where $\phi$ is a convex lower semi-continuous function. Fix $(q, y)\in X^{*}\times X$ and write:
\begin{eqnarray*}
   (L\star M)^{*} (q,y)
  &=&\sup\{\langle q, x\rangle +  \langle y, p\rangle- L(z, p)- M(x-z,p); (z, x, p)\in X\times X\times X^{*}\}\\
&=&  \sup\{\langle q, v+z\rangle +  \langle y, p\rangle- L(z, p)- M(v,p); (z, v, p)\in X\times X\times X^{*}\}\\
 &=&\sup_{p\in X^{*}} \left\{\langle y,p\rangle +\sup\limits_{(z, v)\in X\times X}\{\langle q,v+z\rangle -L(z,p)-\phi (v)\}\}- \phi^*(-p)\right\}\\
&=&\sup_{p\in X^{*}}\left\{\langle y,p\rangle +\sup_{z\in X}\{\langle q,z\rangle -L(z,p)\} +\sup_{v\in X} \{\langle q,v\rangle -\phi (v)\}-\phi^*(-p)\right\}\\
&=&\sup_{p\in X^{*}} \left\{\langle y,p\rangle +\sup_{z\in X}\{\langle q,z\rangle -L(z,p)\} +\phi^{*}(q)-\phi^*(-p)\right\}\\
&=&\sup_{p\in X^{*}} \sup_{z\in X}\left\{\langle y,p\rangle +\langle q,z\rangle -L(z,p)-\phi^{*}(-p)\right\} +\phi^{*}(q)\\
&=& (L+T)^{*}(q,y)+\phi^{*}(q)
\end{eqnarray*}
where $T(z,p):= \phi^{*}(-p)$ for all $(z,p)\in X\times X^{*}$. Note now that 
\begin{eqnarray*}
T^*(q,y)=\sup_{z,p}\left\{\langle q,z\rangle +\langle y,p\rangle -  \phi^{*}(-p)\right\}=\left\{\begin{array}{lll}+\infty &\hbox{if }&q\ne 0\\ \phi (-y)&\hbox{if }&q=0\end{array} \right.
\end{eqnarray*}
in such a way that by using the duality between sums and convolutions in both variables, we get
\begin{eqnarray*}
(L+T)^{*}(q,y)&=&{\rm conv}(L^*,T^*)(q,y)\\
&=&\inf_{r\in X^{*},z\in X}\left\{ L^*(r,z)+T^*(-r+q,-z+y)\right\}\\
&=&\inf_{z\in X}\left\{ L^*(q,z)+\phi (z-y)\right\}\\
\end{eqnarray*}
and  finally 
\begin{eqnarray*}
(L\star M)^{*} (q,y)&=&(L+T)^{*}(q,y)+ \phi^{*}(q)\\
&=&\inf_{z\in X}\left\{ L^*(q,z)+\phi (z-y)\right\}+ \phi^{*}(q)\\
&=&\inf_{z\in X}\left\{ L^*(q,z)+\phi^{*}(q)+\phi (z-y)\right\}\\
  &=&(L^{*}\oplus M^{*})(q,y).
\end{eqnarray*}
The rest follows in the same way. For (4) write
\begin{eqnarray*}
(L^{*}\oplus M^{*})(q,y)&=&  \inf\limits_{w\in X}\left\{L^{*}(q,w)+M^*(q, y-w)\right\}\\
&=&  \inf\limits_{w\in X}\left\{L(-w,-q)+M(w-y,-q)\right\}\\
&=& (L\star M)(-y, -q).
 \end{eqnarray*} 
The following proposition summarizes some of the remarkable permanence properties of ASD Lagrangians.
\begin{proposition} Let $X$ be a reflexive Banach space, then the following holds:
   \begin{enumerate}
   \item If $L$ is in $\ASD$, then $L^{*} \in {\cal L}_{\rm AD}(X^{*})$, and if $\lambda >0$, then  $\lambda {\bf \cdot} L$  also belong to $\ASD$.
\item  If $L$ and $M$ are in $ \ASD$ and one of them is basic, then the Lagrangians $L\oplus M$, and   $L\star M$  also belong to $\ASD$.  
\item If $L_{i}\in {\cal L}_{\rm AD}(X_{i})$ where $X_{i}$ is a reflexive Banach space for each $i\in I$, then $\Sigma_{i\in I} L_{i}$   is  in ${\cal L}_{\rm AD}(\Pi_{i\in I} X_{i})$.
   \item If $L\in \ASD$ and $\Lambda :X\to X^{*}$ is a skew-adjoint bounded linear operator  (i.e., $\Lambda^{*}=-\Lambda$), then the Lagrangian $L_{\Lambda}$ is also in $\ASD$.
 \item If $L\in \ASD$ and if $\Lambda :X\to X^{*}$ is an invertible  skew-adjoint operator, then the Lagrangian ${}_{\Lambda}L$  is also in $\ASD$. 
 \item If $L \in \ASD$ and $M\in {\cal L}_{\rm AD}(Y)$, then for any bounded linear operator $A: X\to Y^{*}$,  the Lagrangian $L\oplus_{A}M$ belongs to ${\cal L}_{\rm AD}(X\times Y)$
 \item  If $\phi$ is a proper convex lower semi-continuous function on $X\times Y$ and $A$ is any bounded linear operator $A: X\to Y^{*}$,  then $\phi \oplus_{\rm as}{A}$ belongs to ${\cal L}_{\rm AD}(X\times Y)$
 \end{enumerate}
\end{proposition}
  \noindent{\bf Proof:} (1) and the stability by  multiplication with a scalar is straightforward. (2) follows from the above lemma and (3) is obvious.  
  To show (4)  fix $(q, y)\in X^{*}\times X$, set $r=\Lambda x+p$ and write:
   \begin{eqnarray*}
  L_{\Lambda}^{*} (q,y)&=&\sup\{\langle q, x\rangle +  \langle y, p\rangle- L(x, \Lambda x+p); (x,p)\in X\times X^{*}\}\\
  &=& \sup\{\langle q, x\rangle +  \langle y, r-\Lambda x\rangle- L(x, r); (x,r)\in X\times X^{*}\}\\
  &=& \sup\{\langle q+\Lambda y, x\rangle +  \langle y, r\rangle- L(x, r); (x,r)\in X\times X^{*}\}\\
  &=&  L^{*}(q+\Lambda y, y)=  L(-y, -\Lambda y-q)\\
  &=& L_{\Lambda}(-y, -q).
 \end{eqnarray*}
  For (5) let $r=x-\Lambda^{-1}p$ and   $s=\Lambda x$ and write
 \begin{eqnarray*}
  {}_{\Lambda}L^{*} (q,y)&=&\sup\{\langle q, x\rangle +  \langle y, p\rangle- L(x-\Lambda^{-1}p, \Lambda x); (x,p)\in X\times X^{*}\}\\
  &=& \sup\{\langle q, \Lambda^{-1}s\rangle +  \langle y,s-\Lambda  r\rangle- L(r, s); (r,s)\in X\times X^{*}\}\\
  &=& \sup\{\langle -\Lambda^{-1}q+ y, s\rangle +  \langle \Lambda y, r\rangle- L(r, s); (r,s)\in X\times X^{*}\}\\
  &=&  L^{*}(\Lambda y, -\Lambda^{-1}q+ y)=  L(-y+\Lambda^{-1}q, -\Lambda y)\\
  &=& {}_{\Lambda}L(-y, -q).
 \end{eqnarray*}
 For (6), it is enough to notice that for $(\tilde x, \tilde p) \in (X\times Y)\times (X^{*}\times Y^{*})$, we can write 
 \[
 L\oplus_{A}M (\tilde x, \tilde p)=(L+ M) (\tilde x, \tilde A \tilde x +\tilde p)
 \] where $\tilde A: X\times Y \to X^{*}\times Y^{*}$ is the skew-adjoint operator defined by 
 \[
 \tilde A (\tilde x)=\tilde A ((x,y))= (A^{*}y, -Ax).
 \]
 Assertion  (7) follows from (4) since 
 \[
\phi\oplus_{\rm as}A((x,y), (p,q))=\phi (x, y)+\phi^{*}(-A^{*}y-p, Ax-q)= L_{\tilde A} ((x,y), (p,q))
 \]
where $L((x,y), (p.q)):=\phi (x,y)+\phi^{*}(-p,-q)$ is obviously in ${\cal L}_{\rm AD}(X\times Y)$ and where $\tilde A: X\times Y \to X^{*}\times Y^{*}$ is again the skew-adjoint operator defined by 
 $ \tilde A ((x,y))= (A^{*}y, -Ax).$

\begin{remark}\rm The proof of (4) and (5) above clearly shows that $L_{\Lambda}$ (resp., ${}_{\Lambda}L$) is partially anto-selfdual if and only if $L$ is anti-selfdual on the graph of $\Lambda$. 
\end{remark} 
\begin{remark}\rm
An important use of the above proposition is when  $M_{\lambda}(x, p)=\frac{\|x\|^{2}}{2\lambda^{2}}+\frac{\lambda^{2}\|p\|^{2}}{2}$, then $L_{\lambda}=L\star M_{\lambda}$ is a $\lambda$-regularization of the Lagrangian $L$, which is reminescent of the Yosida theory for operators and for convex functions. This will be most useful in \cite{G3} and \cite{GT2}.  
\end{remark}

\begin{remark} \rm Denote by $\ASDP$ the cone of {\it sub-ASD Lagrangians}: i.e.,  those $L$ in $\LX$ such that 
\begin{equation}
L^{*}(p,x)\geq L(-x, -p)\quad {\rm for\, all}\quad (x,p)\in X\times X^{*}.
\end{equation}
A typical example is a Lagrangian of the form $L(x,p)=\phi (Cx)+\phi^{*}(-p)$ where either $C$ is a surjective operator from $X$ onto itself, or when $C$ has a dense range and $\phi$ is continuous.  It is easy to see that $\ASDP$ also satisfies the following permanence properties:
 \begin{enumerate}
\item  If $L$ is in $ \ASDP$, $M$ is a basic $ASD$-Lagrangian and $\lambda >0$,  then the Lagrangians $L+M$,  $L\star M$ and $\lambda {\bf \cdot} L$  also belong to $\ASDP$.  
\item If $L_{i}\in {\cal L}^{+}_{\rm AD}(X_{i})$ where $X_{i}$ is a reflexive Banach space for each $i\in I$, then $\oplus_{i\in I} L_{i}$   is  in ${\cal L}^{+}_{\rm AD}(\Pi_{i\in I} X_{i})$.
   \item If $L\in \ASDP$ and $\Lambda :X\to X^{*}$ is  skew-adjoint then  $L_{\Lambda}$ is also in $\ASDP$.
 \item If $L\in \ASDP$ and if $\Lambda :X\to X^{*}$ is an invertible  skew-adjoint operator, then ${}_{\Lambda}L$  is also in $\ASDP$. 
 \item If $L \in \ASDP$ and $M\in {\cal L}^{+}_{\rm AD}(Y)$, then for any bounded linear operator $A: X\to Y^{*}$,  the Lagrangian $L\oplus_{A}M$ belongs to ${\cal L}^{+}_{\rm AD}(X\times Y)$
\end{enumerate}
\end{remark}

\section{ASD Lagrangians in  variational  problems with no boundary constraint}
An immediare corollary of Theorem 2.4 in the special case of ASD Lagrangians is the following result which will be used repeatedly in the sequel. 

 \begin{theorem} Let $\Lambda :X\to X^{*}$ be a bounded linear skew-adjoint operator on a reflexive Banach space $X$, and let $L$ be an anti-self dual Lagrangian on the graph of   $\Lambda$. Assume one of the following hypothesis:
\begin{trivlist}
\item {\rm (A)} $\lim\limits_{\|x\|\to \infty}\frac{L(x,\Lambda x)}{\|x\|}=+\infty$, or
 
\item {\rm (B)} The operator $\Lambda$ is invertible and the map $x\to L(x, 0)$ is bounded above on a neighborhood of the origin of $X$. 
\end{trivlist}
Then there exists $\bar x\in X$, such that: 
 \begin{equation}
 \left\{ \begin{array}{lcl}
\label{eqn:existence.2}
  L( \bar x, \Lambda \bar x)&=&\inf\limits_{x\in X} L(x,\Lambda x)=0.\\
 \hfill (-\Lambda \bar x, -\bar x) &\in & \partial L (\bar x,\Lambda \bar x).
\end{array}\right.
 \end{equation}
\end{theorem}
 \noindent{\bf Proof:} It suffices to apply Theorem 2.4 in the case where $R(x)=-x$.  
  In the  case  where $\Lambda$ is also invertible, then we directly apply  Proposition 2.1 to the Lagrangian ${}_{\Lambda}L (x,p)=L(x+\Lambda^{-1}p, \Lambda x)$ which is partially anti-selfdual. \\ 
 
We note that in view of Remark 2.7, it is sufficient to have a Lagrangian $L$ in $\ASDP$ that  is non-negative on the graph of $\Lambda$, that is if $L(x, \Lambda x)\geq 0$ for all $x\in X$.

\subsubsection*{Example 1: A variational formulation for the Lax-Milgram theorem}

Given a bilinear continuous functional $a$ on a Banach space $X$, and assuming that $a$ is coercive: i.e.,  for some $\lambda >0$, we have that $a(v,v)\geq \lambda \|v\|^{2}$ for every $v\in X$. It is well known that if $a$ is symmetric, then for any  $f\in X^{*}$, we can use a variational approach to find $u\in X$, such that for every $v\in X$, we have $
a(u,v) = \langle v,f\rangle$. The procedure amounts to minimize on $H$ the convex functional 
$
\psi (u)=\frac{1}{2}a(u,u)-\langle u,f\rangle.$

The theorem of Lax-Milgram deals with the case when $a$ is not symmetric, for which the above variational argument does not work. Theorem 4.1 however yields the following  variational formulation and proof of the original Lax-Milgram theorem. 

\begin{corollary} Let $a$ be a coercive continuous bilinear form  on $X\times X$. For any $f\in X^{*}$, consider the functional 
\[
I(v)=\psi (v) +\psi^*(-\Lambda v)
\]
where $\psi (v)=\frac{1}{2}a(v,v)-\langle v,f\rangle$, $\psi^{*}$ its Legendre conjugate and where $\Lambda :X\to X^{*}$ is the skew-adjoint operator  defined by 
$
\langle \Lambda v,w\rangle=\frac{1}{2}(a(v,w)-a(w,v)).
$
Then, there exists $u\in X$, such that 
\[
I( u)=\inf_{v\in H}I(v)=0\quad \hbox{\rm
  and 
 \quad $a(u,v) = \langle v,f\rangle$ for every $v\in X$}.
 \]
\end{corollary}

\noindent {\bf Proof:}  
Consider the Lagrangian  $L(x,p)=\psi (x) +\psi^{*}(-p)$ which is clearly anti-self dual. Apply Theorem 4.1 and note  that:
\[
L(u, \Lambda u)=0 \quad \hbox{\rm if and only if \quad $\psi (u) +\psi^{*}(-\Lambda u)=0=-\langle \Lambda u,u\rangle$,} 
\]
 which means that $-\Lambda u\in \partial \psi (u)$. In other words, we have for every $v\in X$
\[
-\frac{1}{2}(a(u,v)-a(v,u))=\frac{1}{2}(a(u,v)+a(v,u))- \langle v,f\rangle
\]
which yields our claim. 

\subsubsection*{Example 2: Inverting variationally a non-selfadjoint matrix}
 An immediate finite dimensional application  of the above corollary is the following variational solution for  the linear equation $Ax=y$ where $A$ is an $n\times n$-matrix and $y\in \R^{n}$. It then suffices to minimize
\[
I(x)=\frac{1}{2}\langle Ax,x\rangle +\frac{1}{2}\langle A^{-1}_{s}(y-A_{a}x),y-A_{a}x\rangle -\langle y,x\rangle.
\]
on $\R^{n}$, where $A_{a}$ is the anti-symmetric part of $A$ and $A^{-1}_{s}$ is the inverse of the symmetric part. If $A$ is coercive, i.e., $
\langle Ax,x\rangle \geq c|x|^{2}$ for all $x\in {\bf R}^{n}$, 
then there is a solution $\bar x \in {\bf R}^{n}$ to the equation obtained as
$
I(\bar x)=\inf_{x\in {\bf R}^{n}}I(x)=0.
$

\subsection{ASD Lagrangians as representations of certain maximal monotone operators}
As noted above, the basic examples of anti-selfdual Lagrangians are of the form 
\begin{equation}
L(x,p)=\varphi (x) +\varphi^{*}(-Bx-p)
\end{equation}
 where $\varphi$ is a convex and lower semi-continuous function on $X$, $\varphi^{*}$ is its Legendre conjugate on $X^{*}$ and where $B:X\to X^{*}$ is skew-symmetric.  This suggests that ASD Lagrangians are natural extensions of  operators of the form $A+\partial \phi$,  where $A$ is positive and $\phi$ is convex. This is an important subclass of maximal monotone operators which can now be resolved variationally. 

Indeed, first consider the cone $\CX$ of all bounded below, proper convex lower semi-continuous functions on $X$, and let $\PX$ be the cone of all positive bounded linear operators from $X$ into $X^{*}$ (i.e., $\langle Ax, x\rangle \geq 0$ for all $x\in X$). Consider also the subclasses 
\[
{\cal C}_{0}(X) =\{\phi \in \CX; \inf_{x\in X} \phi (x) =0\}\quad \hbox{\rm and \quad ${\cal A}_{0}(X) =\{A \in \PX; A^{*}=-A\}$}.
\] 
\begin{proposition} (1) There is a  projection  $\Pi: (\CX, \PX) \to ({\cal C}_{0}(X), {\cal A}_{0}(X))$ such that if $(\phi_{0}, A_{0})$ is the image of 
$(\phi, A)$ by $\Pi$, then a pair $(x, f)\in X\times X^{*}$ satisfies $(A+\partial \phi)(x)=f$ if and only if $(A_{0}+\partial \phi_{0})(x)=f$. \\
(2) For any pair $(\phi, A)\in \CX \times \PX$ there exists a Lagrangian $L_{_{(\phi, A)}}\in \ASD$ such that the equation $(A+\partial \phi)(x)=0$ has a solution ${\bar x}\in X$  if and only if  the functional $I(x)=L_{_{(\phi, A)}}(x,0)$ attains its infimum.
\end{proposition}
\noindent{\bf Proof:} (1) Define the projection as follows: 
 For $(\phi, A)\in (\CX, \PX)$,  decompose $A$ into a symmetric $A^{s}$ and an anti-symmetric part  $A^{a}$, by simply writing $A^{s}=\frac{1}{2}(A +A^{*})$ and $A^{a}=\frac{1}{2}(A -A^{*})$. Let $\phi_{0}$ be the convex functional $\psi +\psi^{*}(0)$, where $\psi (x)=\frac{1}{2} \langle Ax, x\rangle +\varphi (x)$, and define the projection as $\Pi (\phi, A)=(\phi_{0}, A^{a}).$
 
  (2) Associate to each pair $(\phi, A)\in \CX \times \PX$, the anti-selfdual Lagrangian 
\[
L_{(\phi, A)}(x, p)=L_{(\phi_{0}, A_{a})}(x, p)=\phi_{0}(x)+\phi_{0}^{*}(-A^{a}x -p)\quad \hbox{\rm for any $(x,p)\in X\times X^{*}$},
\]
 where $(\phi_{0}, A^{a})$ is the projection of $(\phi, A)$. The fact that the minimum of $I(x)= \phi_{0} (x) +\phi_{0}^{*}(- A^{a}x)$ is equal to $0$ and is attained  at some $\bar x \in X$ means that
\[
\phi_{0} (\bar x) +\phi_{0}^{*}(-A^{a}\bar x)=0=-\langle A^{a}\bar x,\bar x\rangle
\] 
which yields, in view of Legendre-Fenchel duality that $
- A^{a}\bar x \in  \partial \phi_{0} (\bar x)= A^{s}\bar x  +\partial \varphi (\bar x) $, 
hence $\bar x$ satisfies $-Ax\in \partial \varphi (x)$.  

\begin{remark} \rm We note the following relations between classical operations on functions and operators and the operations on ASD Lagrangians. 
\begin{itemize}

\item For $\lambda >0$ and $\phi \in {\cal C}(X)$ , we have $\lambda {\bf \cdot} L_{(\phi, A)}=L_{(\lambda^{2}\phi (\frac{\cdot}{\lambda}), A)}$. 

\item $L_{(\phi_{1}, A_{1})}+L_{(\phi_{2}, A_{2})}=L_{(\phi_{1}+\phi_{2}, A_{1}+A_{2})}$.

\item   $L_{(\phi_{\lambda}, 0)}= L_{(\phi, 0)}\star L_{(\frac{\|x\|^{2}}{2\lambda^{2}}, 0)}$ where  $\phi_{\lambda}$ is the Yosida regularization of $\phi$. 

\item More generally, $L_{(\phi_{1}, 0)}\star L_{(\phi_{2}, 0)}=L_{(\phi_{1}\star\phi_{2}, 0)}$

\end{itemize}

\end{remark}
 In the sequel,  whenever $\phi$ is a functional on $X$ and $f\in X^{*}$,  we shall  denote by $\phi +f$ the functional  defined for $x\in X$ by $\phi (x) +\langle f, x\rangle$. Now we can a variational resolution to the following  nonlinear Lax-Milgram type result.

\begin{corollary}  Assume one of the following conditions on a pair $(\phi, A)\in \CX \times \PX$:
 \begin{trivlist}
 \item {\rm (A)} $\lim\limits_{\|x\|\to \infty}\|x\|^{-1}(\varphi (x) + \frac{1}{2}\langle Ax, x\rangle) =+\infty$, or  
 \item {\rm (B)} The operator $A^{a}=\frac{1}{2}(A -A^{*}):X\to X^{*}$ is onto and $\phi$ is bounded above on the bounded sets of $X$.
\end{trivlist}
Then, there exists for any $f\in X^{*}$, a solution $\bar x \in X$ to the equation $-Ax+f\in \partial \varphi (x)$ that can be obtained as a minimizer of the problem:
\begin{equation}
 \label{min10}
\inf_{x\in X}\left\{\psi (x) +\psi^{*}(- A^{a}x) \right\}=0  
\end{equation}
where $\psi$ is the convex functional $\psi (x)=\frac{1}{2} \langle Ax, x\rangle +\varphi (x) -\langle f, x\rangle$.    \end{corollary}
\noindent{\bf Proof:} Associate to $(\psi, A)\in \CX \times \PX$, the anti-selfdual Lagrangian 
\[
L_{(\psi, A)}(x, p)=L_{(\psi, A^{a})}(x, p)=\psi (x)+\psi^{*}(-A^{a}x -p),\quad \hbox{\rm for $(x,p)\in X\times X^{*}$},
\]
 The fact that the minimum in (\ref{min10}) is attained at some $\bar x \in X$, follows from Theorem 4.1. It means that
$
\psi (\bar x) +\psi^{*}(-A^{a}\bar x)=-\langle A^{a}\bar x,\bar x\rangle$ which yields, in view of Legendre-Fenchel duality that $
- A^{a}\bar x \in  \partial \psi (\bar x)= A^{s}\bar x  +\partial \varphi (\bar x) -f$, 
hence $\bar x$ satisfies $-Ax+f\in \partial \varphi (x)$.  

\begin{remark} All what is needed in the above proposition is that the function $\phi$ be $A$-convex for some operator $A$, meaning that $\psi (x)=\frac{1}{2} \langle Ax, x\rangle +\varphi (x)$ is convex and lower semi-continuous. 
\end{remark}
  
\subsubsection*{Example 3: A variational principle for a non-symmetric Dirichlet problem}
Let ${\bf a}:\Omega \to  {\bf R^{n}}$   be a smooth function on a bounded domain $\Omega$ of  $\bf R^{n}$, and consider the first order linear operator
\[
Av={\bf a}\cdot \nabla v=\Sigma_{i=1}^{n}a_{i}\frac{\partial v}{\partial x_{i}} \]
Assume   that the vector field $\Sigma_{i=1}^{n}a_{i}\frac{\partial v}{\partial x_{i}} $ is actually the restriction of a smooth vector field $\Sigma_{i=1}^{n}{\bar a_{i}}\frac{\partial v}{\partial x_{i}}$ defined on an open neighborhood $X$ of $\bar \Omega$ and that each  ${\bar a_{i}}$ is a $C^{1,1}$ function on $X$. 
Consider the Dirichlet problem:
\begin{equation}
\label{Ex1}
 \left\{ \begin{array}{lcl}
    \hfill Ê\Delta u+\Sigma_{i=1}^{n}a_{i}\frac{\partial u}{\partial x_{i}}&=& |u|^{p-2}u +f  \hbox{\rm \, on \, $\Omega$}
\\
 \hfill  u &=& 0 \quad \quad \quad \quad \hbox{\rm on \quad $\partial \Omega$. }\nonumber \\
       \end{array}  \right.
   \end{equation}
 If $a_{i}=0$, then to find a solution, it is sufficient to minimize the functional
 \[
 \Phi (u)=\frac{1}{2} \int_{\Omega} | \nabla  u |^{2}dx +\frac{1}{p}\int_{\Omega}|u|^{p} dx+\int_{\Omega}fu dx
 \]
and get the solution of $\partial \Phi (u)=0$.\\
 However,  if the non self-adjoint term $a$ is not zero,  we can use the above to get
 \begin{theorem} Assume  ${\rm div }({\bf a})\geq 0$ on $\Omega$, and consider on $H^{1}_{0}(\Omega)$, 
 the functional 
\[
I(u)=\Psi (u) +\Psi^{*}({\bf a}.\nabla u+\frac{1}{2}{\rm div }({\bf a}) \, u )
\]
where
\[
 \Psi (u)=\frac{1}{2} \int_{\Omega} | \nabla  u |^{2}dx +\frac{1}{p}\int_{\Omega}|u|^{p} dx+\int_{\Omega}fu dx +\frac{1}{4}\int_{\Omega}{\rm div }({\bf a}) \, |u|^{2} dx
 \]
 and $\Psi^{*}$ is its Legendre transform. Then, there exists $\bar u\in H^{1}_{0}(\Omega)$ such that:
\[
I(\bar u)=\inf \{I(u);  u\in H^{1}_{0}(\Omega)\}=0, 
\]
and $\bar u$ is a solution of $(\ref{Ex1})$.
 \end{theorem}
\noindent{\bf Proof:}  Indeed, $\Psi$ is clearly convex and lower semi-continuous on $H^{1}_{0}(\Omega)$ while the operator 
$\Lambda u=-{\bf a}.\nabla u-\frac{1}{2}{\rm div }({\bf a}) \, u$ is skew-adjoint, since 
\[
\int_{\Omega} ({\bf a}.\nabla u) u+\frac{1}{2}{\rm div }({\bf a}) \, u^{2} dx =0.
\]
 Again the functional $ I(u)=\Psi (u) +\Psi^{*}({\bf a}.\nabla u+\frac{1}{2}{\rm div }({\bf a}) \, u$ is given by a self-dual Lagrangian $L(u, \Lambda u)$ where $L(u, v)= \Psi (u) +\Psi^{*}(v)$ is defined on $H^{1}_{0}(\Omega)\times H^{-1}(\Omega)$. The existence follows from Theorem  4.1, since  $\Psi$ is clearly coercive. Note that $\bar u$ then satisfies 
 \[
{\bf a}.\nabla \bar u+\frac{1}{2}{\rm div }({\bf a}) \, \bar u =\partial \Psi (\bar u)=-\Delta \bar u+ \bar u^{p-1} +f +\frac{1}{2}{\rm div }({\bf a}) \, \bar u
\]
and therefore $\bar u$ is a solution for $(\ref{Ex1})$.

 \subsubsection*{Example 4: A variational solution for variational inequalities} 
Given again a bilinear continuous functional $a$ on $X\times X$, and $\varphi:X\to {\bf R}$ a convex l.s.c, then solving the corresponding  variational inequality amounts to constructing for any $f\in X^{*}$, 
a point $y\in X$ such that for all $z\in X$,
\begin{equation}
\label{varineq.0}
a(y,y-z)+\varphi (y)-\varphi (z) \leq \langle y-z, f\rangle.
\end{equation}
It is well known that this problem can be rewritten as
\[
f\in Ay+\partial \varphi (y)
\]
where $A$ is the bounded linear operator from $X$ into $X^{*}$ defined by $a(u,v)=\langle Au, v \rangle$. This means that the variational inequality (\ref{varineq.0}) can be rewritten and solved using the variational principle (21). For example, we can solve variationally the following ``obstacle'' problem.

\begin{corollary} Let $a$ be bilinear continuous functional $a$ on a reflexive Banach space $X\times X$ so that $a(v,v)\geq \lambda \|v\|^{2}$, and let $K$ be a convex closed subset of $X$. Then, for any $f\in X^{*}$, there is $\bar x\in K$  such that 
\begin{equation}
\label{varineq2}
a(\bar x,\bar x-z) \leq \langle \bar x-z, f\rangle \quad \hbox{\rm for all $z\in K$}.
\end{equation}
The point $\bar x$  can be obtained as  a minimizer of the following problem:
\[ 
\inf_{x\in X}\left\{\varphi (x) +(\varphi+\psi_{K})^{*}(-\Lambda x) \right\}=0  
\]
where $\varphi (u)=\frac{1}{2}a(u,u)-\langle f, x\rangle$, $\Lambda :X\to X^{*}$ is the skew-adjoint operator  defined by 
\[
\langle \Lambda u,v\rangle=\frac{1}{2}(a(u,v)-a(v,u)).
\]
 and where $\psi_{K} (x)=0$ on $K$ and $+\infty$ elsewhere.
\end{corollary}
\subsection{ASD Lagrangians and anti-Hamiltonian systems}
Recall that an important class of Hamiltonian systems can be written as
\[
(A^{*}y,Ax) \in \partial H(x, y)
\]
where $A:X\to Y^{*}$ is a --normally symmetric-- operator and $H$ is a convex (Hamiltonian) on $X\times Y$.  The next proposition show however that the theory of ASD-Lagrangians is more suited for ``Anti-Hamiltonian'' systems of the form
\[
(-A^{*}y,Ax) \in \partial H(x, y).
\]
 \begin{proposition} Let $\phi$ be any coercive and proper convex lower semi-continuous function  on $X\times Y$ with $(0,0)\in {\rm dom}(\phi)$, and let $A:X\to Y^{*}$ be any bounded linear operator. Assume $B_{1}:X\to X^{*}$ (resp., $B_{2}:Y\to Y^{*}$) are skew-adjoint operators,  then there exists $(\bar x,\bar y) \in X\times Y$ such that
\begin{equation}
(-A^{*}\bar y+B_{1}\bar x,A\bar x+B_{2}\bar y) \in \partial \phi (\bar x, \bar y).
\end{equation}
The solution is obtained as a minimizer on $X\times Y$ of the functional 
 \[
I(x,y)=\phi (x, y)+\phi^{*}(-A^{*}y+B_{1}x, Ax+B_{2}y).
\]
 \end{proposition} 
\noindent {\bf Proof:} It is enough to apply Theorem 4.1 to the ASD Lagrangian
 \[
L((x,y), (p,q))=\phi(x, y)+\phi^{*}(-A^{*}y+B_{1}x-p, Ax+B_{2}y-q).
\]
obtained by shifting to the right the ASD Lagrangian $\phi\oplus_{\rm as} A$ by the skew-adjoint operator $(-B_{1}, -B_{2})$. This yields that $I(x,y)=L((x,y), (0, 0))$ attains its minimum at some $(\bar x,\bar y) \in X\times Y$ and that the minimum is actually $0$. In other words, 
\begin{eqnarray*}
0&=&I(\bar x,\bar y)=\phi(\bar x, \bar y)+\phi^{*}(-A^{*}\bar y+B_{1}\bar x, A\bar x+B_{2}\bar y)\\
&=&\phi (\bar x, \bar y)+\phi^{*}(-A^{*}\bar y+B_{1}\bar x, A\bar x+B_{2}\bar y)-\langle (\bar x, \bar y),  (-A^{*}\bar y+B_{1}\bar x, A\bar x+B_{2}\bar y)\rangle
\end{eqnarray*}
from which the equation follows. 

\begin{corollary} Given positive operators $B_{1}:X\to X^{*}$, $B_{2}:Y\to Y^{*}$ and convex functions $\phi_{1}$ in $\CX$ and  $\phi_{2}$ in $\CY$ having $0$ in their respective domains,  we consider the convex functionals $\psi_{1} (x)=\frac{1}{2} \langle B_{1}x, x\rangle +\varphi_{1} (x)$ and  $\psi_{2} (x)=\frac{1}{2} \langle B_{2}x, x\rangle +\varphi_{2} (x)$. Let $A:X\to Y^{*}$ be a bounded linear operator such that
    \[
 \lim\limits_{\|x\|+\|y\|\to \infty}\frac{\psi_{1}(x)
 +\psi_{2}(y)}
 {\|x\|+\|y\|}=+\infty,
 \]
    Then, for any $(f, g)\in X^{*}\times Y^{*}$ and any $c\in \R$,  there exists a solution $(\bar x, \bar y) \in X\times Y$ to the system of equations 
      \begin{equation}
 \left\{ \begin{array}{lcl}
\label{eqn:existence}
-A^{*}y-B_{1}x+f &\in& \partial \varphi_{1} (x)\\
\hfill c^{2}Ax-B_{2}y+g&\in& \partial \varphi_{2} (y).
\end{array}\right.
 \end{equation}
 It can be obtained as a minimizer of the problem:
 \begin{equation}
 \label{min1}
\inf_{x,y\in X\times Y}\left\{\chi_{1} (x) +\chi_{1}^{*}(-B_{1}^{a}x-A^{*}y) +\chi_{2} (y)+\chi_{2}^{*}(-B_{2}^{a}y+c^{2}Ax)\right\}=0  
\end{equation}
where $B_{1}^{a}$ (resp., $B_{2}^{a}$) are the skew-symmetric parts of $B_{1}$ and $B_{2}$ and where $\chi_{1}(x)=\psi_{1} (x) -\langle f, x\rangle$ and $\chi_{2}(x)=\psi_{2} (x) -\langle g, x\rangle$
  \end{corollary}
\noindent{\bf Proof:} This follows by applying the above proposition to the convex function 
$\phi (x,y)=\chi_{1}(x) +\chi_{2}(y)$ and the skew-symmetric operators $-B_{1}^{a}$ and $-B_{2}^{a}$.  Note that the operator ${\tilde A}: X\times Y\to X^{*}\times Y^{*}$ defined by  ${\tilde A} (x,y)=(A^{*}y, -c^{2}Ax)$ is skew adjoint once we equip $X\times Y$ with the scalar product 
\[
\langle (x,y), (p,q)\rangle= \langle x,p\rangle +c^{-2}\langle y, q\rangle.
\]
We then get  
\begin{equation}
 \left\{ \begin{array}{lcl}
\label{eqn:existence}
-A^{*}y-B^{a}_{1}x+f &\in& \partial \varphi_{1} (x)+ B_{1}^{s}(x)\\
\hfill c^{2}Ax-B^{a}_{2}y+g&\in& \partial \varphi_{2} (y) + B_{2}^{s}(y).
\end{array}\right.
\end{equation}
which gives the result. \\
Another approach consists of associating to the pairs $(\phi_{1}, B_{1})$ and  $(\phi_{2}, B_{1})$ the anti-selfdual Lagrangians 
\[
L(x, p)= \psi_{1} (x) -\langle f, x\rangle +\psi_{1}^{*}(-B_{1}^{a}x +f -p),\quad \hbox{\rm for $(x,p)\in X\times X^{*}$},
\]
and
\[
M(y, q)= \psi_{2} (y)-\langle g, y\rangle+\psi_{2}^{*}(-B_{2}^{a}y +g -q),\quad \hbox{\rm for $(y,q)\in Y\times Y^{*}$},
\] 
Now apply Theorem 4.1  to the twisted-sum Lagrangian $L\oplus_{A}M$.  

\subsubsection*{Example 5: A variational principle for coupled equations}
Let ${\bf b_{1}}:\Omega \to  {\bf R^{n}}$ and  ${\bf b_{2}}:\Omega \to  {\bf R^{n}}$ be two smooth vector fields on a bounded domain $\Omega$ of  $\bf R^{n}$, verifying the conditions in example 3 and let $B_{1}v={\bf b_{1}}\cdot \nabla v$ and  $B_{2}v={\bf b_{2}}\cdot \nabla v$ be the corresponding first order linear operators.   Consider the Dirichlet problem:
\begin{equation}
\label{Ex1.500}
 \left\{ \begin{array}{lcl}
    \hfill Ê\Delta (v+u) +{\bf b_{1}}\cdot \nabla u &=& u^{p-1} +f  \hbox{\rm \, on \, $\Omega$}
\\
  \hfill Ê\Delta (v-c^{2}u) +{\bf b_{2}}\cdot \nabla v&=& v^{q-1} +g  \hbox{\rm \, on \, $\Omega$}\\
 \hfill  u=v &=& 0 \quad \quad \quad \quad  \hbox{\rm on  $\partial \Omega$. }\nonumber \\
       \end{array}  \right.
   \end{equation}
 We can use the above to get
 \begin{theorem} Assume  ${\rm div }({\bf b_{1}})\geq 0$ and ${\rm div }({\bf b_{2}})\geq 0$ on $\Omega$, $1<p,q\leq \frac{n+2}{n-2}$  and consider on $H^{1}_{0}(\Omega)\times H^{1}_{0}(\Omega)$ 
 the functional 
\[
I(u, v)=\Psi (u) +\Psi^{*}({\bf b_{1}}.\nabla u +\frac{1}{2}{\rm div }({\bf b_{1}}) \,  u+\Delta v )+\Phi (v) +\Phi^{*}({\bf b_{2}}.\nabla v +\frac{1}{2}{\rm div }({\bf b_{2}}) \,   v -c^{2}\Delta u )
\]
where
\[
 \Psi (u)=\frac{1}{2} \int_{\Omega} | \nabla  u |^{2}dx +\frac{1}{p}\int_{\Omega}|u|^{p} dx+\int_{\Omega}fu dx +\frac{1}{4}\int_{\Omega}{\rm div }({\bf b_{1}}) \, |u|^{2} dx, 
 \]
 \[
 \Phi (v)=\frac{1}{2} \int_{\Omega} | \nabla  v |^{2}dx +\frac{1}{q}\int_{\Omega}|v|^{q} dx+\int_{\Omega}gv dx +\frac{1}{4}\int_{\Omega}{\rm div }({\bf b_{2}}) \, |v|^{2} dx
 \]
and $\Psi^{*}$ and $\Phi^{*}$ are their Legendre transforms. Then there exists $(\bar u, \bar v)\in H^{1}_{0}(\Omega)\times H^{1}_{0}(\Omega)$ such that:
\[
I(\bar u, \bar v)=\inf \{I(u, v);  (u, v)\in H^{1}_{0}(\Omega)\times H^{1}_{0}(\Omega)\}=0, 
\]
and $(\bar u, \bar v)$ is a solution of $(\ref{Ex1.500})$.
 \end{theorem}
 We can also reduce general minimization problems of functionals of the form $I(x)=\phi (x) +\psi (Ax)$ to the much easier problem of minimizing ASD Lagrangians. Indeed we have
 
 \begin{proposition} Let $\phi$ (resp., $\psi$)  be a convex lower semi-continuous function on a reflexive Banach space $X$ (resp. $ Y^{*})$ and let $A:X\to Y^{*}$ be a bounded linear operator. To minimize the functional $I(x)=\phi (x) +\psi (Ax)$ on $X$, we consider on $X\times Y$ the functional 
 \[
 I(x,y)=\phi (x) +\psi^{*} (y)+\phi^{*}(-A^{*}y)+\psi (Ax). 
 \]
Assuming $\lim\limits_{\|x\|+\|y\|\to \infty} I(x,y)=+\infty$, then the infimum of $I$ is zero and is attained at a point $(\bar x, \bar y)$ which determines the extremals of the min-max problem:
\[
\sup\{-\psi^{*}(y)-\phi^{*}(-A^{*}y)  =\inf\{\phi (x)+\psi (Ax); x\in X\}.
\] 
 \end{proposition}
They also satisfy the system: 
     \begin{equation}
 \left\{ \begin{array}{lcl}
\label{eqn:existence}
-A^{*}y &\in& \partial \varphi (x)\\
\hfill Ax&\in& \partial \psi^{*} (y).
\end{array}\right.
 \end{equation}
\noindent{\bf Proof:} It is sufficient to note that $I(x,y)=L((x,y), (0,0)$ where $L$ is an anti-self dual Lagrangian defined on $X\times Y$ by:
\[
L((x,y), (p,q))=\phi (x) +\psi^{*} (y)+\phi^{*}(-A^{*}y-p)+\psi (Ax-q). 
\]
 By considering more general twisted sum Lagrangians,  we obtain the following application

\begin{theorem} Let $X$ and $Y$ be two reflexive Banach spaces and let $A: X\to Y^{*}$ be any bounded linear operator.  Assume  $L \in \ASD$ and $M\in {\cal L}_{\rm AD}(Y)$ are such that
 \[
 \lim\limits_{\|x\|+\|y\|\to \infty}\frac{L(x,A^{*}y)+M(y, -Ax)}{\|x\|+\|y\|}=+\infty,
 \]
 Then there exists $(\bar x, \bar y)\in X\times Y$, such that:
 \begin{equation}
L(\bar x,A^{*}\bar y)+M(\bar y, -A\bar x)=\inf\limits_{(x,y)\in X\times Y)}L(x,A^{*}y)+M(y, -Ax)=0.
\end{equation}
Moreover, we have
    \begin{equation}
 \left\{ \begin{array}{lcl}
\label{eqn:existence}
\hfill L(\bar x, A^{*}\bar y)+\langle \bar x, A^{*}\bar y\rangle&=&0\\ 
M(\bar y, -A\bar x)+\langle \bar y, -Ax\rangle&=&0\\
\hfill (-A^{*} \bar y, -\bar x) &\in&  \partial L (\bar x, A^{*} \bar x)\\
\hfill (A\bar x, -\bar y) &\in&  \partial M (\bar y,-A \bar x)
 \end{array}\right.
 \end{equation}
   \end{theorem}
\noindent{\bf Proof:} It is sufficient to apply Theorem 4.1 to the ASD Lagrangian  $L\oplus_{A}M$.  
  \section{ASD Lagrangians associated to boundary value problems}
  For problems involving boundaries, we may start with an ASD Lagrangian $L$, but the operator $\Lambda$ may be skew-adjoint modulo a term involving the boundary. Assuming we can represent this term by a pair of operators 
$(b_{1}, b_{2})$ from  $X$ into a Hilbert space $H_{1}\times H_{2}$  which correspond to an adequate splitting of the boundary, then we may try to recover anti-selfduality by adding a correcting term via a boundary Lagrangian $\ell$.  In this section, we look into frameworks where Lagrangians of the form 
 \[
M (x,p)=L(x, \Lambda x+p) +\ell (b_{1}(x), b_{2}(x))
 \]
can be made anti-selfdual.
  
  \subsection{Anti-selfduality  involving boundary Lagrangians}
\begin{definition}\rm  (1) A {\it boundary operator} will be any surjective  continuous linear map $(b_{1}, b_{2}): X\to H_{1}\times H_{2}$ from $X$ onto the product of Hilbert spaces $H_{1}\times H_{2}$. \\  
(2) An operator $\Lambda: X\to X^{*}$ is said to be {\it skew-symmetric modulo the boundary} operator $(b_{1}, b_{2})$, if for every $x,y \in X$, 
\begin{equation}
\langle \Lambda x , y\rangle_{_{(X,X^{*})}} = -\langle \Lambda y , x\rangle_{_{(X,X^{*})}} + \langle b_{2}(x), b_{2}(y) \rangle_{_{H_{2}}} -\langle b_{1}(x), b_{1}(y)\rangle _{_{H_{1}}}
\end{equation}
That is, $\Lambda^{*}=-\Lambda +b_{2}^{*}b_{2} -b_{1}^{*}b_{1}$ which means that the operator $\Lambda - \frac{1}{2}(b_{2}^{*}b_{2} -b_{1}^{*}b_{1})$ is skew symmetric.  We shall then say that we have a skew symmetric triplet $(\Lambda, b_{1}, b_{2})$. 
\end{definition}
   We also consider a {\it Boundary Lagrangian} $\ell :H_{1}\times H_{2} \to \R \cup \{+\infty\}$ which is also proper convex and lower semi-continuous, and its Legendre transform in both variable,
 \[
    \ell^*( h_{1},h_{2})= \sup \{ \braket{k_{1}}{h_{1}} +  \braket{k_{2}}{h_{2}} - \ell(k_{1},k_{2});\, k_{1}\in H_{1},k_{2} \in H_{2}  \}
\]
\begin{definition} We say that $\ell$ is a {\it self-dual boundary Lagrangian} if
\begin{equation}
\ell^*(-h_{1},h_{2})=\ell (h_{1},h_{2}) \quad \hbox{\rm for all $(h_{1},h_{2}) \in H_{1}\times H_{2}$}.
\end{equation} 
\end{definition}
It is easy to see that such a boundary Lagrangian will always satisfy the inequality
\begin{equation}
\label{little.positivity}
\hbox{$\ell (r,s)\geq \frac{1}{2}(\|s\|^2-\|r\|^2)$ for all $(r,s)\in H_1\times H_2$.}
\end{equation}
The basic example of a self dual boundary Lagrangian is given by a function $\ell$ on $H_{1}\times H_{2}$, of the form $\ell (r,s) =\psi_{1}(r) +\psi_{2}(s)$, with  $\psi_{1}^{*}(r)=\psi_{1}(-r)$ and $ \psi_{2}^{*}(s)=\psi_{2}(s)$. Here the choices for $\psi_{1}$ and $\psi_{2}$ are rather limited and the typical sample is:
\[
\psi_{1} (r)= \frac{1}{2}\|r\|^{2} -2\langle a, r\rangle +\|a\|^{2}, \quad \hbox{\rm and \quad $\psi_{2}(s)=\frac{1}{2}\|s\|^{2}$.}
\]
where  $a$ is given in $H_{1}$. \\
Boundary operators allow us to build new ASD Lagrangians. We shall present several ways to do so, which correspond to various conditions that $\Lambda$ and the Lagrangian $L$ may or may not satisfy in applications. 

\begin{proposition} Let $\ell$ be a self dual boundary Lagrangian on the Hilbertian product  $H_{1}\times H_{2}$, and let $(\Lambda, b_{1}, b_{2}): X\to X^{*}\times H_{1}\times H_{2}$ be a skew symmetric triplet where $X$ is a reflexive Banach space. Suppose $H$ is a  linear subspace  of $X^{*}$ containing ${\rm Range}(\Lambda)$ such that $X_{0}={\rm Ker}(b_1,b_2)$ is dense in $X$ for the $\sigma (X, H)$-topology, and consider  $L$ to be a  Lagrangian on $X$ such that for each $p\in X^{*}$, the map  $x\to L(x, p)$ is continuous  for the $\sigma (X, H)$-topology. 

\begin{enumerate}
\item If $L$ is anti-self dual on the graph of $\Lambda$,  then the Lagrangian 
 \[
 M(x, p)=L(x, \Lambda x+p) +\ell (b_{1}(x), b_{2}(x))
\]
  is partially anti-self dual.  
\item If $L$ is anti-selfdual on the elements of $X\times H$, then $M$ 
 is also anti-self dual on the elements of $X\times H$.   
   \end{enumerate}
      \end{proposition}
      
 \noindent{\bf Proof:}   \noindent{\bf Proof:} Fix $(q, y) \in X^{*}\times X$ and calculate
 \begin{eqnarray*}
 M^{*} (q,y)&=&\sup\{\langle q, x\rangle_{X} +  \langle y, p\rangle_{X}- M(x, p); (x,p)\in X\times X^{*}\}\\
 &=&\sup\left\{ \langle q, x\rangle_{X} +  \langle y, p\rangle_{X} -L(x, \Lambda x +p) -\ell  (b_{1}(x), b_{2}(x)); (x,p)\in X\times X^{*}\right\} 
 \end{eqnarray*}
 Setting $r=\Lambda x+p$, we obtain 
 \begin{eqnarray*}
 M^{*} (q,y)&=& \sup \left\{\langle x, q\rangle + \langle y, r-\Lambda x\rangle  -L(x, r)-\ell (b_{1}(x), b_{2}(x)); (x,r)\in X\times X^{*}\right\} \\
 &=& \sup \{\langle x, q\rangle   + \langle b_{1}(y), b_{1}(x)\rangle -\langle b_{2}(y), b_{2}(x)\rangle  + \langle \Lambda y,  x\rangle +\langle y, r\rangle\\
 && \quad \quad \quad  -L(x, r) -\ell (b_{1}(x), b_{2}(x)); (x,r)\in X\times X^{*}\}\\
 &=& \sup \{\langle x, q+\Lambda y\rangle   +\langle y, r\rangle-L(x, r)\\
 && \quad \quad   + \langle b_{1}(y), b_{1}(x)\rangle -\langle b_{2}(y), b_{2}(x)\rangle   -\ell (b_{1}(x), b_{2}(x)); (x,r)\in X\times X^{*}\}\\
  &=&\sup \{\langle x, q+\Lambda y\rangle   +\langle y, r\rangle-L(x, r)
   + \langle b_{1}(y), b_{1}(x+x_{0})\rangle -\langle b_{2}(y), b_{2}(x+x_{0})\rangle\\
   && \quad \quad  \quad   -\ell (b_{1}(x+x_{0}), b_{2}(x)+x_{0}); (x_{0}, x,r)\in X_{0}\times X\times X^{*}\}\\
   &=&\sup \{\langle w-x_{0}, q+\Lambda y\rangle   +\langle y, r\rangle-L(w-x_{0}, r)
   + \langle b_{1}(y), b_{1}(w)\rangle -\langle b_{2}(y), b_{2}(w)\rangle\\
     && \quad \quad  \quad   -\ell (b_{1}(w), b_{2}(w); (x_{0}, w,r)\in X_{0}\times X\times X^{*}\}
 \end{eqnarray*}
 Now suppose  $(q, y) \in H\times X$, and use the fact that $X_0$ is $\sigma (X, H)$ dense in $X$, that $Range \Lambda \subset H$ and the continuity of $x\to L(x, p)$ in that topology to obtain
   \begin{eqnarray*}
 M^{*} (q,y)
 &=&\sup \{\langle z, q+\Lambda y\rangle   +\langle y, r\rangle-L(z, r)
   + \langle b_{1}(y), b_{1}(w)\rangle -\langle b_{2}(y), b_{2}(w)\rangle\\
   && \quad \quad  \quad   -\ell (b_{1}(w), b_{2}(w); (z, w,r)\in X\times X\times X^{*}\}\\
   &=&\sup \{\langle z, q +\Lambda y\rangle  +\langle y, r\rangle  - L(z, r); (z,r)\in X\times X^{*}\}\\
 &&+\sup\{ \langle b_{1}(y), b_{1}(w)\rangle -\langle b_{2}(y), b_{2}(w)\rangle -\ell (b_{1}(w), b_{2}(w));  w\in X\}\\
 &=&\sup \{\langle z, q +\Lambda y\rangle  +\langle y, r\rangle  - L(z, r); (z,r)\in X\times X^{*}\}\\
 &&+\sup\{ \langle b_{1}(y), a\rangle -\langle b_{2}(y), b\rangle -\ell (a, b);  (a,b)\in H_{1}\times H_{2}\}\\
 &=&L^{*}(q+\Lambda y, y)+\ell^{*} (b_{1}(y), -b_{2}(y)) \\
 &=&L(-y, -q-\Lambda y)+\ell (-b_{1}(y), -b_{2}(y))\\
 &=&M(-y, -q)
 \end{eqnarray*}
  Here is another situation that occurs in certain applications.
\begin{definition} \rm Say that $(b_{1}, b_{2})$ is a  {\it regular boundary operator} if there is a projection $\Pi: X \to X_0:=Ker (b_{1}, b_{2})$ so that the bounded linear map $(\Pi, b_{1}, b_{2}):X\to Ker (b_{1}, b_{2})\oplus H_{1}\oplus H_{2}$ is an isomorphism.
\end{definition}
Denote by $K: X \to X_{0}$ the projection in such a way that the bounded linear map $(K, b_{1}, b_{2}):X\to X_{0}\oplus H_{1}\oplus H_{2}$ is an isomorphism. We can identify $X^{*}$ with the space $X_{0}^{*} \oplus H_{1}\oplus H_{2}$ in such a way that the duality between $X$ and $X^{*}$ is given by:
  \[
 \langle x, p \rangle = \langle x, (p_{0}, p_{1}, p_{2}) \rangle=\langle x, K^{*}p_{0}\rangle_{X, X^*} +\langle b_{1}(x), p_{1}\rangle_{H_1} +\langle b_{2}(x), p_{2}\rangle_{H_2}.
  \]

\begin{proposition}  Let $\ell$ be a self dual boundary Lagrangian on the Hilbertian product $H_{1}\times H_{2}$, and let $(\Lambda, b_{1}, b_{2}): X\to X^{*}\times H_{1}\times H_{2}$ be a regular skew symmetric triplet where $X$ is a reflexive Banach space. Consider  $L$ to be a  Lagrangian on $X$ such that for each $x\in X$, the map $p\to L(x,p)$ is continuous on $X^*$. 
\begin{enumerate}
\item If $L$ is a Lagrangian on $X\times X^{*}$ that is anti-self dual on the graph of $\Lambda$,  then the Lagrangian 
 \[
 M(x, p)=L(x, \Lambda x+K^{*}p_{0}) +\ell (b_{1}(x)+p_{1}, b_{2}(x)-p_{2})
\]
  is partially anti-self dual. Here $X^{*}$ is identified with $X_{0}^{*} \oplus H_{1}\oplus H_{2}$. 
\item If $L$ is anti-self dual  on $X$,  then  $M$  is anti-self dual on $X\times X_{0}^{*}$.   
   \end{enumerate}
     \end{proposition}
     
  \noindent{\bf Proof:} Fix $(q, y) \in X^{*}\times X$, with $q=(q_{0}, 0,0)$  and calculate
\begin{eqnarray*}
M^{*} (q,y)&=&\sup\{\langle q, x\rangle_{X} +  \langle y, p\rangle_{X}- M(x, p); (x,p)\in X\times X^{*}\}\\
&=&\sup\left\{\langle x, K^{*}q_{0}\rangle   +\langle b_{1}(y), p_{1}\rangle +\langle b_{2}(y), p_{2}\rangle+\langle y, K^*p_{0}\rangle \right. \\
 &&\left.-L(x, \Lambda x +K^{*}p_{0}) -\ell (b_{1}(x)+p_{1}, b_{2}(x)-p_{2}); x\in X, p_{0}\in X_{0}^{*}, p_{1}\in H_{1}, p_{2}\in H_{2}\right\} 
\end{eqnarray*}
where $(p_{0}, p_{1}, p_{2})$ represent $p\in X^{*}$. 

Since the operator $A=\Lambda - \frac{1}{2}(b_{2}^{*}b_{2} -b_{1}^{*}b_{1})$ is skew-adjoint on $X$, we can apply the results of the last section to get that $(A+\epsilon I)$ is onto for each $\epsilon>0$. In other words, $A$ has dense range in $X^*$, which yields that ${\rm Range} (\Lambda) +K^*(X^*_0)={\rm Range} (\Lambda) + {\rm Ker}(b_1,b_2)^{\perp}\supset {\rm Range}(\Lambda - \frac{1}{2}(b_{2}^{*}b_{2} -b_{1}^{*}b_{1})={\rm Range} (A)$ is also dense in $X^*$. 

Setting $r=\Lambda x+K^{*}p_{0}$, $f_{1}=b_{1}(x)+p_{1}$ and $f_{2}=b_{2}(x)-p_{2}$, we obtain that
\begin{eqnarray*}
M^{*} (q,y)&=& \sup\{\langle x, K^{*}q_{0}\rangle  + \langle y, r-\Lambda x\rangle +\langle b_{1}(y), f_{1}-b_{1}(x)\rangle +\langle b_{2}(y), b_{2}(x)-f_{2}\rangle. \\
&&\quad \quad \left.-L(x, r) -\ell (f_{1}, f_{2}); x\in X, r\in {\rm Range} (\Lambda) +K^*(X^*_0), f_{1}\in H_{1}, f_{2}\in H_{2}\right\} \\
&=& \sup\{\langle x, K^{*}q_{0}\rangle  + \langle y, r-\Lambda x\rangle +\langle b_{1}(y), f_{1}-b_{1}(x)\rangle +\langle b_{2}(y), b_{2}(x)-f_{2}\rangle. \\
&&\quad \quad \left.-L(x, r) -\ell (f_{1}, f_{2}); x\in X, r\in X^{*}, f_{1}\in H_{1}, f_{2}\in H_{2}\right\} \\
&=&  \sup\{\langle x, K^{*}q_{0}\rangle  +\langle b_{1}(y), f_{1}\rangle -\langle b_{2}(y), f_{2}\rangle
  +\langle y, r\rangle +   \langle \Lambda y,  x\rangle   \\
&&\quad \quad \left.-L(x, r) -\ell (f_{1}, f_{2}); x\in X, r\in X^{*}, f_{1}\in H_{1}, f_{2}\in H_{2}\right\}\\
&=&L^{*}(K^{*}q_{0}+\Lambda y, y)+\ell^{*} (b_{1}(y), -b_{2}(y)) \\
&=&L(-y, -K^{*}q_{0}-\Lambda y)+\ell (-b_{1}(y), -b_{2}(y))\\
&=&M(-y, -q)
\end{eqnarray*}
since $(q_{0}, 0,0)$ represents $q$ in $X_{0}^{*}\times H_{1}\times H_{2}$.
   
 In the case where $\Lambda$ is essentially onto (modulo the boundary) we have yet another useful setting.

 \begin{definition} \rm Say that a skew symmetric triplet $(\Lambda, b_{1}, b_{2}): X\to X^{*}\times H_{1}\times H_{2}$ is a {\it nice boundary operator} if the map $(\Lambda, b_{1}): X \to {\rm Range} (\Lambda) \oplus H_{1}$ is an isomorphism.
\end{definition}
In this case, we   identify $X^{*}$ with the space $X_{0} \oplus H_{1}$ where $X_{0}=X/Ker (\Lambda))\sim{\rm Range} (\Lambda)$ in such a way that the duality between $X$ and $X^{*}$ is given by:
  \[
 \langle x, p \rangle = \langle x, (p_{0}, p_{1}) \rangle=\langle \Lambda x, p_{0}\rangle +\langle b_{1}(x), p_{1}\rangle.
  \]

\begin{proposition}  Let $\ell$ be a self dual boundary Lagrangian on the Hilbertian product $H_{1}\times H_{2}$, and let $(\Lambda, b_{1}, b_{2}): X\to X^{*}\times H_{1}\times H_{2}$ be a  skew symmetric triplet on a reflexive Banach space $X$ that is a nice boundary operator. Then  
\begin{enumerate}
\item If $L$ is a Lagrangian on $X\times X^{*}$ that is anti-self dual on the graph of $\Lambda$,  then the Lagrangian 
 \[
 N(x, p)=L(x+p_{0}, \Lambda x) +\ell (b_{1}(x)+p_{1}, b_{2}(x))
\]
  is partially anti-self dual on $X$.  
  \item If $L$ is anti-self dual  on $X$,  then  $N$  is anti-self dual on the elements of $X\times (X_{0}\oplus \{0\})$   
   \end{enumerate}
     \end{proposition}
 \noindent{\bf Proof:} Indeed, fix $((q_{0},0), y) \in (X_{0}\times H_{1})\times X$  and calculate
\begin{eqnarray*}
N^{*} (q,y)&=&\sup\{\langle \Lambda x, q_{0}\rangle  + \langle \Lambda y, p_{0}\rangle +\langle b_{1}(y), p_{1}\rangle\\
 &&\left.-L(x+p_{0}, \Lambda x) +\ell (b_{1}(x)+p_{1}, b_{2}(x)); x\in X, p_{0}\in X_{0}, p_{1}\in H_{1}\right\} 
\end{eqnarray*}
Setting $r= x+p_{0}$, $f_{1}=b_{1}(x)+p_{1}$, we obtain that
\begin{eqnarray*}
N^{*} (q,y)&=& \sup\{\langle \Lambda x, q_{0}\rangle  
  +  \langle \Lambda y, r- x\rangle +\langle b_{1}(y), f_{1}-b_{1}(x)\rangle \\
&&\left.-L(r, \Lambda x) -\ell (f_{1}, b_{2}(x)); x\in X, r\in X, f_{1}\in H_{1}\right\} \\
&=&\sup\{\langle \Lambda x, q_{0}\rangle   
  +  \langle \Lambda y, r\rangle +  \langle y, \Lambda x\rangle +\langle b_{1}(y), f_{1}\rangle -\langle b_{2}(y), b_{2}(x)\rangle\\
&&\left.-L(r, \Lambda x) -\ell (f_{1}, b_{2}(x)); x\in X, r\in X, f_{1}\in H_{1}\right\}
\end{eqnarray*}
Since $X$ can be identified with $X_{0}\oplus H_{2}$ via the correspondence $x\to (\Lambda x, b_{2}(x))$, we obtain:
\begin{eqnarray*}
N^{*} (q,y)&=&\sup\{\langle s, q_{0}+y\rangle   
  +  \langle \Lambda y, r\rangle   +\langle b_{1}(y), f_{1}\rangle -\langle b_{2}(y), f_{2}\rangle\\
&&\left.-L(r, s) -\ell (f_{1}, f_{2}); s\in X_{0}, r\in X, f_{1}\in H_{1}, f_{2}\in H_{2}\right\}\\
&=&\sup\{\langle s, q_{0}+y\rangle   
  +  \langle \Lambda y, r\rangle -L(r, s); s\in X_{0}, r\in X\}  \\
&&+\sup\{\langle b_{1}(y), f_{1}\rangle -\langle b_{2}(y), f_{2}\rangle -\ell (f_{1}, f_{2});  f_{1}\in H_{1}, f_{2}\in H_{2}\}\\
&=&L^{*}(\Lambda y, q_{0}+ y )+\ell^{*} (b_{1}(y), -b_{2}(y)) \\
&=&L(-q_{0}- y, -\Lambda y)+\ell (-b_{1}(y), -b_{2}(y))\\
&=&N(-y, -q).
\end{eqnarray*} 

\subsection{Variational properties of ASD Lagrangians with boundary terms}
One can now deduce the following
 \begin{theorem} Let $(\Lambda, b_{1}, b_{2}): X\to X^{*}\times H_{1}\times H_{2}$ be a  skew symmetric triplet,  $\ell$  a self dual boundary Lagrangian on $H_{1}\times H_{2}$ and let $L:X\times X^{*}\to {\bf R}\cup\{+\infty\}$ be anti-self dual on the graph of $\Lambda$. Assume one of the following hypothesis:
 \begin{trivlist}
 \item  {\rm (A)} The boundary operator $(b_{1},b_{2})$ is regular and  for every $x\in X$, the map $p\to L(x,p)$  
 is bounded on the bounded sets of $X^*$.
 
  \item {\rm (B)} The triplet $(\Lambda, b_{1},b_{2})$ is a nice boundary operator and  the map $x \to L (x, 0)$ is bounded  on the bounded sets of $X$.  
 \end{trivlist}
 Then, there exists $\bar x \in X$ such that:
  \begin{equation}
\label{eqn:zero.2}
 L(\bar x, \Lambda \bar x) +  \ell (b_{1}\bar x, b_{2}\bar x)=\inf_{x\in X} \left\{L(x, \Lambda x) +  \ell (b_{1}x, b_{2}x)\right\}= 0.
 \end{equation}
 Moreover, we have
  \begin{equation}
 \left\{ \begin{array}{lcl}
\label{eqn:application.1}
  L(\bar x, \Lambda \bar x) + \langle \bar x, \Lambda \bar x\rangle
 &=&0.\\
  \hfill (-\Lambda \bar x, -\bar x) &\in & \partial L (\bar x,\Lambda \bar x)\\
\hfill \ell( b_{1}(\bar x), b_{2}(\bar x))&=& \frac{1}{2}(\|b_{2}\bar x\|^{2} -\|b_{1}\bar x\|^{2}).
 \end{array}\right.
 \end{equation}
  In particular, for any $a\in H_{1}$ there exists  $\bar x \in X$ such that  $b_{1}(\bar x)=a$ and satisfying (\ref {eqn:application.1}). It is obtained as a minimizer on $X$ of the functional 
\[
I(x)=L(x, \Lambda x)+\frac{1}{2}\|b_{1}(x)\|^{2} -2\langle a, b_{1}(x)\rangle +\|a\|^{2} +\frac{1}{2}\|b_{2}(x)\|^{2}.
\] 
 \end{theorem}
 \noindent{\bf Proof:} Under case (A), we use proposition 3.2 to get that the Lagrangian 
 \[
 M(x, p)=L(x, \Lambda x+p_0) +\ell (b_{1}(x)+p_1, b_{2}(x)-p_2)
\]
  is partially anti-self dual.  
    In case (B), we use Proposition 5.3 to conclude that the Lagrangian 
 \[
 N(x, p)=L(x+p_{0}, \Lambda x) +\ell (b_{1}(x)+p_{1}, b_{2}(x))
\]
 is partially anti-self dual on $X$.   \\
 In both cases,  the hypothesis  implies that $M(0,p)$ (resp., $N(0,p)$) is bounded above on the bounded sets of $X^{*}$. Theorem 4.1 then applies to yield $\bar x \in X$ such that (\ref{eqn:zero.2}) is satisfied.\\
  To establish (\ref{eqn:application.1}), write 
 \begin{eqnarray*}
  L(x, \Lambda x) +  \ell (b_{1}x, b_{2} x)&=& L( x, \Lambda x)+ \langle  x, \Lambda x\rangle- \langle x, \Lambda  x\rangle + \ell (b_{1} x, b_{2} x)\\
  &=&L(x, \Lambda x)+ \langle x, \Lambda x\rangle- \frac{1}{2}(\|b_{2}x\|^{2} -\|b_{1}x\|^{2})  + \ell (b_{1} x, b_{2} x). 
\end{eqnarray*}  
Since $L(x,p)\geq -\langle x, p\rangle$ and $\ell (r,s)\geq \frac{1}{2}(\|s\|^2-\|r\|^2)$, we immediately obtain (\ref{eqn:application.1}).\\
  In particular, for any $a\in H_{1}$,  consider the  boundary Lagrangian,  
\[
\ell (r,s)= \frac{1}{2}\|r\|^{2} -2\langle a, r\rangle +\|a\|^{2} +\frac{1}{2}\|s\|^{2}. 
\]
which is clearly self-dual. We then get 
\begin{eqnarray*}
  L(  x, \Lambda   x) +  \ell (b_{1} x, b_{2}  x)&=&L(  x, \Lambda   x)+ \langle   x, \Lambda   x\rangle- \frac{1}{2}(\|b_{2}x\|^{2} -\|b_{1}x\|^{2})  + \ell (b_{1}x, b_{2}x)\\&=&
  L(x, \Lambda x) + \langle  x, \Lambda  x\rangle +\|b_{1}(x)-a\|^{2}.
\end{eqnarray*}  
In other words, $\bar x \in X$  is a solution of  $\inf\limits_{x\in X}\left\{L(x, \Lambda x) + \langle  x, \Lambda  x\rangle +\|b_{1}(x)-a\|^{2}\right\}=0$, and since $L(x,p)\geq -\langle x, p\rangle$, we obtain:
  \begin{equation}
 \left\{ \begin{array}{lcl}
\label{eqn:application}
  L(\bar x, \Lambda \bar x) + \langle \bar x, \Lambda \bar x\rangle
 &=&0.\\
\hfill  b_{1}(\bar x)&=&a. 
 \end{array}\right.
 \end{equation}
\subsection{Variational principle for operators which are positive modulo a boundary}
   Consider again  $(b_{1}, b_{2}): X\to H_{1}\times H_{2}$ to be a regular boundary operator. 
\begin{definition}   Say that  $A:X\to X^{*}$ is {\it positive modulo the boundary operator $(b_{1}, b_{2})$} if the operator $A-\frac{1}{2}(b_{2}^{*}b_{2} -b_{1}^{*}b_{1})$  is positive. 
 \end{definition}

\begin{corollary} Let $A:X\to X^{*}$ be positive modulo a boundary operator $(b_{1}, b_{2})$
 and set $\Lambda=A^{a}+\frac{1}{2}(b_{2}^{*}b_{2} -b_{1}^{*}b_{1})$. Let $\phi$ be a convex function in $\CX$ with $0$ in its domain and such that one of the following conditions holds:
\begin{trivlist}
\item {\rm (A)} The boundary operator $(b_{1}, b_{2})$ is regular and 
\[
 \lim\limits_{\|x\|\to \infty}\|x\|^{-1} \big\{\phi (x)+ \frac{1}{2}\langle Ax, x\rangle -\frac{1}{4}(  \|b_{2}x\|^{2} -\|b_{1}x\|^{2})\big\}=+\infty.
\]
  \item {\rm (B)} The triplet $(\Lambda, b_{1},b_{2})$ is a nice boundary operator and $\phi$ is bounded on the bounded sets of $X$. 
 \end{trivlist}
Then for any $a \in H_{1}$ and any $f\in X^{*}$,  the equation 
\begin{equation}
 \left\{ \begin{array}{lcl}
\label{eqn:laxmilgram3}
 \hfill  -Ax&\in &\partial \varphi (x)+f\\
\hfill  b_{1}(x)&=&a \\
\end{array}\right.
\end{equation}
has a solution $\bar x \in X$ that is a minimizer of the problem:
\[ 
I(x)=\psi (x) +\psi^{*}(- A^{a}x-\frac{1}{2}b_{2}^{*}b_{2}x +\frac{1}{2}b_{1}^{*}b_{1} x)+\frac{1}{2}(\|b_{1}(x)\|^{2} +\|b_{2}(x)\|^{2} )-2\langle a, b_{1}(x)\rangle +\|a\|^{2}   
\]
where $\psi (x)=\phi (x)+ \frac{1}{2}\langle Ax, x\rangle -\frac{1}{4}(  \|b_{2}x\|^{2} -\|b_{1}x\|^{2})+ \langle f, x\rangle$.   \end{corollary}
\noindent{\bf Proof:} Let $B= A-\frac{1}{2}(b_{2}^{*}b_{2} -b_{1}^{*}b_{1})$ and decompose it into its symmetric $B^{s}=A^{s}-\frac{1}{2}(b_{2}^{*}b_{2} -b_{1}^{*}b_{1})$ and its anti-symmetric part $B^{a}=A^{a}$, by simply writing $B^{s}=\frac{1}{2}(B +B^{*})$ and $B^{a}=\frac{1}{2}(B -B^{*})$. We can then write $A=B^{s}+\Lambda$ where $\Lambda :=A^{a}+\frac{1}{2}(b_{2}^{*}b_{2} -b_{1}^{*}b_{1})=B^{a}+\frac{1}{2}(b_{2}^{*}b_{2} -b_{1}^{*}b_{1})$ which means that $\Lambda$ is skew-symmetric modulo the boundary operator $(b_{1}, b_{2})$. 
 For $f\in X^{*}$,  consider  the convex functional 
 \[
 \psi (x)=\varphi (x)+\frac{1}{2} \langle B^{s}x, x\rangle  +\langle f, x\rangle=\phi (x)+ \frac{1}{2}\langle Ax, x\rangle -\frac{1}{4}(  \|b_{2}x\|^{2} -\|b_{1}x\|^{2})+\langle f, x\rangle. 
 \]
 The proposed minimization problems amounts to applying Theorem 5.5 to 
the anti-self dual Lagrangian  $L(x,p)=\psi (x) +\psi^{*}(-p)$, the operator $\Lambda$ and the boundary Lagrangian $\ell (r,s)= \frac{1}{2}\|r\|^{2} -2\langle a, r\rangle +\|a\|^{2} +\frac{1}{2}\|s\|^{2}$. 
Note that 
\[
 I(x)=\psi (x) +\psi^{*}(- \Lambda x)+\langle x, \Lambda x \rangle +\|b_{1}(x)-a\|^{2}.
 \]
The fact that the minimum is attained at some $\bar x$ and is equal to $0$,  implies that $b_{1}(\bar x)=a$ and that $
\psi (\bar x) +\psi^{*}(-\Lambda \bar x)=-\langle \Lambda\bar x,\bar x\rangle$ 
which means that 
\[
-A^{a}({\bar x})-\frac{1}{2}(b_{2}^{*}b_{2} -b_{1}^{*}b_{1})({\bar x})=-\Lambda \bar x \in  \partial \psi (\bar x)=  \partial \varphi (\bar x) + A^{s}({\bar x})-\frac{1}{2}(b_{2}^{*}b_{2} -b_{1}^{*}b_{1})({\bar x})+f
\]
 and therefore $-A{\bar x}\in \partial \varphi ({\bar x})+f$.  
 
  \begin{remark} \rm  Again the above applies to functions $\phi$ that are {\it A-convex modulo a boundary $(b_1, b_2)$} meaning those functions $\phi$ such that there exists an operator $A$ such that  $\psi (x) =\phi (x)+ \frac{1}{2}\langle Ax, x\rangle -\frac{1}{4}(  \|b_{2}x\|^{2} -\|b_{1}x\|^{2})$ is convex and lower semi-continuous. 
\end{remark}  
    
\subsubsection*{Example 6: A variational principle for non-linear transport equations} 
As in example 3,  Let ${\bf a}:\Omega \to  {\bf R^{n}}$ and $a_{0}:\Omega \to  {\bf R}$ be two smooth functions on a bounded domain $\Omega$ of  $\bf R^{n}$, and consider the first order linear operator
\[
Av={\bf a}\cdot \nabla v=\Sigma_{i=1}^{n}a_{i}\frac{\partial v}{\partial x_{i}}\quad \hbox {and \quad $\Lambda v= {\bf a}\cdot \nabla v +a_{0}v.$}
\]
As in \cite{Ba}, we shall assume throughout that the vector field $\Sigma_{i=1}^{n}a_{i}\frac{\partial v}{\partial x_{i}} $ is actually the restriction of a smooth vector field $\Sigma_{i=1}^{n}{\bar a_{i}}\frac{\partial v}{\partial x_{i}}$ defined on an open neighborhood $X$ of $\bar \Omega$ and that each  ${\bar a_{i}}$ is a $C^{1,1}$ function on $X$. We  also assume that the boundary of $\Omega$ is piecewise $C^{1}$, in such a way that the outer normal ${\bf n}$ is defined almost everywhere on $\partial \Omega$. In this case, if we denote by 
\[
\Sigma_{-}=\{x\in \partial \Omega; \, {\bf n(}x)\cdot {\bf a}(x) <0\}\, \hbox{\rm 
  and  $\Sigma_{+}=\partial \Omega \setminus \Sigma_{-}=\{x\in \partial \Omega; \, {\bf n}(x)\cdot {\bf a}(x) \geq 0\}$},
\]
then a trace $u_{|_{\Sigma_{-}}}$ makes sense in $L^{2}_{\rm loc}(\Sigma_{-})$ as soon as $u\in L^{2}(\Omega)$ and $\Lambda u\in L^{2}(\Omega)$. \\
Let now $\beta :{\bf R}\to {\bf R}$ be a continuous nondecreasing function so that its antiderivative $j$ is convex, and let $f\in L^{2}(\Omega)$. We are interested in finding variationally solutions for the nonlinear transport equation:
\begin{equation}
\label{transport.eq.1}
 \left\{ \begin{array}{lcl}
    \hfill -{\bf a}\cdot \nabla u -a_{0}u&=& \beta (u) +f  \hbox{\rm \, on \, $\Omega$}
\\
 \hfill  u(x) &=& u_{0} \quad \quad \quad \quad \hbox{\rm on \quad $\Sigma_{-}$. }\nonumber \\
       \end{array}  \right.
   \end{equation}
 First, we identify the appropriate underlying space. Consider the space
\[
H^{1}(\Omega)=\{u\in L^{2}(\Omega);\, A u\in L^{2}(\Omega)\}.
\]
equipped with the norm $\|u\|_{H^{1}}=\|u\|_{2}+\|A u\|_{2}$. 
As noticed in \cite{Ba},  that a function $u$ belongs to $H^{1}(\Omega)$ does not necessarily  guarantee that its trace $u_{|_{\Sigma_{-}}}$ is in the space
\[
L^{2}_{A}(\Sigma_{-})=\left\{ u\in L^{2}_{\rm loc}(\Sigma_{-}); 
\int_{{\Sigma_{-}}}|u(x)|^{2}|{\bf n}(x)\cdot {\bf a}(x)|  d\sigma <+\infty \right\}.
\]
However, if $u\in H^{1}(\Omega)$ and  $u_{\Sigma_{-}} \in L^{2}_{A}(\Sigma_{-})$, then necessarily $u_{\Sigma_{+}} \in L^{2}_{A}(\Sigma_{+})$. The appropriate space for our setting is therefore
\[
H^{1}_{A}(\Omega)=\{u\in H^{1}(\Omega); u_{|_{\Sigma_{-}}} \in L^{2}_{A}(\Sigma_{-})\}.
\]
equipped with the norm $\|u\|_{H_{A}^{1}}=\|u\|_{2}+\|Au\|_{2} +\|u_{|_{\Sigma_{-}}}\|_{_{L^{2}_{A}(\Sigma_{-})}}$. 

To define appropriate boundary spaces, we follow \cite{DL} and consider for each open subset $\Gamma$ of $\partial \Omega$, the space 
\[
H^{^{1/2}}_{_{00}}(\Gamma)=\left\{ v\in L^{2}_{A}(\Gamma); \hbox{ $\exists w\in H^{1}(\Omega)$, $w=0$ on $\partial \Omega \setminus \Gamma$, and $w=v$ on $\Gamma$}\right\}
\]
A {\it trace theorem} (\cite{LM}, Vol III. p. 307) or \cite{BC}) yields that the restriction mapping $u \to u_{\Gamma}$ is a continuous surjective map from $V=\{v\in H^{1}(\Omega); v_{|_{\partial \Omega \setminus \Gamma}}=0\}$ onto $H^{^{1/2}}_{_{00}}(\Gamma)$. It follows that  there is a continuous surjection from $H^{1}_{A}(\Omega)$ onto $H^{1}_{0}(\Omega)\oplus H^{^{1/2}}_{_{00}}(\Sigma_{-})\oplus H^{^{1/2}}_{_{00}}(\Sigma_{+})$ via the map
\[
T: H^{1}_{A}(\Omega) \to  H^{1}_{0}(\Omega)\oplus H^{^{1/2}}_{_{00}}(\Sigma_{-})\oplus H^{^{1/2}}_{_{00}}(\Sigma_{+}), 
\]
given by $Tu =(Ku, u_{|_{}{\Sigma_{-}}}, u_{|_{\Sigma_{+}}})$, where $K: H^{1}_{A}(\Omega) \to  H^{1}_{0}(\Omega)$ is the operator that associates to $u\in H^{1}_{A}(\Omega)$ the unique function 
$w\in H^{1}_{0}(\Omega)$ such that $\Delta w=\Delta u$ and $w=0$ on $\partial \Omega$.

If now $a_{0}(x)- \frac{1}{2}{\rm div} a (x) \geq 0$ on $\Omega$, 
then $\Lambda$ is positive modulo the boundary operators $u\to (u_{|_{\Sigma_{-}}}, u_{}{|_{\Sigma_{+}}} ) \in L^{2}(\Sigma_{-})\times L^{2}(\Sigma_{+})$, since 
\[
\int_{\Omega}u\Lambda u dx =\int_{\Omega} (({\bf a \cdot } \nabla u) u +a_{0}|u|^{2})dx=
\int_{\Omega}(a_{0}-\frac{1}{2}{\rm div \, {\bf a} })|u|^{2}dx+ \int_{\partial \Omega}|u|^{2}{\bf n}\cdot {\bf a}  d\sigma.
\]
and the operator 
\[
\Lambda_{1}(u) := {\bf a}\cdot \nabla u+\frac{1}{2}{\rm div} ({\bf a})u=\Lambda (u) -(a_{0}-\frac{1}{2}{\rm div \, {\bf a} })u
\]
is therefore skew-adjoint modulo that boundary since then
\begin{equation}
\label{antisymmetry}
\int_{\Omega} v\Lambda_{1}u \, dx = -\int_{\Omega}u \Lambda_{1}v \, dx + \int_{\partial \Omega}uv\, {\bf n}\cdot {\bf a}  d\sigma.
\end{equation}
We can now state:
 \begin{theorem}  Assume the coercivity condition
$a_{0}(x)- \frac{1}{2}{\rm div} a (x) \geq \alpha >0$ on $\Omega$. For any $f\in L^{2}(\Omega)$ and $u_{0}\in L^{2}_{A}(\Sigma_{-})$, consider  the following functional on the space $H^{1}_{A}(\Omega)$
\begin{eqnarray}
\label{functional}
{I}(u)&=&\psi (u)+\psi^{*}(-{\bf a}\cdot \nabla u-\frac{1}{2}({\rm div}{\bf a}) u)+\frac{1}{2}\int_{\Sigma_{+}}|u(x)|^{2}{\bf n}(x)\cdot {\bf a}(x)\, d\sigma
-\frac{1}{2}\int_{\Sigma_{-}}|u(x)|^{2}{\bf n}(x)\cdot {\bf a}(x)\,  d\sigma \nonumber \\ &&+2\int_{\Sigma_{-}}u(x)u_{0}(x){\bf n}(x)\cdot {\bf a}(x)\,  d\sigma
-\int_{\Sigma_{-}}|u_{0}(x)|^{2}{\bf n}(x)\cdot {\bf a}(x)\,  d\sigma.\nonumber
\end{eqnarray}
 and where $\psi$ is the convex functional  on $L^{2}(\Omega)$ defined by
\[
\psi (u)=\int_{\Omega}\left\{j(u(x))+f(x)u(x)+\frac{1}{2}(a_{0}-\frac{1}{2}{\rm div \, a }) u^{2})\right\}dx
\]  
and where $\psi^{*}$ is its Legendre conjugate.
 
Then there exists a solution $\bar u$ for (\ref{transport.eq.1}) that is obtained as a minimizer of the problem:
\[
I(\bar u)=\inf\{ I (u);u\in H^{1}_{A}(\Omega)\}=0.
\]
\end{theorem}
  \noindent{\bf Proof:} The only problem remaining is the fact that the convex functional $\psi$ defined by:
\[
\psi (u)=\int_{\Omega}\left\{j(u(x))+f(x)u(x)+\frac{1}{2}(a_{0}-\frac{1}{2}{\rm div \, a }) u^{2})\right\}dx
\]  
is not necessarily coercive  on $H^{1}_{A}(\Omega)$, so we consider instead for each $\epsilon >0$, the functional
\[
\varphi_{\epsilon}(u)=\psi (u)+\frac{\epsilon}{2}\int_{\Omega}|\nabla u|^{2}dx 
\] 
which obviously is.  Assuming without loss that $u_{0}=0$ and 
setting
\begin{equation}
\label{functional}
I_{\epsilon}(u)=\varphi_{\epsilon} (u)+\varphi_{\epsilon}^{*}(-\Lambda_{1}u) +\frac{1}{2}\int_{\Sigma_{+}}|u(x)|^{2}{\bf n}\cdot {\bf a}\, d\sigma-\frac{1}{2}\int_{\Sigma_{-}}|u(x)|^{2}{\bf n}\cdot {\bf a}\, d\sigma.
\end{equation}
The above lemma now applies and we get $u_{\epsilon}\in 
H^{1}_{A}(\Omega)$ such that
\[
\inf_{u\in H^{1}_{A}(\Omega)}I_{\epsilon} (u)=I_{\epsilon}(u_{\epsilon})=0,
\]
This means that $u_{\epsilon}$ belongs to $ {\rm Dom} (\partial \varphi_{\epsilon})$ and satisfies
$
-\Lambda_{1}u_{\epsilon} \in \partial \varphi_{\epsilon} (u_{\epsilon})$, 
 which implies
 \[
-\Lambda_{1}u_{\epsilon}=\beta (u_{\epsilon}) +f + (a_{0}-\frac{1}{2}{\rm div \, a })u_{\epsilon}-\epsilon \Delta u_{\epsilon}. 
\]
 In other words, we have for each $\epsilon >0$, 
\begin{equation}
\label{transport.eq}
 \left\{ \begin{array}{lcl}
    \hfill \epsilon \Delta u_{\epsilon}-a_{0}  u_{\epsilon}-\Sigma_{i=1}^{n}a_{i}\frac{\partial u_{\epsilon}}{\partial x_{i}} &=& \beta (u_{\epsilon}) +f  \hbox{\rm \, on \, $\Omega$}
\\
 \hfill  u_{\epsilon} &=& 0 \quad \quad \quad \quad \hbox{\rm on \quad $\Sigma_{-}$, }\nonumber \\
\hfill \frac{\partial u_{\epsilon}}{\partial n} &=&0 \quad \quad \quad \quad \hbox{\rm on \quad $\partial \Omega \setminus \Sigma_{-}$. }\nonumber
       \end{array}  \right.
   \end{equation}
It is now standard to show that,  as $\epsilon \to 0$, $u_{\epsilon}$ converges in $L^{2}(\Omega)$ to a solution $u$ of (\ref{transport.eq.1}). For details, see Bardos \cite{Ba}. 

\subsection{ASD Lagrangians on intermediate Hilbert spaces}
As one can see in the previous example, it is more desirable to have coercivity on the space $L^{2}(\Omega)$ and therefore we need to ``extend''  anti-selfduality from the Banach space $H^{1}_{A}(\Omega)$ to the ambient Hilbert space $L^{2}(\Omega)$.\\
 This situation is common in applications to partial differential equations, where an ambient Hilbert space $H$ is usually present in such a way that $X$ is a dense subset of $H$, and the identity injection $i: X\to H$ is continuous. The scalar product and the norm of $H$ are denoted by $(u,v)$ and $|\, |$ respectively. By duality, the adjoint $i^{*}: H^{*}\to X^{*}$ is also one-to-one with dense range. One often identifies $H$ with its dual $H^{*}$, in such a way that we have a representation of the form: $X\subset H\equiv H^{*}\subset X^{*}$. In this representation, we have $\langle h,x\rangle=(h,x)$ whenever $h\in H$ and $x\in X$.  The pair $(X, H)$ is sometimes called an {\it evolution pair}. We shall need the following notion.
 
 \begin{definition} Let  $(\Lambda, b_{1}, b_{2}): X\to X^{*}\times H_{1}\times H_{2}$ be a  skew symmetric triplet on a reflexive Banach space,  and let $H$ be a Hilbert space so that $(X, H)$ is  an evolution pair. We say that $(X, H, \Lambda)$  is a {\it maximal evolution triple} if $X_{0}={\rm Ker}(b_{1},b_{2})$ is dense in $H$,  $\Lambda$ maps $X$ into $H$  and if 
 \[
 X=\{x\in H; \,  \sup\{\langle x, \Lambda y\rangle_{H}-\frac{1}{2}(\|b_1(x)\|_{H_1}^2+\|b_2(x)\|_{H_2}^2); \, y\in X, \|y\|_{H}\leq 1\}<+\infty\}.
 \]
 \end{definition}

 \begin{proposition} Let $\ell$ be a self dual boundary Lagrangian on the Hilbertian product  $H_{1}\times H_{2}$, let $(\Lambda, b_{1}, b_{2}): X\to X^{*}\times H_{1}\times H_{2}$ be a regular skew symmetric triplet  on a reflexive Banach space $X$, and let $H$ be a  Hilbert space such that $(X, H, \Lambda)$  is a {\it maximal evolution triple}.
 
  If $L$ is anti-self dual on  $H$ such that for each $p\in H$, the map $x\to L(x, p)$ is continuous on $H$, then the Lagrangian 
 \[
M (x,p)= \left\{ \begin{array}{lcl}
  \hfill L(x, \Lambda x+p) +\ell (b_{1}(x), b_{2}(x))\, &{\rm if}& {x\in X} \\
  +\infty &\,& {\rm otherwise} \\
\end{array}\right.
\]
 is also anti-self dual on  $H\times H$.   
      \end{proposition}
       \noindent{\bf Proof:} 
       Fix $(q, y) \in H\times X$ and calculate
 \begin{eqnarray*}
 M^{*} (q,y)&=&\sup\{\langle q, x\rangle_{H} +  \langle y, p\rangle_{H}- M(x, p); (x,p)\in H\times H\}\\
 &=&\sup\left\{ \langle q, x\rangle_{X} +  \langle y, p\rangle_{X} -L(x, \Lambda x +p) -\ell  (b_{1}(x), b_{2}(x)); (x,p)\in X\times H\right\} 
 \end{eqnarray*}
 Setting $r=\Lambda x+p$, we obtain since $y\in X$, that
 \begin{eqnarray*}
 M^{*} (q,y)&=& \sup \left\{\langle x, q\rangle + \langle y, r-\Lambda x\rangle  -L(x, r)-\ell (b_{1}(x), b_{2}(x)); (x,r)\in X\times H\right\} \\
 &=& \sup \{\langle x, q\rangle   + \langle b_{1}(y), b_{1}(x)\rangle -\langle b_{2}(y), b_{2}(x)\rangle  + \langle \Lambda y,  x\rangle +\langle y, r\rangle\\
 && \quad \quad \quad  -L(x, r) -\ell (b_{1}(x), b_{2}(x)); (x,r)\in X\times H\}\\
 &=& \sup \{\langle x, q+\Lambda y\rangle   +\langle y, r\rangle-L(x, r)\\
 && \quad \quad   + \langle b_{1}(y), b_{1}(x)\rangle -\langle b_{2}(y), b_{2}(x)\rangle   -\ell (b_{1}(x), b_{2}(x)); (x,r)\in X\times H\}\\
  &=&\sup \{\langle x, q+\Lambda y\rangle   +\langle y, r\rangle-L(x, r)
   + \langle b_{1}(y), b_{1}(x+x_{0})\rangle -\langle b_{2}(y), b_{2}(x+x_{0})\rangle\\
   && \quad \quad  \quad   -\ell (b_{1}(x+x_{0}), b_{2}(x)+x_{0}); (x_{0}, x,r)\in X_{0}\times X\times H\}\\
   &=&\sup \{\langle w-x_{0}, q+\Lambda y\rangle   +\langle y, r\rangle-L(w-x_{0}, r)
   + \langle b_{1}(y), b_{1}(w)\rangle -\langle b_{2}(y), b_{2}(w)\rangle\\
     && \quad \quad  \quad   -\ell (b_{1}(w), b_{2}(w); (x_{0}, w,r)\in X_0\times X\times H\}
 \end{eqnarray*}
 Now  use the fact that $X_0$ is  dense in $H$,  and the continuity of $x\to L(x, p)$ on $H$, to obtain
   \begin{eqnarray*}
 M^{*} (q,y)
 &=&\sup \{\langle z, q+\Lambda y\rangle   +\langle y, r\rangle-L(z, r)
   + \langle b_{1}(y), b_{1}(w)\rangle -\langle b_{2}(y), b_{2}(w)\rangle\\
   && \quad \quad  \quad   -\ell (b_{1}(w), b_{2}(w); (z, w,r)\in H\times X\times H\}\\
   &=&\sup \{\langle z, q +\Lambda y\rangle  +\langle y, r\rangle  - L(z, r); (z,r)\in H\times H\}\\
 &&+\sup\{ \langle b_{1}(y), b_{1}(w)\rangle -\langle b_{2}(y), b_{2}(w)\rangle -\ell (b_{1}(w), b_{2}(w));  w\in X\}\\
 &=&\sup \{\langle z, q +\Lambda y\rangle  +\langle y, r\rangle  - L(z, r); (z,r)\in H\times H\}\\
 &&+\sup\{ \langle b_{1}(y), a\rangle -\langle b_{2}(y), b\rangle -\ell (a, b);  (a,b)\in H_{1}\times H_{2}\}\\
 &=&L^{*}(q+\Lambda y, y)+\ell^{*} (b_{1}(y), -b_{2}(y)) \\
 &=&L(-y, -q-\Lambda y)+\ell (-b_{1}(y), -b_{2}(y))\\
 &=&M(-y, -q)
 \end{eqnarray*}

      If now $(q, y) \in H\times (H \setminus X)$, then 
 \begin{eqnarray*} 
M^{*} (q,y)&\geq &\sup \left\{\langle x, q\rangle + \langle y, r-\Lambda x\rangle  -L(x, r)-\ell (b_{1}(x), b_{2}(x)); (x,r)\in X\times H\right\} \\
&\geq &\sup \left\{-\|x\|_{H}\|q\|_{H}  + \langle y, \Lambda x\rangle  - L(x, 0)-\ell (b_{1}(x), b_{2}(x)); x\in X\right\}\\
&\geq &\sup \left\{-\|q\|_{H}  + \langle y, \Lambda x\rangle  -C-\frac{1}{2}(\|b_1(x)\|_{H_1}^2+\|b_2(x)\|_{H_2}^2); x\in X, \|x\|_{H}\leq 1\right\}\\
&=&+\infty, 
 \end{eqnarray*}
since otherwise 
$y\in X$ which is a contradiction. 
 
\begin{remark} \rm In practice, the Hilbert space is usually given and $X$ is usually obtained from the domain of some unbounded operator on $H$. The full scope of this setting is developed in \cite{GT2}. For now, we give the following illustrative example
\end{remark}

\subsubsection*{Example 7: More general transport equations}
 Consider the following general transport equation 
\begin{equation}
\label{transport.eq.gen}
 \left\{ \begin{array}{lcl}
    \hfill -\Lambda u&=& |u|^{p-2}u + Bu + f  \hbox{\rm \, on \, $\Omega$}
\\
   u(x) &=& u_{0} (x)\quad \quad \quad \quad \quad \quad \hbox{\rm on \, $\Sigma_{-}$. }\nonumber \\
       \end{array}  \right.
   \end{equation}
where $B: L^{2}(\Omega) \to L^{2}(\Omega)$ is a positive bounded linear operator, $f\in L^{2}(\Omega)$ and $u_{0}\in L^{2}_{A}(\Sigma_{-})$. 

We again decompose $B$ into a symmetric and an anti-symmetric part, $B_{s}$ and $B_{a}$, by writing $B_{s}=\frac{1}{2}(B +B^{*})$ and $B_{a}=\frac{1}{2}(B -B^{*})$, and we consider the convex functional defined on $L^{2}$ by:
\[
\psi (u)=\frac{1}{2} \int_{\Omega}(\frac{1}{p}|u(x)|^{p}+u(x)(Bu)(x) +f(x)u(x)+\frac{1}{2}(a_{0}-\frac{1}{2}{\rm div \, {\bf a} }) u^{2})dx
\]  
and its conjugate $\psi^{*}$. Let again $\Lambda_{1}$ be the operator
\[
\Lambda_{1}(u) = \Sigma_{i=1}^{n}{a_{i}}\frac{\partial u}{\partial x_{i}}+\frac{1}{2}{\rm div} ({\bf a})u=\Lambda (u) -(a_{0}-\frac{1}{2}{\rm div \, {\bf a} })u.
\]
 The functional on $L^{2}(\Omega)$, is now defined as 
\begin{eqnarray}
\label{functional}
{\tilde I}(u)&=&\psi (u)+\psi^{*}(-\Lambda_{1}u - B_{a}u)\nonumber \\ &&+\frac{1}{2}\int_{\Sigma_{+}}|u(x)|^{2}{\bf n}(x)\cdot {\bf a}(x)\, d\sigma
-\frac{1}{2}\int_{\Sigma_{-}}|u(x)|^{2}{\bf n}(x)\cdot {\bf a}(x)\,  d\sigma \nonumber \\ &&+2\int_{\Sigma_{-}}u(x)u_{0}(x){\bf n}(x)\cdot {\bf a}(x)\,  d\sigma
-\int_{\Sigma_{-}}|u_{0}(x)|^{2}{\bf n}(x)\cdot {\bf a}(x)\,  d\sigma.\nonumber
\end{eqnarray}
if $u\in H^{1}_{A}(\Omega)$ and $+\infty$ elsewhere. 
On can then verify the following.
\begin{theorem} If $1<p\leq2$, then there exists $\bar u\in H^{1}_{A}(\Omega)$ such that
\begin{equation}
\label{zero.bis}
{\tilde I}(\bar u)=\inf\{ {\tilde I} (u);u\in H^{1}_{A}(\Omega)\}=0.
\end{equation}
and $\bar u$ solves equation (\ref{transport.eq.gen}). 
 \end{theorem}
 
 \subsection{ASD Lagrangians for coupled equations with prescribed boundaries}
 
 Assume  $L \in \ASD$ and $M\in {\cal L}_{\rm AD}(Y)$ where  $X$ and $Y$ are two reflexive Banach spaces and let $A: X\to Y^{*}$ be any bounded linear operator. Let $(\Lambda, b_{1}, b_{2}): X\to X^{*}\times H_{1}\times H_{2}$ (resp. $(\Gamma, c_{1}, c_{2}): Y\to Y^{*}\times K_{1}\times K_{2}$) be  skew symmetric triplets, and let $\ell$ (resp., $m$ be a self dual boundary Lagrangian on $H_{1}\times H_{2}$ (resp., $K_{1}\times K_{2}$), in such a way that the Lagrangians 
 \[
 L_{\Lambda}(x,p)=L(x, \Lambda x+p) +\ell (b_{1}(x), b_{2}(x))\quad \hbox{\rm and \quad $M_{\Gamma}(y,q)=L(y, \Gamma y+q) +\ell (c_{1}(y), b_{2}(y))$}
 \]
   are ASD and therefore the Lagrangian 
   \[
   L_{\Lambda}\oplus_{_{A}}M_{\Gamma} ((x,y), (p,q))=L( x,A^{*} y+\Lambda x+p)+M( y, -A x+\Gamma y+q) +\ell (b_{1} x, b_{2} x)+m (c_{1} y, c_{2} y)
   \]
 is also anti-selfdual. Consider the functional $I(x,y)=L_{\Lambda}\oplus_{_{A}}M_{\Gamma} ((x,y), (0,0))$, that is
 \[
 I(x,y):=L( x,A^{*} y+\Lambda x)+M( y, -A x+\Gamma y) +\ell (b_{1} x, b_{2} x)+m (c_{1} y, c_{2} y).
 \]
 We can now state
     \begin{theorem} Assume that $
 \lim\limits_{\|x\|+\|y\|\to \infty}\frac{I(x,y)}{\|x\|+\|y\|}=+\infty.$
  Then there exists $(\bar x, \bar y)\in X\times Y$ such that:
  \begin{equation}
\label{eqn:zero.200}
I(\bar x, \bar y)=\inf_{(x, y)\in X\times Y}I(x,y) = 0.
 \end{equation}
  In particular, for any $a\in H_{1}$ and $b\in K_{1}$, there exists  $(\bar x, \bar y)\in X\times Y$ such that:
  \begin{equation}
 \left\{ \begin{array}{lcl}
\label{eqn:application.100}
 L(\bar x, A^{*}\bar y+\Lambda \bar x)+\langle \bar x, A^{*}\bar y+\Lambda \bar x\rangle&=&0\\
 M(\bar y, -A\bar x+\Gamma \bar y)+\langle \bar y, -A\bar x+\Gamma \bar y\rangle&=&0\\ 
\hfill (-A^{*} \bar y-\Lambda \bar x, -\bar x) &\in&  \partial L (\bar x, A^{*} \bar x+\Lambda \bar x)\\
\hfill (A\bar x-\Gamma \bar y, -\bar y) &\in&  \partial M (\bar y,-A \bar x+\Gamma \bar y)\\
\hfill  b_{1}(\bar x)&=&a \\
\hfill c_{1}(\bar y)&=&b
 \end{array}\right.
 \end{equation}
It is obtained as a minimizer on $X\times Y$ of the functional 
\begin{eqnarray*}
I(x,y)&=&L(x, A^{*}y+\Lambda x)+\frac{1}{2}\|b_{1}(x)\|^{2} -2\langle a, b_{1}(x)\rangle +\|a\|^{2} +\frac{1}{2}\|b_{2}(x)\|^{2}\\
&&+M(y, -Ax+\Gamma y)+\frac{1}{2}\|c_{1}(y)\|^{2} -2\langle b, c_{1}(y)\rangle +\|b\|^{2} +\frac{1}{2}\|c_{2}(y)\|^{2}. 
\end{eqnarray*}
 \end{theorem}
\noindent{\bf Proof:} Note that we can rewrite
\begin{eqnarray*}
I(x,y)&=&L(x, A^{*}y+\Lambda x) + \langle  x, A^{*}y+\Lambda  x\rangle\\
&& + M(y, -Ax+\Gamma y)  +\langle  y, -A x+\Gamma  y\rangle\\
 && +\|b_{1}(x)-a\|^{2}+\|c_{1}(x)-b\|^{2},
\end{eqnarray*}
in such a way that if $I(\bar x,\bar y)=0$, then the fact that the sum of each two consecutive terms constituting $I$ above is non-negative, prove our claim (\ref{eqn:application.100}).
\begin{corollary}
Let $B_{1}:X\to X^{*}$ (resp., $B_{2}:Y\to Y^{*}$) be positive operators modulo a regular boundary $(b_{1}, b_{2}):X\to H_{1}\times H_{2}$ (resp., $(c_{1}, c_{2}):Y\to K_{1}\times K_{2}$),  
  let $\phi_{1}$ (resp $\phi_{2}$) be a convex function in $\CX$ (resp. in $\CY$) and consider the convex functions 
\[
\psi_{1} (x)=\frac{1}{2} \langle B_{1}x, x\rangle +\varphi_{1} (x) -\frac{1}{4}(  \|b_{2}x\|^{2} -\|b_{1}x\|^{2})\]
 \[
 \psi_{2} (x)=\frac{1}{2} \langle B_{2}x, x\rangle +\varphi_{2} (x)-\frac{1}{4}(  \|c_{2}x\|^{2} -\|c_{1}x\|^{2})
 \]
 Let $A:X\to Y^{*}$ be a bounded linear operator such that
$ \lim\limits_{\|x\|+\|y\|\to \infty}\frac{\psi_{1}(x)
  +\psi_{2}(y)
 }{\|x\|+\|y\|}=+\infty.$\\
 Then, for any $(a,b) \in H_{1}\times K_{1}$, any $(f,g)\in X^{*}\times Y^{*}$ and any $\alpha\in R$,  there exists a solution $(\bar x, \bar y) \in X\times Y$ to the system of equations 
      \begin{equation}
 \left\{ \begin{array}{lcl}
\label{eqn:existence}
-A^{*}y-B_{1}x+f &\in& \partial \varphi_{1} (x)\\
\hfill \alpha^{2}Ax-B_{2}y+g&\in& \partial \varphi_{2} (y)\\
\hfill b_{1}(\bar x)&=&a\\
\hfill b_{2}(\bar y)&=&b
\end{array}\right.
 \end{equation}
 It is obtained as a minimizer  
 on $X\times Y$ of the functional:
\begin{eqnarray*}
\label{min100}
I(x,y)&=&\chi_{1} (x) +\chi_{1}^{*}(-B_{1}^{a}x-\frac{1}{2}b_{2}^{*}b_{2}x +\frac{1}{2}b_{1}^{*}b_{1} x-A^{*}y)+\frac{1}{2}(\|b_{1}x\|^{2} +\|b_{2}x\|^{2} )-2\langle a, b_{1}(x)\rangle +\|a\|^{2}   \\
&& +\chi_{2} (y)+\chi_{2}^{*}(-B_{2}^{a}y-\frac{1}{2}c_{2}^{*}c_{2}y +\frac{1}{2}c_{1}^{*}c_{1} y +\alpha^{2}Ax)+\frac{1}{2}(\|c_{1}y\|^{2} +\|c_{2}y\|^{2} )-2\langle b, c_{1}y\rangle +\|b\|^{2}   
\end{eqnarray*}
where $\chi_{1} (x)=\psi_{1} (x)- \langle f, x\rangle$ and $\chi_{2} (y)=\psi_{2} (y)- \langle g, y\rangle$.   
\end{corollary}
 \noindent{\bf Proof:} Associate  the following anti-selfdual Lagrangians on $X\times X^{*}$  and $Y\times Y^{*}$ respectively,
 \[
L(x, p)= \chi_{1} (x) +\chi_{1}^{*}(-B_{1}^{a}x-\frac{1}{2}b_{2}^{*}b_{2}x +\frac{1}{2}b_{1}^{*}b_{1} x  -p)+\frac{1}{2}(\|b_{1}x\|^{2} +\|b_{2}x\|^{2} )-2\langle a, b_{1}(x)\rangle +\|a\|^{2}  
\] 
\[
M(y, q)= \chi_{2} (y)+\psi_{2}^{*}(-B_{2}^{a}y-\frac{1}{2}c_{2}^{*}c_{2}y +\frac{1}{2}c_{1}^{*}c_{1} y  -q)+\frac{1}{2}(\|c_{1}y\|^{2} +\|c_{2}y\|^{2} )-2\langle b, c_{1}y\rangle +\|b\|^{2}  
\] 
Now apply the preceeding corollary to these two ASD Lagrangians, to the operator $\alpha^{2}A: X\to Y^{*}$ and to the product $X\times Y$ equipped with the scalar product $
\langle (x,y), (p,q)\rangle= \langle x,p\rangle +\alpha^{-2}\langle y, q\rangle.$ We get  
\begin{equation}
 \left\{ \begin{array}{lcl}
\label{eqn:existence}
-A^{*}y-B^{a}_{1}x+f &\in& \partial \varphi_{1} (x)+ B_{1}^{s}(x)\\
\hfill \alpha^{2}Ax-B^{a}_{2}y+g&\in& \partial \varphi_{2} (y) + B_{2}^{s}(y)\\
\hfill b_{1}(\bar x)&=&a\\
\hfill b_{2}(\bar y)&=&b.
\end{array}\right.
\end{equation}
which gives the result. 

\subsubsection*{Example 8: A variational principle for a coupled system with prescribed boundary conditions}
Let ${\bf a}:\Omega \to  {\bf R^{n}}$ and  ${\bf b}:\Omega \to  {\bf R^{n}}$ be two smooth vector fields on a bounded domain $\Omega$ of  $\bf R^{n}$, verifying the conditions in example 6 and consider their corresponding first order linear operator $B_{1}u={\bf a}\cdot \nabla u $ and  $B_{2}v={\bf b}\cdot \nabla v$.  Let 
\[
\Sigma_{-}^{1}=\{x\in \partial \Omega; {\bf a\cdot n}(x)<0\}\quad {\rm and}\quad  \Sigma_{-}^{2}=\{x\in \partial \Omega; {\bf b\cdot n}(x)<0\}.
\]
For  $u_{0}\in L^{2}_{B_{1}}(\Sigma^{1}_{-})$ and  $v_{0}\in L^{2}_{B_{2}}(\Sigma^{2}_{-})$, 
 consider the Dirichlet problem:
\begin{equation}
\label{Ex1.50}
 \left\{ \begin{array}{lcl}
    \hfill Ê\Delta v -{\bf a}\cdot \nabla u -a_{0}u&=& |u|^{p-2}u +f  \hbox{\rm \, on \, $\Omega$}
\\
  \hfill Ê-\alpha^{2}\Delta u -{\bf b}\cdot \nabla v-b_{0}v&=& |v|^{q-2}v +g  \hbox{\rm \, on \, $\Omega$}\\
 \hfill  u&=& u_{0} \quad \quad \quad  \quad \,\, \hbox{\rm on  $\Sigma_{-}^{1}$ }\nonumber \\
  \hfill  v&=& v_{0} \quad \quad \quad  \quad \quad \hbox{\rm on  $\Sigma_{-}^{2}$. }\nonumber \\
        \end{array}  \right.
   \end{equation}

We can use the above to get
 \begin{theorem} Assume  $ a_{0}(x)- \frac{1}{2}{\rm div} {\bf a} (x) \geq \alpha >0$   and   $b_{0}(x)- \frac{1}{2}{\rm div} {\bf b} (x) \geq \alpha >0$ on $\Omega$,  $2<p,q\leq \frac{2n}{n-2}$. For any $f, g\in L^{2}(\Omega)$ and $(u_{0}, v_{0})\in L^{2}_{A}(\Sigma^{1}_{-})\times  L^{2}_{B}(\Sigma^{2}_{-})$, consider on $H^{1}_A(\Omega)\times H^{1}_B(\Omega)$ the functional
 \begin{eqnarray}
\label{functional}
I(u, v)&=&\Psi (u) +\Psi^{*}(-{\bf a}.\nabla u -\frac{1}{2}{\rm div }({\bf a}) \,  u+\Delta v )\nonumber \\
&&+\frac{1}{2}\int_{\Sigma^{1}_{+}}|u(x)|^{2}{\bf n}(x)\cdot {\bf a}(x)\, d\sigma
-\frac{1}{2}\int_{\Sigma^{1}_{-}}|u(x)|^{2}{\bf n}(x)\cdot {\bf a}(x)\,  d\sigma \nonumber \\ 
&&+2\int_{\Sigma^{1}_{-}}u(x)u_{0}(x){\bf n}(x)\cdot {\bf a}(x)\,  d\sigma
-\int_{\Sigma^{1}_{-}}|u_{0}(x)|^{2}{\bf n}(x)\cdot {\bf a}(x)\,  d\sigma\nonumber\\
&&+\Phi (v) +\Phi^{*}(-{\bf b}.\nabla v -\frac{1}{2}{\rm div }({\bf b})\, v-\alpha^{2}\Delta u )\\
&&+\frac{1}{2}\int_{\Sigma^{2}_{+}}|v(x)|^{2}{\bf n}(x)\cdot {\bf b}(x)\, d\sigma
-\frac{1}{2}\int_{\Sigma^{2}_{-}}|v(x)|^{2}{\bf n}(x)\cdot {\bf b}(x)\,  d\sigma \nonumber \\ 
&&+2\int_{\Sigma^{2}_{-}}v(x)u_{0}(x){\bf n}(x)\cdot {\bf b}(x)\,  d\sigma
-\int_{\Sigma^{2}_{-}}|v_{0}(x)|^{2}{\bf n}(x)\cdot {\bf b}(x)\,  d\sigma, \nonumber
\end{eqnarray} 
where
\[
 \Psi (u)= \frac{1}{p}\int_{\Omega}|u|^{p} dx+\int_{\Omega}fu dx +\frac{1}{2}\int_{\Omega}(a_{0}-\frac{1}{2}{\rm div }({\bf a})) \, |u|^{2} dx, 
 \]
 \[
 \Phi (v)= \frac{1}{q}\int_{\Omega}|v|^{q} dx+\int_{\Omega}gv dx +\frac{1}{2}\int_{\Omega}(b_{0}-\frac{1}{2}{\rm div }({\bf b})) \, |v|^{2} dx, 
 \]
and $\Psi^{*}$ and $\Phi^{*}$ are their Legendre transforms. The infimum is zero and there exists a minimizer $(\bar u, \bar v)\in H^{1}_{A}(\Omega)\times H^{1}_{B}(\Omega)$  that is a solution of $(\ref{Ex1.50})$.
 \end{theorem}
 The conditions on ${\bf a}$ and ${\bf b}$ insure that   the   first order linear operators
$B_{1}u:= {\bf a}\cdot \nabla u +a_{0}u $ (resp., $B_{2}v:= {\bf b}\cdot \nabla v +b_{0}v
$) 
are positive modulo the boundary operators $u\to (u_{|_{\Sigma^{1}_{-}}}, u_{|_{\Sigma^{1}_{+}}} ) \in L_{A}^{2}(\Sigma^{1}_{-})\times L_{A}^{2}(\Sigma^{1}_{+})$ (resp., $v\to (v_{|_{\Sigma^{2}_{-}}}, v_{|_{\Sigma^{2}_{+}}} ) \in L_{B}^{2}(\Sigma^{2}_{-})\times L_{B}^{2}(\Sigma^{2}_{+})$). Apply now the above with $A=\alpha^{2}\Delta$.
    
 \section{Time dependent anti-self dual Lagrangians}
 
Let $H$ be a Hilbert space with $\braket \,{\,}\,{\,}$ as scalar product and let $[0,T]$ be a fixed real interval  where $(0<T< +\infty)$.   Consider the classical space 
$L^{2}_{H}$ of Bochner integrable functions from $[0,T]$ into $H$ with norm denoted by $\|\cdot\|_2$, as well as the Hilbert space
 \[
A^{2}_{H} = \{ u:[0,T] \rightarrow H; \,\dot{u} \in L^{2}_{H}  \}
\]
consisting of all absolutely continuous arcs $u: [0,T]\to H$, 
equipped with the norm
\[
     \|u\|_{A^{^{2}}_{H}} = (\|u(0)\|_{H}^{2} +
     \int_0^T \|\dot{u}\|^{2} dt)^{\frac{1}{2}}.
\]
   \begin{definition} Let $L:[0,T]\times H\times H\to \R \cup \{+\infty\}$ be  measurable with respect to   the  $\sigma$-field  generated by the products of Lebesgue sets in  $[0,T]$ and Borel sets in $H\times H$. 
  We say that  $L$ is an {\it anti-self dual Lagrangian} (ASD) on $[0,T]\times H\times H$  if  for any $t\in [0,T]$, the map $L_t:(x,p)\to L(t, x,p)$ is in $\ASDH$: that is if 
  \[
  L^*(t, p,x)=L(t, -x, -p) \quad \hbox{\rm for all $(x,p)\in H\times H$}. 
  \]
   where here $L^*$ is the Legendre transform in the last two variables.
  \end{definition}
The most basic time-dependent $ASD$-Lagrangians are again of the form 
\[
L(t,x,p)=\varphi (t, x) +\varphi^{*}(t, -p)
\]
where for each $t$, the function $x\to \varphi (t, x)$ is convex and lower semi-continuous.  We now show how this property naturally ``lifts'' to path space. 

\subsection{ASD Lagrangians on path spaces}

\begin{proposition} 
Suppose that $L$ is an anti-self dual Lagrangian on $[0,T]\times H\times H$, then for each $\omega \in {\bf R}$, the Lagrangian  $M (u,p):=\int_0^T e^{2wt}L(t, e^{-wt}u(t), e^{-wt}p (t)) dt$ is anti-self dual on $L^{2}_{H}$.
\end{proposition}
 \noindent{\bf Proof:} It is sufficient to show that for any Lagrangian $L(t, x, p)$, we have the formula: 
 \begin{eqnarray*}
M^{*}(p,s) &:=& \sup \left\{ \int_0^T ( \braket{p(t)}{u(t)} + \braket{s(t)}{v(t)} -
     L(t, u(t),v(t))) dt\ ;  \ (u,v) \in L_H^{2} \times L_H^{2} \right\}\\
 &=&\int_0^T L^{*}(t, p(t),s(t))dt.
 \end{eqnarray*}
 For that, first note that  for all $u, v \in L_H^{2}$  and
     $p,s \in L_{H}^{^{2}}$, we have:
\[
     \int_0^T L(t, u(t),v(t))dt + \int_0^T L^{*}(t,p(t),s(t))dt \geq \int_0^T (\braket{p(t)}{u(t)} +
     \braket{s(t)}{v(t)} dt,
\]
which implies
\begin{eqnarray*}
\int_0^T L^{*}(t, p(t),s(t))dt &\geq& \sup \left\{ \int_0^T ( \braket{p(t)}{u(t)}+ \braket{s(t)}{v(t)} -
L(t,u(t),v(t) )dt; (u,v) \in L_H^{2} \times L_H^{2} \right\}\\
   &=& M^{*}(p,s)
   \end{eqnarray*}
 For the reverse  inequality, assume 
   $M^{*}(p,s)<\int_0^T L^{*}(t, p(t), s (t)) dt$  for some $(p,s)$  in $L_{H}^{2}\times
L_{H}^{2}$, and let
 $\mu(t)$ be such that $\mu(t) < L^{*}(t, p(t),s(t))$
for all $t$ while $\int_0^T \mu(t)dt >M^{*}(p,s)$. We then have for all $t$,
\[
   -\mu(t) > -L^{*}(t, p(t),s(t)) = \inf \{
L(t, u,v)-\braket{u}{p(t)}  - \braket{v}{s(t)};  (u,v)
\in H\times H \}.
\]
By a standard measurable selection theorem (see\cite{Ca}), there exists
a  measurable pair $(u_1,u_2)
\in L_H^{2}\times L_H^{2}$ such that $
     -\mu(t) \geq L(t, u_1(t),u_2(t)) - \braket{u_1(t)}{p(t)}
     -\braket{u_2(t)}{s(t)}$. 
   Therefore
\begin{eqnarray*}
M^{*}(p,s) &<& \int_0^T \mu(t) dt\leq \int_0^T -L(t, u_1(t),u_2(t)) +
\braket{u_1(t)}{s(t)}
     +\braket{u_2(t)}{p(t)}  dt \nonumber\\
& \leq& \sup \left\{ \int_0^T
(\braket{p(t)}{u(t)} + \braket{s(t)}{v(t)} -
     L(t, u(t),v(t))) dt\ ;(u,v) \in L_H^{2} \times L_H^{2} \right\} \nonumber\\
   &=&M^{*}(p,s)
\end{eqnarray*}
 which is a contradiction.\\
  
 \noindent{\bf A representation of $(A^{^{2}}_{H})$:}  One way to represent the space  $A^{^{2}}_{H}$ is to identify it with the product space
$H \times L^{2}_{H}$, in such a way that its dual  $(A^{2}_{H})^*$ can also be
identified with $ H \times
L^{2}_{H}$  via the formula:
\[
     \braket{u}{(p_{1},p_{0})}_{_{A^{2}_{H},H \times L_{H}^{2}}} =
      \braket{u(0)}{p_{1}}_{_H} + \int_0^T \braket{\dot{u}(t)}{p_{0}(t)}dt.
\]
where $u\in A^{2}_{H}$ and $(p_{1},p_{0})\in H\times L_{H}^{2}$. 
 \begin{proposition} 
Suppose $L$ is an anti-self dual Lagrangian on $[0,T]\times H\times H$  and that  $\ell$ is  a self-dual boundary Lagrangian on $H\times H$, then the Lagrangian
 defined on $A^{2}_{H}\times (A^{2}_{H})^{*}=A^{2}_{H}\times (H \times L_{H}^{2})$ by 
  \[
{ N} (u,p)= \int_0^T L(t, u(t)+p_{0}(t),\, \dot{u}(t)) dt + \ell (u(0)+p_{1},\, u(T)) 
\]
is anti-selfdual on $A^{2}_{H}\times (L^{2}_{H}\times \{0\})$. 
  \end{proposition}
   \noindent {\bf Proof:}  For $(v, q) \in A^{2}_{H}\times (A^{2}_{H})^{*}$ with $q$ represented by $(q_{0}(t), 0)$ write:
 \begin{eqnarray*}
    { N}^{*}(q,v) &=& \sup_{p_{1} \in H} \sup_{p_{0} \in L_H^{2}} \sup_{u \in
     A^{2}_{H}} \{ \langle p_{1},v(0) \rangle + \int_0^T \left[ \braket{p_{0}(t)}{\dot{v}(t)}+\braket{q_{0}(t)}{\dot{u}(t)}
     - L(t, u(t)+p_{0}(t),\dot{u}(t))\right] dt\\
 && \quad \quad \quad \quad \quad \quad \quad -\ell (u(0)+p_{1},u(T)) \}.
\end{eqnarray*}
  Making a substitution $u(0) + p_{1} = a \in H$ and $u(t)+p_{0}(t)=y(t) \in L_H^2$,
we obtain
  \begin{eqnarray*}
    { N}^{*}(q,v) &=& \sup_{a \in H} \sup_{y \in L_H^{2}} \sup_{u \in
     A^{2}_{H}} \{ \langle a-u(0),v(0) \rangle -\ell (a,u(T))\\
  &&\quad \quad \quad \quad \quad \quad+ \int_0^T \left[ \braket{y(t)-u(t)}{\dot{v}(t)}+\braket{q_{0}(t)}{\dot{u}(t)}
     - L(t, y(t),\dot{u}(t))\right] dt 
   \end{eqnarray*}
Since $\dot u$ and $\dot v\in L_H^2$, we have:
\[
\int_0^T \braket{u}{\dot{v}}=- \int_0^T \braket{\dot u}{v} +
\braket{v(T)}{u(T)} - \braket{v(0)}{u(0)},
\]
    which implies
 \begin{eqnarray*}
    { N}^{*}(q,v) &=&\sup_{a \in H} \sup_{y \in L_H^{2}} \sup_{u \in
     A^{2}_{H}} \{ \langle a,v(0) \rangle -\braket{v(T)}{u(T)}-\ell (a,u(T))\\
  &&\quad \quad \quad \quad \quad + \int_0^T \left[ \braket{y(t)}{\dot{v}(t)}+\braket{v(t)+q_{0}(t)}{\dot{u}(t)} 
     - L(t, y(t),\dot{u}(t))\right] dt \}.
  \end{eqnarray*}
Identify now $A_{H}^{^2}$ with $H\times L_H^2$
via the correspondence:
  \begin{eqnarray*}
     (b,r) \in H\times L_H^2 &\mapsto & b+\int_t^T r(s)\, ds\in
A_{H}^{^2}\\
     u\in A_{H}^{^2}   &\mapsto & \big( u(T),-\dot u(t)\big)\in H\times
L_H^2 .
\end{eqnarray*}
We finally obtain
  \begin{eqnarray*}
    { N}^{*}(q,v) &=&\sup_{a \in H}\sup_{b \in H}  \{ \langle a,v(0) \rangle -\braket{v(T)}{b}-\ell (a,b)\\
  &+& \sup_{y \in L_H^{2}} \sup_{r \in
     L^{2}_{H}}\int_0^T \left[ \braket{y(t)}{\dot{v}(t)}+\braket{v(t)+q_{0}(t)}{r(t)} 
     - L(t, y(t),r(t))\right] dt\\
  &=&\int_{0}^{T}L^{*}(t, \dot{v}(t),v(t)+q_{0}(t)) dt +\ell^{*}(v(0), -v(T))\\
   &=&\int_{0}^{T}L(t, -v(t)-q_{0}(t), -\dot{v}(t)) dt +\ell (-v(0), -v(T))\\
   &=& { N}(-v, -q). 
  \end{eqnarray*}
  \subsection{ASD Lagrangians in the calculus of variations}
  
    \begin{theorem} Suppose $L$ is an anti-self dual Lagrangian on $[0,T]\times H\times H$ and $\ell $ is  a self-dual boundary Lagrangian on $H\times H$, and consider the following functional
    \[
    I_{\ell,L}(u) = \int_0^T L(t, u(t),\dot{u}(t))dt + \ell (u(0),u(T)).
    \]
     Suppose  there exists $C>0$ such that for all $x\in L^{2}_{H}$, $
     \int_{0}^{T}L(t,x(t),0)dt \leq  C(1+\|x\|^{2}_{L^{2}_{H}}).$
    Then there exists $v \in A_{H}^2$ such that
$\big(v(t),\dot{{v}} (t)\big)\in \mbox{\rm Dom} (L)$ for almost all 
$t\in [0,T]$    and 
\[ I_{\ell,L}( v)=\inf\limits_{u\in A_{H}^2}I_{\ell, L}(u)=0.
\]
  In particular, for every $v_{0}\in H$ the following functional 
\[
 I_{\ell,L}(u)=  \int_0^T L(t, u(t),\dot{u}(t))dt +\frac{1}{2}\|u(0)\|^{2} -2\langle v_0, u(0)\rangle +\|v_0\|^{2} +\frac{1}{2}\|u(T)\|^{2}
\]
has minimum equal to zero on $A_{H}^2$. It  is attained at a unique path $v$ which then satisfies:
\begin{equation}
v(0)=v_{0}\,\, {\rm and }\,\, \big(v(t),\dot{v} (t)\big)\in \mbox{\rm Dom} (L) \hbox{\rm \quad for almost all
$t\in [0,T]$},
\end{equation}
\begin{equation}
\label{eqn:2.10}
\frac{d}{dt}\partial_{p}L(t, v(t), {\dot v}(t)) =\partial_{x} L (t, v(t), {\dot v}(t))
\end{equation}
\begin{equation}
\label{eqn:2.20}
 (-{\dot v}(t), -v(t))\in \partial L(t, v(t), {\dot v} (t)),
\end{equation}
  \begin{equation}
\label{eqn:2.30}
 \|v(t)\|_{H}^{2}=\|v_{0}\|^{2}-2\int_{0}^{t}L(s, v(s), {\dot v}(s)) ds \quad \hbox{\rm for every $t\in [0,T]$.}
\end{equation}
If $L$ is autonomous and $v\in C^1([0,T], H)$, then for all $t\in [0,T]$, we have:
\begin{equation}
\label{eqn:2.31}
\|\dot v (t)\|\leq \|\dot v (0)\|.
\end{equation}
 \end{theorem} 
 \noindent{\bf Proof:} Apply Proposition 6.2 to get that  
 \[
{ N} (u,p)= \int_0^T L(t, u(t)+p_{0}(t),\, \dot{u}(t)) dt + \ell (u(0)+p_{1},\, u(T)) 
\]
is partially anti-self dual on $A^{2}_{H}$. It now suffices to apply Theorem 4.1 since in this case
$ N(0, p) =\int_0^T L(t, p_{0}(t),0) dt + \ell (p_{1},\, 0) \leq  C_{2}(1+\|p_{0}\|^{2}_{L^{2}_{H}})+\|p_{1}\|_H^{2}$, which means that $N(0, p)$ is bounded on the bounded sets of  $(A_{H}^2)^{*}$. \\
For a given $v_{0}\in H$, use the boundary Lagrangian 
\[
\ell (r,s)= \frac{1}{2}\|r\|^{2} -2\langle v_{0}, r\rangle +\|v_{0}\|^{2} +\frac{1}{2}\|s\|^{2}. 
\]
which is clearly self-dual. We then get
 \[
I_{\ell,L}(u) = \int_0^T \left[ L(t, u(t),\dot{u}(t))+\langle u(t), {\dot u}(t)\rangle \right] dt +  \|u(0)-v_{0}\|^{2}.
\]
 Since $ L(t,x,p)\geq - \langle x,p\rangle $ for all $(t,x,p)\in [0,T]\times H\times H$,
  the fact that $ I_{\ell,L}( v)=\inf\limits_{u\in A_{H}^2}I_{\ell, L}(u)=0$,  then yields
$v(0)=v_{0}$ and that
 \begin{equation}
 \label{weird}
  L(s, v(s),\dot{v}(s)+\langle v(s), {\dot v}(s)\rangle=0 \quad \hbox {\rm for almost all $s\in [0,T]$.}
  \end{equation}
     This clearly yields (\ref{eqn:2.30}), since we then have:
   \[
   \frac{d (|v(s)|^{2})}{ds}=-2L(s, v(s),{\dot v}(s)). 
   \]
 To prove (\ref{eqn:2.20}), use (\ref{weird}) and the fact that $L$ is anti-selfdual to write:
 \[
  L(s, v(s),\dot{v}(s)+L^{*}(s, -\dot{v}(s), -v(s)) +\langle (v(s), {\dot v}(s)), ({\dot v}(s), v(s) )\rangle=0.
 \]   
Now apply Legendre-Fenchel  duality in the space $H\times H$.  The uniqueness and (\ref{eqn:2.31}) follow from the following observation.
  
\begin{lemma} Suppose $L(t,\, , \,)$ is convex on $H\times H$ for each $t\in [0, T]$, and that $x(t)$ and $v(t)$ are two paths in $C^1([0,T], H)$ satisfying $x(0)=x_0$, $v(0)=v_0$, $-(\dot x, x)\in \partial L(t, x,\dot x)$ 
and $-(\dot v, v)\in \partial L(t, v,\dot v)$. 
Then $\| x(t)-v(t)\| \leq \| x(0)-v(0)\|$ for each $t\in [0, T]$. 
\end{lemma} 
\noindent{\bf Proof:} Estimate $\alpha (t)=\frac{d}{dt}\frac{{\| x(t)-v(t)\|}^2}{2}$ as follows:
{\small
\begin{eqnarray*}
\alpha (t)
&=&\langle v(t)-x(t),\dot v(t)-\dot x(t)\rangle \\
&=&\frac{1}{2}\langle v(t)-x(t),\dot v(t)-\dot x(t)\rangle +
\frac{1}{2}\langle\dot v(t)-\dot x(t),v(t)-x(t)\rangle\\
&=&\frac{1}{2}\left(\langle (v(t)-x(t),\dot v(t)-\dot x(t)),
(\dot v(t)-\dot x(t),v(t)-x(t))\rangle_{H\times H}\right)\\
&=&\frac{1}{2}\left(\langle (v(t)-x(t),\dot v(t)-\dot x(t)),(L_x(t, x(t),\dot x(t))-
L_x(t, v(t),\dot v(t)),L_y(t, x(t),\dot x(t))-L_y(t, v(t),\dot v(t))\rangle\right)\\
&=&\frac{1}{2}\left(\langle (v(t),\dot v(t))-(x(t),\dot x(t)),
\big( L_x(t, x(t),\dot x(t)),L_y(t, x(t),\dot x(t))\big) -
\big( L_x(t, v(t),\dot v(t)),L_y(t, v(t),\dot v(t))\big)\rangle\right)\\
&=&\frac{1}{2}\left(\langle (v(t),\dot v(t))-(x(t),\dot x(t)),
\partial L(t, x(t),\dot x(t))-\partial L(t, v(t),\dot v(t))\rangle\right)\\
 &\leq& 0 
\end{eqnarray*}
}
in view of the convexity of $L$.

 It then follows that $\| x(t)-v(t)\| \leq \| x(0)-v(0)\|$ for all $t>0$. Now if $L$ is autonomous, $v(t)$ and $x(t)=v(t+h)$ are solutions for any $h>0$, so that (\ref{eqn:2.31}) follows from the above. 

\subsection{ASD Lagrangians associated to gradient flows}

The most basic example of a self-dual Lagrangian already provides a variational formulation and proof of existence for gradient flows. The following extends some of the resuts in \cite{GT1}.
\begin{theorem} \label{evolving_descent}
Let $\phi :[0,T]\times H \rightarrow {\bf R} \cup \{+\infty\}$
be a measurable function with respect to   the  $\sigma$-field in
$[0,T]\times H$ generated by the products of Lebesgue sets in  $[0,T]$
and
Borel sets in $H$. Assume that for every $t\in [0,T]$, the function $\phi(t,\cdot)$ is
convex and lower semicontinuous on $H$, and $A_{t}$ is a bounded linear positive operator on $H$ such that for some positive functions $\gamma, \beta^{-1} \in L^{\infty}[0,T]$, we have 
\begin{equation}
\label{bounds}
 \beta (t)\|x\|^{p}\leq \phi (t,x) +\frac{1}{2}\langle A_{t}x, x\rangle\leq  \gamma (t)\|x\|^{q}.
\end{equation}
Then, for any $u_{0}\in H$, the functional
\begin{equation}
\label{functional2}
 I(u)=\frac{1}{2}(|u(0)|^{2} + |u(T)|^{2})
-2\langle u(0), u_{0} \rangle +|u_{0}|^{2} +
\int_{0}^{T} \left[\psi (t, u(t)) + \psi^{*}(t,-A^{a}_{t}u(t)-\dot{u}(t))\right]dt
\end{equation}
where $\psi $ is the convex functiona $\psi (t,x)=\phi (t,x)+\frac{1}{2}\langle A_{t}x, x\rangle$ has  a unique minimizer $v$ in $\AH$ such that: 
\begin{equation}
\label{zero}
    I(v)=\inf\limits_{u \in \AH} I(u)=0.
\end{equation}
Among the paths in $\AH$, $v$ is the unique solution to
\begin{equation}
\label{floweq}
\left\{ \begin{array}{lcl}
             -A_{t}u(t)- \dot{v}(t) &\in& \partial \phi (t, v(t))  \quad{\rm
a.e.\quad
on}\quad
[0,T]\\
\hfill v(0) &=& u_{0}.
         \end{array}  \right.
\end{equation}
\end{theorem}

\noindent {\bf Proof:}  This follows directly from Theorem 6.2 applied to the anti-selfdual Lagrangian $L(t, x, p)=\psi (t,x) +\psi^{*}(t,-A_{t}^{a}x-p)$. Note that the conditions (\ref{bounds}) yield that $\int_{0}^{T}L(t, x(t), 0) dt=\int_{0}^{T}\psi (t,x(t))+\psi^{*}(t,A^{a}_{t}x(t)) dt $ is bounded on the bounded sets of $L^{2}_{H}$. 
 
 \subsection{Variational resolution of parabolic equations with prescribed boundaries} 
Suppose now that for each $t\in [0,T]$,  $(b^t_{1}, b^t_{2}): X_{t}\to H^t_{1}\times H^t_{2}$ are {\it regular boundary operators} from a reflexive Banach space $X_{t}$ into Hilbert spaces $H^t_{1}, H^t_{2}$, and that there are operators 
$\Lambda_t: X_{t}\to X_{t}^{*}$ which are skew-adjoint modulo the boundary   $(b^t_{1}, b^t_{2})$, that is for every $x,y \in X_{t}$, we have:
\[
\langle \Lambda_t x , y\rangle_{X_{t}} = -\langle \Lambda_t y , x\rangle_{X_{t}}  + \langle b^t_{2}(x), b^t_{2}(y) \rangle_{_{H_{2}^t}} -\langle b^t_{1}(x), b^t_{1}(y)\rangle_{_{H_{1}^t}}.
\]
Suppose $H$ is a Hilbert space such that for each $t$, $(X_{t}, \Lambda_{t}, H)$ is a maximal evolution triple, in particular $X_{t}\subset H \subset X_{t}^{*}$ and $\Lambda_{t}:X_{t}\to H$. Now starting with a time-dependent ASD Lagrangian $L$ on $H$, and self-dual {\it state-boundary Lagrangians} $m_t :H^t_{1}\times H^t_{2} \to \R \cup \{+\infty\}$, we get by Proposition 5.4 that 
   \[
M (t, x,p)= \left\{ \begin{array}{lcl}
  \hfill L(t, x, \Lambda_{t} x+p) +m_{t} (b^{t}_{1}(x), b^{t}_{2}(x))\, &{\rm if}& {x\in X_{t}} \\
  +\infty &\,& {\rm otherwise} \\
\end{array}\right.
\]
 is also anti-self dual on  $H\times H$ for each $t\in [0,T]$.
 
 If now $\ell$ is a self-dual time-boundary Lagrangian on $H$, then   
 \[
{\tilde M} (u,p)= 
 \int_0^T \left\{M(t, u(t),p(t)+\dot{u}(t)) \right\}dt + \ell (u(0),u(T)) 
  \]
 is partially anti-self dual Lagrangian on $A^{2}_{H}$, and Theorem 6.2 then applies to get that 
 \[
I(u)={\tilde M} (u,0)=\int_0^T \left\{L(t, u(t),\Lambda_{t} u (t)+\dot{u}(t)) + m_{t} (b^{t}_{1}u(t), b^{t}_{2}u(t))\right\}dt + \ell (u(0),u(T)) 
\]
has a minimum at ${\bar v} (t)$, and that the minimal value is zero. Applying the theorem with the time boundary Lagrangian on $H$,
\[
\ell (x,p) =  \frac{1}{2}\|x\|^{2} -2\langle v_{0}, x\rangle +\|v_{0}\|^{2}+\frac{1}{2}\|p\|^{2} 
\]
where  $v_{0}$ is a given initial value in $H$, and with a state boundary Lagrangian  
\[
m_{t}(x,p) =  \frac{1}{2}\|x\|^{2} -2\langle \gamma (t), x\rangle +\|\gamma (t)\|^{2}+\frac{1}{2}\|p\|^{2}, 
\]
where  $\gamma (t)$ is prescribed in $H^{t}_{1}$ for each $t$, we get that ${\bar v} (t)$ satisfies:
  \begin{equation}
 \left\{ \begin{array}{lcl}
\label{eqn:?}
 \hfill  L(t, v(t),\Lambda_{t} v (t)+\dot{v}(t)) +\langle v(t), \Lambda_{t} v(t)+\dot{v}(t) \rangle &= &0 \quad \quad \hbox {\rm a.e. $t\in [0, T]$}\\
\hfill  (-\Lambda_{t} v (t)-{\dot v}(t), -v(t))&\in& \partial L(t, v(t), {\dot v} (t))\\
\hfill b^{t}_{1}(v(t))&=&\gamma (t) \quad \hbox{\rm a.e $t\in [0, T]$}\\
\hfill v(0)&=&v_{0}
\end{array}\right.
\end{equation}
By starting with the most basic Lagrangian 
$
L(t,x,p)=\varphi (t, x) +\varphi^{*}(t, -p)
$
 we get
 \begin{theorem}
 Under the above conditions on $(X_{t}, H, H^{t}_{1}, H^{t}_{2}, b^{t}_{1}, b^{t}_{2})$, consider  bounded linear operators $A_{t}:X_{t}\to X_{t}^{*}$ such that $A_{t}-\frac{1}{2}((b^{t}_{2})^{*}b^{t}_{2} -(b^{t}_{1})^{*}b^{t}_{1})$ is positive and denote by  $\Lambda_{t}$ the operator $
 \Lambda_{t}=\frac{1}{2}(A_{t} -A^{*}_{t})+\frac{1}{2}((b^{t}_{2})^{*}b^{t}_{2} -(b^{t}_{1})^{*}b^{t}_{1})
$ which is skew-adjoint modulo the boundary. \\
For each $t\in [0, T]$, suppose $(X_{t}, H, \Lambda_{t})$ is a {\it maximal evolution triple}  and that $\varphi (t, \cdot)$ is a convex continuous function on  $H$. For $f\in L^{2}([0,T];H)$, $v_{0}\in H$ and $\gamma (t)\in H^{1}_{t}$ consider the following functional on $A_{H}^2$, 
\begin{eqnarray*}
 I(u) &=& \int_0^T\left\{\psi(t, u(t)) +\psi^{*}(t, -\Lambda_{t} u (t)-\dot{u}(t))+
 \frac{1}{2} (|b^{t}_{1}u(t)|^{2} + |b^{t}_{2}u(t)|^{2})-2\langle \gamma (t), u(t)\rangle +|\gamma (t)|^{2})\right\} dt\\
 && + \frac{1}{2}(|u(0)|^{2}+|u(T)|^{2})-2\langle u(0), v_{0} \rangle +|v_{0}|^{2}, 
\end{eqnarray*}
where  $
\psi (t, x)=\varphi (t, x) +\frac{1}{2} \langle A_{t}x, x\rangle  - \frac{1}{4}( \|b^{t}_{2}x\|^{2} -\|b^{t}_{1}x\|^{2})+\langle f(t), x\rangle$. 
 Suppose  there is $C>0$ so that for every $x\in A^{2}_{H}$, 
 \[
 \int_{0}^{T}\psi (t,x(t))+\psi^{*}(t,-\Lambda_{t} x (t) dt \leq C (1+\|x\|_{L^{2}_{H}}^2).
 \]
  Then there exists
$v \in A_{H}^2$ such that 
 $ I ( v)=\inf\limits_{u\in A_{H}^2}I (u)=0.$
 Moreover, $v$ solves 
 \begin{equation}
 \left\{ \begin{array}{lcl}
\label{eqn:?}
 \hfill   -A_{t} v(t)-\dot{v}(t)  &\in &\partial \varphi (t, v(t)) +f(t) \quad \hbox {\rm a.e. $t\in [0, T]$}\\
\hfill b^{t}_{1}(v(t))&=&\gamma (t) \quad \quad \quad \quad \quad \quad \quad  \hbox{\rm a.e $t\in [0, T]$}\\
\hfill v(0)&=&v_{0}.
\end{array}\right.
\end{equation}
\end{theorem}
 
\subsubsection*{Example 9: Non linear Transport evolutions} With the notation of Example 6, we consider the equation
\begin{equation}
\label{transport.eq.10}
 \left\{ \begin{array}{lcl}
    \hfill -Ê\frac {\partial u}{\partial t}-\Sigma_{i=1}^{n}a_{i}\frac{\partial u}{\partial x_{i}} -a_{0}u&=& \beta (u) +f  \hbox{\rm \, on \, $[0,T]\times \Omega$}
\\
 \hfill  u(t, x) &=& \gamma (t,x) \quad \quad  \hbox{\rm on \quad $[0,T]\times \Sigma_{-}$. }\nonumber \\
 \hfill u(0,x)&=&u_{0}(x) \quad \quad \hbox{\rm on $\Omega$} 
       \end{array}  \right.
   \end{equation}
 where $u_{0}\in H_{A}^{1}(\Omega)$, $f \in L^{2}(\Omega)$ and where $\gamma (t)\in L^{2}(\Sigma_-, |{\bf n}\cdot {\bf a}|dx)$ for each $t\in [0, T]$. 
Let 
\[
\psi (u)=\int_{\Omega}\left\{j(u(x))+f(x)u(x)+\frac{1}{2}(a_{0}-\frac{1}{2}{\bf div \, a }) |u|^{2})\right\}dx
 \]
 \begin{theorem}  Assume  
 $a_{0}(x)- \frac{1}{2}{\bf div} a (x) \geq \alpha >0$ on $\Omega$, 
and consider  the following functional on the space $X:=A^{2}([0,T]; H^{1}_{A}(\Omega))$.
{\small 
\begin{eqnarray*}
{I}(u)&=&\int_{0}^{T}\left\{\psi (u(t))+\psi^{*}(-{\bf a}\cdot \nabla_{x} u(t)- \frac{1}{2}({\rm div \, a })\, u (t)-\dot u(t)) \right\} dt\\
&&+\int_{0}^{T}\left\{\frac{1}{2}\int_{\Sigma_{+}}|u(t, x)|^{2}{\bf n}\cdot {\bf a}\, d\sigma-\frac{1}{2}\int_{\Sigma_{-}}|u(t,x)|^{2}{\bf n}\cdot {\bf a}\,  d\sigma+\int_{\Sigma_-}(|\gamma (t,x)|^{2}-2\gamma(t,x)u(t,x))|{\bf n}\cdot {\bf a}|\,d\sigma\right\} dt
\\
&& + \int_{\Omega}\left\{\frac{1}{2}(|u(0,x)|^{2}+|u(x, T)|^{2})-2\langle u(0,x), u_{0} (x)\rangle +|u_{0}(x)|^{2}\right\}dx.
\end{eqnarray*}
}
There exists $\bar u\in X$ such that
 $I(\bar u)=\inf_{u\in X}I(u)=0
$
and which solves equation  (\ref{transport.eq.10}).
\end{theorem}
These results will be improved in \cite{GT2}. 
\subsubsection*{Variational resolution for parabolic-elliptic  variational inequalities} 

Consider for each time $t$, a bilinear continuous functional $a_t$ on a Hilbert space $H\times H$  and a convex l.s.c function $\varphi (t, \cdot ) :H\to {\bf R}\cup \{+\infty\}$.  Solving the corresponding  parabolic variational inequality amounts to constructing for a given $f\in L^2([0,T];H)$ and $x_0\in H$, 
a path ${ x (t)} \in A^2_H([0,T])$ such that for all $z\in H$, 
\begin{equation}
\label{varineq1}
\langle {\dot x}(t), { x(t)} -z)+a_t({ x}(t) ,{ x(t)} -z)+\varphi (t, { x} (t))-\varphi (t, z) \leq \langle{ x}(t) -z, f (t)\rangle.
\end{equation}
for almost all  $t\in [0,T]$.  This problem can be rewritten as: $
f(t)\in {\dot y}(t) +A_ty(t)+\partial \varphi (t, y)$, 
where $A_t$ is the bounded linear operator on $H$ defined by $a_t(u,v)=\langle A_tu, v \rangle$. This means that the variational inequality (\ref{varineq1}) can be rewritten and solved using the variational principle in Theorem 6.5  For example, one can then solve variationally the following ''obstacle '' 
problem.

\begin{corollary} Let $(a_t)_t$ be bilinear continuous functionals on $H\times H$ satisfying:
\begin{itemize}
\item  For some $\lambda>0$, $a_t(v,v)\geq \lambda \|v\|^{2}$ on $H$ for every $t\in [0,T]$. 

\item  The map $u\to \int_{0}^{T}a_{t}(u(t), u(t)) dt$ is continuous on $L^{2}_{H}$.
\end{itemize}
 If $K$ is a convex closed subset of $H$, then for any $f\in L^2([0,T];H)$ and any $x_0\in K$,  there exists a path $x\in A^2_H([0,T])$ such that $
  x(0)=x_0, \, \hbox{\rm  $ x (t)\in K$  for almost all $t\in [0, T]$ and}$ 
  \[
\langle {\dot x}(t), { x(t)} -z \rangle +a_t({ x}(t) ,{ x(t)} -z) \leq \langle  x(t)-z, f\rangle \quad \hbox{\rm for all $z\in K$}.
\]
The path $ x (t)$  is  obtained as  a minimizer of the following  functional  on $A^2_H([0,T])$:
{\small
\[
I(y)=\int_0^T\left\{\varphi (t, y(t)) +(\varphi (t, \cdot)+\psi_{K})^{*}(-{\dot y}(t)-\Lambda_t y(t)) \right\}dt + \frac{1}{2}(|y(0)|^{2}+|y(T)|^{2})-2\langle y(0), x_0 \rangle +|x_0|^{2}. 
\]
}
 Here $\varphi (t, y)=\frac{1}{2}a_t(y,y)-\langle f(t), y\rangle$ and $\psi_{K} (y)=0$ on $K$ and $+\infty$ elsewhere, while $\Lambda_t:H\to H$ is the skew-adjoint operator  defined by $
\langle \Lambda_t u,v\rangle=\frac{1}{2}(a_t(u,v)-a_t(v,u))$.
\end{corollary}

\section{Semi-groups associated to autonomous anti-selfdual Lagrangians}

When the Lagrangian $L(x, p)$ is autonomous, the situation is much nicer since we can associate a flow without stringent boundedness or coercivity  conditions. Indeed, we can then use a Yosida-type regularization of ASD-Lagrangian reminiscent of the standard theory for operators and for convex functions. Let us define the {\it Partial Domain} of $\partial L$ to be the set:
\[
{\rm Dom}_{1}(\partial L)=\{x\in X;\hbox{\rm  there exists $p,q\in X^{*}$ such that $(p,0)\in \partial L(x, q)$\}.}
\]
Note that if $L(x,p)=\phi (x) +\phi^{*}(-p)$ with $0$ assumed to be in the domain of $\partial \phi$, then $x_0$ belongs to  ${\rm Dom}_{1}(\partial L)$ if and only if it belongs to the domain of $\partial \phi$. 

We then obtain the following result.

\begin{theorem}   Let $L$ be an anti-selfdual Lagrangian on a Hilbert space $H$ that is uniformly convex in the first variable. Assuming ${\rm Dom}_{1}(\partial L)$ is non-empty, then there exists a semi-group of 1-Lipschitz maps  $(T_{t})_{t\in {\bf R}^{+}}$ on $H$  such that $T_{0}=Id$ and for any $x_0\in {\rm Dom}_{1}(\partial L)$, the path $x(t)=T_tx$ satisfies the following:
\begin{equation}
\label{eqn:2.1}
\frac{d}{dt}\partial_{p}L(x(t), {\dot x}(t)) =\partial_{x} L (x(t), {\dot x}(t))
\end{equation}
\begin{equation}
\label{eqn:2.2}
 (-{\dot x}(t), -x(t))\in \partial L(x(t), {\dot x} (t))
\end{equation}
   and
\begin{equation}
\label{eqn:2.3}
 \|x(t)\|_{H}^{2}=\|x\|^{2}-2\int_{0}^{t}L(x(s), {\dot x}(s)) ds \quad \hbox{\rm for every $t\in [0,T]$.}
\end{equation}
The path $x=(x(t))_t=(T_tx)_t$ is obtained as a minimizer on $A_{H}^2$ of the functional
\[
 I(u)=  \int_0^T L(u(t),\dot{u}(t))dt +\frac{1}{2}\|u(0)\|^{2} -2\langle x, u(0)\rangle +\|x\|^{2} +\frac{1}{2}\|u(T)\|^{2}, 
\]
where $I(x)=\inf\limits_{u\in A^2_H}I(u)$
 \end{theorem} 
 As mentioned above, we can associate to the Lagrangian $L(x, p)$ its $\lambda$-regularization by considering $L_{\lambda}=L\star T_{\lambda}$ where 
$T_{\lambda}(x, p)=\frac{\|x\|^{2}}{2\lambda^{2}}+\frac{\lambda^{2}\|p\|^{2}}{2}$. 
Then $L_{\lambda}$ satisfies the hypothesis of Theorem 6.2, and we can then find 
for each initial point $v\in H$, a path $v_{\lambda} \in A_{H}^2$, with  $v_{\lambda}(0)=v$,  which verify the above properties.  \\
The uniform convexity of $L$ in the first variable  insures that the regularization $L_\lambda$ is uniformly convex in both variables which then yield $C^1$-solutions. The 1-Lipschitz property follows from Lemma 4.3, since in the autonomous case, we can apply it to a solution $u(t)$ and its translate $v(t)=u(t+h)$ to  get $\| u(t+h)-u(t)\| \leq \| u(h)-u(0)\|$ for all $t$, which yields 
that $\lim\limits_{h\to 0}\frac{\| u(t+h)-u(t)\|}{h} \leq \lim\limits_{h\to 0}\frac{\| u(h)-u(0)\|}{h}$. 
The rest of the argument amounts to analyzing what happens when  $\lambda \to 0$.  The details will be given in \cite{GT2}. \\

We can also deal with the following situation which can sometimes do away with coercivity assumptions and to also cover the case of semi-convex potentials.

\begin{theorem} Let $L$ be an anti-selfdual Lagrangian on a Hilbert space $H$ that is uniformly convex in the first variable. Assuming ${\rm Dom}_{1}(\partial L)$ is non-empty, then for any $\omega \in {\bf R}$ there exists a semi-group of  maps $(T_{t})_{t\in {\bf R}^{+}}$ on $H$  such that:  $T_{0}=Id$ and $\|T_tx-T_ty\|\leq e^{-\omega t}\|x-y\|$ for any $x,y\in H$. Moreover, for any $x_0\in {\rm Dom}_{1}(\partial L)$ the path $x(t)=T_tx_0$ satisfies the following:
 \begin{eqnarray}
\label{eqn:2.2}
 -({\dot x}(t)+\omega x(t), x(t))&\in& \partial L(x(t), 
{\dot x} (t)+\omega x(t))\\
x(0)&=&x_0. \nonumber
\end{eqnarray}
 The path $x(t)$  is obtained as a minimizer on $A_{H}^2$ of the functional
\[
\tilde  I(u)=  \int_0^T e^{2\omega t} L(u(t), \omega u (t)+ {\dot u}(t))dt +\frac{1}{2}\|u(0)\|^{2} -2\langle x_0, u(0)\rangle +\|x_0\|^{2} +\frac{1}{2}\|e^{\omega T}u(T)\|^{2}
\]
in such a way that $\tilde I(x)=\inf\limits_{u\in A^2_H}\tilde I(u)=0.$
 \end{theorem} 
 \noindent{\bf Proof:} We associate to $L$, the anti-selfdual Lagrangian 
 \[
 L_{\omega}(t, x,p):=(e^{\omega t} {\bf \cdot} L)(x,p)=e^{2\omega t} L(e^{-\omega t} x, e^{-\omega t} p).
 \]  
Note that if $y(t)$ satisfies:
\begin{equation}
\label{eqn:007}
 (-{\dot y}(t), -y(t))\in \partial L_{\omega}(t, y(t), {\dot y} (t))
\end{equation}
then $x(t)=e^{-\omega t} y(t)$ satisfies
\begin{equation}
\label{eqn:008}
- ({\dot x}(t)+\omega x(t), x(t))\in \partial L(x(t), {\dot x} (t)+\omega x (t))
\end{equation}
However, we cannot apply Theorem 7.1 directly to the Lagrangian $L_{\omega}$ because the latter is not autonomous. However, we shall see in \cite{GT2} that the Yosida regularization argument still works in this case, since we have the following property:
\[
(e^{\omega t}{\bf \cdot} L) \star M_{\lambda}=e^{\omega t}{\bf \cdot} (L\star M_{\lambda}).
\]
Now we can deduce the following which was established in \cite{GT1} in the case of gradient flows of convex potentials (i.e., when $A=0$ and $\omega=0$), and in \cite{GM} in the case of gradient flows of semi-convex functions (i.e., when $A=0$ and $\omega >0$). 
\begin{theorem}
  Let $\phi$ be a proper,  bounded below, convex lower semi-continuous functional on $H$ such that $0\in {\rm Dom} \partial \phi$ and let $A$ be a positive bounded linear operator  on $H$. For any $\omega\in {\bf R}$ and $x_0\in {\rm Dom} \partial \phi$, consider the following functional on  $A^2_H$:
  \begin{eqnarray*}
 I(u)&=&\int_0^Te^{2\omega t}\left\{\psi (u(t))+\psi^*(-A^au(t) -\omega u(t)-{\dot u(t))}\right\}dt\\
 && +\frac{1}{2}\|u(0)\|^{2} -2\langle x_0, u(0)\rangle +\|x_0\|^{2} +\frac{1}{2}\|e^{\omega T}u(T)\|^{2}
 \end{eqnarray*}
where $A^a$ is the anti-symmetric part of $A$, and $\psi (u)=\phi(u)+\frac{1}{2}\langle Au, u\rangle$.
The minimum of $I$ is then zero and is attained at a path $x(t)$ which  is a solution of 
\begin{equation}
 \left\{ \begin{array}{lcl}
\label{eqn:?}
 \hfill   -A x(t) -\omega x(t) -\dot{x}(t)  &\in &\partial \varphi (x(t)) \quad \hbox {\rm a.e. $t\in [0, T]$}\\
 \hfill x(0)&=&x_{0}.
\end{array}\right.
\end{equation}
  \end{theorem} 
\subsection{Nonlinear parabolic equations}  
\subsubsection*{Example 10: Quasi-linear parabolic equations}

Let $\Omega$ be a smooth bounded domain in $\R^n$.
For $p \geq \frac{n-2}{n+2}$,  the Sobolev space   $W^{1,p+1}_0(\Omega)\subset H:=L^2(\Omega)$, and so we define on $L^{2}(\Omega)$ the functional
\begin{equation}
\varphi (u)=
\left\{ \begin{array}{lcl}
 \frac{1}{p+1} \int_{\Omega} | \nabla u|^{p+1} 
& {\rm on} \quad  
W^{1,p+1}_0(\Omega)\\
    +\infty   &{\rm elsewhere}
\end{array}\right.
\end{equation}
  Its conjugate is then
\begin{equation}
\varphi^* (v) =
\frac{p}{p+1} \int_{\Omega} |\nabla \Delta_{p}^{-1} v|^{\frac{p+1}{p}} dx.
\end{equation}
 then for any $\omega \in \R$, any $u_0 \in W_{0}^{1,p+1} (\Omega)$ and any $f\in
W^{-1,\frac{p+1}{p}}(\Omega)$, that the infimum of the
     functional
\begin{eqnarray*}
     I(u) &=&\frac{1}{p+1} \int_0^T e^{2 \omega t}\int_{\Omega} \left(
 | \nabla u(t,x)|^{p+1} -(p+1)f(x)u(x,t)\right)dx
         dt \nonumber \\
 && + 
        p\int_0^Te^{2 \omega t} \int_{\Omega} \left(|\nabla \Delta_{p}^{-1}
(f(x)-\omega u(t,x)-\frac{\partial u}{\partial t} (t,x)))|^{\frac{p+1}{p}}\right)dx
         dt \nonumber \\
&&-2\int_{\Omega}u(0,x)u_0(x)\, dx+\int_{\Omega}|u_0(x)|^{2}\, dx+
\frac{1}{2}  \int_\Omega (|u(0,x)|^2+e^{2T}|u(T,x)|^2)dx
   \end{eqnarray*}
   on the space  $A^2_H$  is equal to zero and is attained
uniquely at an $W^{1,p+1}_0(\Omega)$-valued path $u$ 
such that
$\int^T_0
\|\dot{u}(t)\|_2^2 dt < +\infty$ and which is a  solution of the
     equation:
\begin{equation}
   \left\{ \begin{array}{lcl}
   \hfill  \frac{\partial u}{\partial t} &=& \Delta_{p} u +\omega\, u +f  \hbox{\rm \, on \, $\Omega\times
[0,T]$} \\
 \hfill  u(0,x) &=& u_0 \quad {\rm on} \quad \Omega \nonumber \\
 \hfill  u(t,0)&=&0 \quad {\rm on} \quad \partial \Omega. \nonumber
          \end{array}  \right.
   \end{equation}
 Similarly, we can deal with the equation
 \begin{equation}
   \left\{ \begin{array}{lcl}
    \hfill  \frac{\partial u}{\partial t}(t,x) &=& \Delta_{p} u -Au +\omega u(t,x) +f  \hbox{\rm \, on \, $\Omega\times
[0,T]$} \\
 \hfill  u(0,x) &=& u_0 \quad {\rm on} \quad \Omega \nonumber \\
 \hfill  u(t,0)&=&0 \quad {\rm on} \quad \partial \Omega. \nonumber
          \end{array}  \right.
   \end{equation}
whenever $A$ is a positive operator on  $L^{2}(\Omega)$, 
\subsubsection*{Example 11: Porous media equations}
Let $H=H^{-1}(\Omega)$ equipped with the norm induced by the scalar product
\[
     \braket{u}{v} = \int_{\Omega} u(-\Delta)^{-1}v dx =
     \braket{u}{v}_{H^{-1}(\Omega)}.
\]
For $m \geq \frac{n-2}{n+2}$, we have $L^{m+1} (\Omega) \subset H^{-1}$, and so we can consider  the functional
\begin{equation}
\varphi (u)=
\left\{ \begin{array}{lcl}
      \frac{1}{m+1} \int_\Omega |u|^{m+1}
& {\rm on} \quad  
L^{m+1} (\Omega)\\
    +\infty   &{\rm elsewhere}
\end{array}\right.
\end{equation}
and its conjugate
\begin{equation}
   \varphi^* (v) =
   \frac{m}{m+1} \int_{\Omega}
|  \Delta^{-1}v|^{\frac{m+1}{m}} dx.
\end{equation}
Then, for any $\omega \in \R$,  $u_0 \in H^{-1} (\Omega)$ and 
$f \in L^{2}(\Omega)$, the infimum of the
     functional
\begin{eqnarray*}
     I(u) &=&
  \frac{1}{m+1}  \int_0^T e^{2\omega t}\int_{\Omega}  \left(  |
       u(t,x)|^{m+1}dx +
        m    |(-\Delta)^{-1}
(f(x)-\omega u(t,x)- \frac{\partial u}{\partial t}(t,x))|^{\frac{m+1}{m}}\right) dx
         dt \nonumber \\
&& -\int_0^Te^{2\omega t}\int_\Omega
  u(x,t)((-\Delta)^{-1}f)(x)dx dt +\int_{\Omega}|\nabla (-\Delta)^{-1}u_{0}(x)|^{2}\, dx\nonumber\\
&&-2\int_{\Omega}u_0(x)(-\Delta)^{-1}u(0,x)\, dx +
\frac{1}{2}\left(
\|u(0)\|_{_{H^{-1}}}^2 +
e^{2T}\|u(T)\|_{_{H^{-1}}}^2\right)
   \end{eqnarray*}
   on  the space $A^2_H$ is equal to zero and is attained
uniquely at an $L^{m+1}(\Omega)$-valued path $u$ 
such that
$\int^T_0
\|\dot{u}(t)\|_H^2 dt < +\infty$ and which is a  solution of the
     equation:
\begin{equation}
\label{eqn:pme}
   \left\{ \begin{array}{lcl}
   \hfill \frac{\partial u}{\partial t}(t,x) &=& \Delta u^m +\omega u(t,x)+f  \hbox{\rm \, on \, $\Omega\times
[0,T]$}
\\
   u(0,x) &=& u_0 \quad \hbox{\rm on \quad $\Omega$. }\nonumber \\
       \end{array}  \right.
   \end{equation}

\subsection{Variational resolution for coupled flows and wave-type equations}
 Again, ASD Lagrangians are suited to treat variationally coupled evolution equations.
 
\begin{proposition} Let $\phi$ be a proper convex lower semi-continuous function  on $X\times Y$ and let $A:X\to Y^{*}$ be any bounded linear operator. Assume $B_{1}:X\to X$ (resp., $B_{2}:Y\to Y$) are positive operators,  then for any $(x_{0}, y_{0})\in {\rm dom} (\partial \phi)$ and any $(f,g)\in X^*\times Y^*$, there exists a path $(x(t), y(t)) \in A^{2}_{X}\times A^{2}_{Y}$ such that
 \begin{eqnarray*}
-\dot x(t)-A^{*} y(t)-B_{1} x (t) +f &\in& \partial_{1} \phi ( x(t),  y(t))\\
-\dot y(t)+A x(t)-B_{2} y (t) + g&\in& \partial_{2} \phi ( x(t),  y(t))\\
x(0)&=&x_{0}\\
y(0)&=&y_{0}.
\end{eqnarray*}
The solution is obtained as a minimizer on $A^{2}_{X}\times A^{2}_{Y}$ of the following functional 
 \begin{eqnarray*}
I(x,y)&=&\int_{0}^{T}\left\{\psi (x(t), y(t))+\psi^{*}(-A^{*}y(t)-B^{a}_{1}x(t)-\dot x(t), Ax(t)-B^{a}_{2}y(t)-\dot y (t))\right\} dt\\
&& +\frac{1}{2}\|x(0)\|^{2} -2\langle x_0, x(0)\rangle +\|x_0\|^{2} +\frac{1}{2}\|x(T)\|^{2}\\
&& +\frac{1}{2}\|y(0)\|^{2} -2\langle y_0, y(0)\rangle +\|y_0\|^{2} +\frac{1}{2}\|y(T)\|^{2}.
\end{eqnarray*}
whose infimum is zero. Here $B_{1}^{a}$ (resp., $B_{2}^{a}$) are the skew-symmetric parts of $B_{1}$ and $B_{2}$ and 
\[
\psi (x,y)=\phi (x,y)+\frac{1}{2} \langle B_{1}x, x\rangle  -\langle f, x\rangle+ \frac{1}{2} \langle B_{2}y, y\rangle   -\langle g, x\rangle
\]
 \end{proposition} 
\noindent {\bf Proof:} It is enough to apply Theorem 7.1 to the ASD Lagrangian
 \[
L((x,y), (p,q))=\psi(x, y)+\psi^{*}(-A^{*}y-B^{a}_{1}x-p, Ax-B^{a}_{2}y-q).
\]
obtained by shifting to the right the ASD Lagrangian $\phi\oplus_{\rm as} A$ by the skew-adjoint operator $(B^{a}_{1}, B^{a}_{2})$. If $(\bar x(t), \bar y(t))$ is where the infimum is attained, then we get
\begin{eqnarray*}
0&=&I(\bar x,\bar y)\\
&=&\int_{0}^{T}\big\{\psi(\bar x(t), \bar y(t))+
\psi^{*}(-A^{*}{\bar y}(t)-B^{a}_{1}{\bar x}(t)-{\dot {\bar x}}(t), A{\bar x}(t)-B^{a}_{2}{\bar y}(t)-{\dot {\bar  y}}(t)) \\
&& -\langle ({\bar x}(t), {\bar y}(t)),  (-A^{*}{\bar y}(t)-B^{a}_{1}{\bar x}(t)-{\dot {\bar x}}(t), A{\bar x}(t)-B^{a}_{2}{\bar y}(t)-{\dot {\bar  y}}(t))\rangle\big\}dt\\
&&+\|x(0)-x_{0}\|^{2}++\|y(0)-y_{0}\|^{2}
 \end{eqnarray*}
It follows that $\bar x (0)=x_{0}$, $\bar y (0)=0$ and the integrand is zero for almost all $t$ which yields 
\begin{eqnarray*}
-\dot x(t)-A^{*} y(t)-B^{a}_{1} x (t) &\in& \partial_{1} \psi ( x(t),  y(t)) =\partial_{1} \phi ( x(t),  y(t))+B^{s}_{1}x(t)-f\\
-\dot y(t)+A x(t)-B^{a}_{2} y (t) &\in& \partial_{2} \psi ( x(t),  y(t)) = \partial_{2} \phi ( x(t),  y(t))+B_{2}^{s}y(t)-g\\
x(0)&=&x_{0}\\
y(0)&=&y_{0}.
\end{eqnarray*}

Consider now two convex lower semi-continuous $\phi_{1}$ and $\phi_{2}$ on Hilbert spaces $X$ and $Y$ respectively, as well as two positive operators $B_{1}$ on $X$ and $B_{2}$ on $Y$. For any $(f, g)\in X\times Y$, consider the convex functionals $\psi_{1} (x)=\frac{1}{2} \langle B_{1}x, x\rangle +\varphi_{1} (x)$ and  $\psi_{2} (x)=\frac{1}{2} \langle B_{2}x, x\rangle +\varphi_{2} (x)$, and  the anti-selfdual Lagrangians 
\[
L(x, p)= \psi_{1} (x) -\langle f, x\rangle +\psi_{1}^{*}(-B_{1}^{a}x +f -p),\quad \hbox{\rm for $(x,p)\in X\times X$},
\]
and
\[
M(y, q)= \psi_{2} (y)-\langle g, y\rangle+\psi_{2}^{*}(-B_{2}^{a}y +g -q),\quad \hbox{\rm for $(y,q)\in Y\times Y$},
\] 
For  $w, w' \in {\bf R}$, we associate the following time-dependent ASD Lagrangian:
\[
L_{\omega}(t, x, p):=e^{-2wt}L(e^{wt}x, e^{wt}p)\quad{\rm and} \quad M_{\omega'}(t, y, q)= e^{-2w't}M(e^{w't}y, e^{w't}q).
 \]
 Let $A:X\to Y$ be any bounded linear operator and consider for any $c\in \R$ the following twisted ASD Lagrangian on $X\times Y$
\[ 
(L_{\omega}\oplus_{c^{2}A}M_{\omega'}) (t, (x,y), (p.q)):=L_{\omega}(t, x, A^{*}y+p)+M_{\omega'}(t, y, -c^{2}Ax+q). 
\]
where the duality in $X\times Y$ is given by $\langle (x,y), (p,q)\rangle =\langle x,p\rangle +c^{-2}\langle y, q\rangle$.  
  Applying Theorem 7.2, we obtain 
\begin{proposition} Assume $0\in {\rm Dom} (\partial \phi_{1})$ and $0\in {\rm Dom} (\partial \phi_{2})$, and consider the following functional on  $A^2_X\times A^{2}_{Y}$:
  \begin{eqnarray*}
 I(u,v)&=&\int_0^Te^{-2\omega t}\left\{\psi_{1} (e^{\omega t} u(t))+\psi_{1}^*(e^{\omega t}(-A^{*}v(t)-B_{1}^au(t) -{\dot u(t))}\right\}dt\\
&&+\int_0^Te^{-2\omega' t}\left\{\psi_{2} (e^{\omega' t} v(t))+\psi_{2}^*(e^{\omega' t}(c^{2}Au(t)-B_{2}^av(t) -{\dot v(t))}\right\}dt\\
 && +\frac{1}{2}\|u(0)\|^{2} -2\langle x_0, u(0)\rangle +\|x_0\|^{2} +\frac{1}{2}\|u(T)\|^{2}\\
 && +\frac{1}{2}\|v(0)\|^{2} -2\langle y_0, v(0)\rangle +\|y_0\|^{2} +\frac{1}{2}\|v(T)\|^{2}.
 \end{eqnarray*}
The minimum of $I$ is then zero and is attained at a path $(\bar x(t), \bar y(t)$, in such a way that  $x(t)= e^{\omega t} \bar x(t)$ and  $y(t)= e^{\omega' t} \bar y(t)$ form a solution of  the system of equations 
      \begin{equation}
 \left\{ \begin{array}{lcl}
\label{eqn:existence}
-\dot x(t)+\omega x(t)-A^{*}y(t)-B_{1}x(t)+f &\in& \partial \varphi_{1} (x(t))\\
\hfill -\dot y(t)+\omega' y(t)+c^{2}Ax(t)-B_{2}y(t)+g&\in& \partial \varphi_{2} (y(t))\\
\hfill x(0)&=&x_{0}\\
\hfill y(0)&=&y_{0}.
\end{array}\right.
 \end{equation}
   \end{proposition}
 
\subsubsection*{Example 12: A variational principle for coupled equations}
Let ${\bf b_{1}}:\Omega \to  {\bf R^{n}}$ and  ${\bf b_{2}}:\Omega \to  {\bf R^{n}}$ be two smooth vector fields on a bounded domain $\Omega$ of  $\bf R^{n}$, verifying the conditions in example 3 and consider their corresponding first order linear operator $B_{1}v={\bf b_{1}}\cdot \nabla v  \hbox{\rm \, and  $B_{2}v={\bf b_{2}}\cdot \nabla v$.}$ Consider the Dirichlet problem:
\begin{equation}
\label{Ex1.5000}
 \left\{ \begin{array}{lcl}
    \hfill -\frac{\partial u}{\partial t}-Ê\Delta (v- u) +{\bf b_{1}}\cdot \nabla u &=& |u^{p-2}|u +f  \hbox{\rm \, on \, $(0, T]\times\Omega$}
\\
  \hfill -\frac{\partial v}{\partial t}+Ê\Delta (v+c^{2} u)  +{\bf b_{2}}\cdot \nabla v&=& |v^{q-2}|v +g  \hbox{\rm \, on \, $(0, T]\times\Omega$}\\
 \hfill  u(t,x)=v(t,x) &=& 0 \quad \quad \quad \quad  \hbox{\rm on  $(0, T]\times\partial \Omega$. }\nonumber \\
   \hfill  u(0,x)&=& u_{0}(x) \quad \quad \quad \quad  \hbox{\rm for  $x\in  \Omega$. }\nonumber\\
    \hfill  v(0,x)&=& v_{0}(x) \quad \quad \quad \quad  \hbox{\rm for  $x\in  \Omega$. }\nonumber
       \end{array}  \right.
   \end{equation}
 We can use the above to get
 \begin{theorem} Assume  ${\rm div }({\bf b_{1}})\geq 0$ and ${\rm div }({\bf b_{2}})\geq 0$ on $\Omega$, $1<p,q\leq \frac{n+2}{n-2}$  and consider on $A^{2}_{_{H^{1}_{0}(\Omega))}}\times A^{2}_{_{H^{1}_{0}(\Omega)}}$ 
 the functional 
\begin{eqnarray*}
I(u, v)&=&\int_{0}^{T}\left\{\Psi (u(t)) +\Psi^{*}({\bf b_{1}}.\nabla u(t) +\frac{1}{2}{\rm div }({\bf b_{1}}) \,  u(t)-\Delta v (t) -\dot u(t))\right\}dt\\
&& +\int_{0}^{T}\left\{\Phi (v(t))+\Phi^{*}({\bf b_{2}}.\nabla v (t)+\frac{1}{2}{\rm div }({\bf b_{2}}) \,   v(t) +c^{2}\Delta u (t)-\dot v(t))\right\}dt\\
&&+\int_{\Omega}\left\{\frac{1}{2}(|u(0,x)|^{2}+|u(T,x)|^{2}) -2u(0,x)u_{0}(x) +|u_{0}(x)|^{2}\right\}dx\\
&&+\int_{\Omega}\left\{\frac{1}{2}(|v(0,x)|^{2}+|v(T,x)|^{2}) -2v(0,x)v_{0}(x) +|v_{0}(x)|^{2}\right\}dx
\end{eqnarray*}
where
\[
 \Psi (u)=\frac{1}{2} \int_{\Omega} | \nabla  u |^{2}dx +\frac{1}{p}\int_{\Omega}|u|^{p} dx+\int_{\Omega}fu dx +\frac{1}{4}\int_{\Omega}{\rm div }({\bf b_{1}}) \, |u|^{2} dx, 
 \]
 \[
 \Phi (v)=\frac{1}{2} \int_{\Omega} | \nabla  v |^{2}dx +\frac{1}{q}\int_{\Omega}|v|^{q} dx+\int_{\Omega}gv dx +\frac{1}{4}\int_{\Omega}{\rm div }({\bf b_{2}}) \, |v|^{2} dx
 \]
and $\Psi^{*}$ and $\Phi^{*}$ are their Legendre transforms. Then there exists 
$(\bar u, \bar v)\in A^{2}_{_{H^{1}_{0}(\Omega))}}\times A^{2}_{_{H^{1}_{0}(\Omega)}}$ such that:
\[
I(\bar u, \bar v)=\inf \{I(u, v);  (u, v) \in A^{2}_{_{H^{1}_{0}(\Omega))}}\times A^{2}_{_{H^{1}_{0}(\Omega)}}\}=0, 
\]
and $(\bar u, \bar v)$ is a solution of $(\ref{Ex1.5000})$.
 \end{theorem}
   
\subsubsection*{Example 13: Pressureless gaz of sticky particles}

 Motivated by the recent work of  Brenier \cite{Br} we consider equations of the form   
 \begin{equation}
\partial_{tt}X=c^{2}\partial_{yy}X-\partial_{t}\partial_{a}\mu, \quad \partial_{a}X\geq 0,\quad \mu \geq 0.
\end{equation}
where here $X(t):=X(t,a, y)$ is a function  on $K=[0,1]\times {\bf R}/{\bf Z}$, and  $\mu (t, a, y)$ is a nonnegative measure that plays the role of a Lagrange multiplier for the constraint $\partial_{a}X\geq 0$. Following Brenier, we reformulate the problem with the following system:
 \begin{equation}
 \left\{ \begin{array}{lcl}
\label{eqn:existence}
-\dot X(t)-\frac {\partial U}{\partial y}(t) &\in& \partial \varphi_{1} (X(t))\\
\hfill -\dot U(t)+\frac {\partial X}{\partial y}(t) &=& 0\\
\hfill X(0)&=&X_{0}\\
\hfill U(0)&=&U_{0}
\end{array}\right.
 \end{equation}
where $\phi_{1}$ is the convex function defined on $L^{2}(K)$ by 
\begin{equation}
\varphi_{1} (X)=
\left\{ \begin{array}{lcl}
 0
& {\rm if} \quad  \partial_{a}X\geq 0
 \\
    +\infty   &{\rm elsewhere}
\end{array}\right.
\end{equation}
We can solve this system with the above method by first setting $\phi_{2}(U)=0$ for every $U\in L^{2}(K)$ and by considering the Hilbert spaces
$X=Y= H^{2}_{per}(K)$ to be the subspace of $A^{2}_{K}$ consisting of functions that are periodic in $y$. Define on this space the operator  $AX=\frac {\partial X}{\partial y}$ in such a way that $A^{*}=-A$. We consider now the functional 
 \begin{eqnarray*}
 I(X,U)&=&\int_0^T\left\{\phi_{1} (X(t))+\phi_{1}^*(-\frac {\partial U}{\partial y}(t) -{\dot X(t))}\right\}dt\\
&&+\int_0^T\left\{\phi_{2}^*(\frac {\partial X}{\partial y}(t) -{\dot U(t))}\right\}dt\\
 && +\frac{1}{2}\|X(0)\|^{2} -2\langle X_0, X(0)\rangle +\|X_0\|^{2} +\frac{1}{2}\|X(T)\|^{2}\\
 && +\frac{1}{2}\|U(0)\|^{2} -2\langle Y_0, U(0)\rangle +\|Y_0\|^{2} +\frac{1}{2}\|U(T)\|^{2}
 \end{eqnarray*}
 If $(X_{0}, U_{0})$ are such that $\partial_{a}X_{0}\geq 0$, then 
 the minimum of $I$ is then zero and is attained at a path $(\bar X(t), \bar U(t)$, which solves  the above system of equations. 

\section{Variational resolution of certain implicit PDEs}

Motivated by the time dependent case, we briefly describe in this section how ASD Lagrangians can be used to solve variationally certain types of implicit differential equations. Indeed, letting $\Omega$ be a bounded domain in ${\R}^{n}$, we let $\mu=\mu^{+}-\mu^{-}$ be a signed finite measure on the boundary $\partial \Omega$. We consider an equation of the form
 \begin{eqnarray}
 \label{last.eq.0}
\hfill L(x, u(x), \Lambda u (x)) + \langle u(x), (\Lambda u)(x)\rangle_{_{X, X^{*}}} &=&0\quad \quad \quad \quad \quad {\rm on} \quad \Omega \\ 
 u&=& 0\quad  \mu^{-}\, {\rm a.e.\quad on} \quad \Sigma.\nonumber
 \end{eqnarray}
 where $L:\Omega \times X \times X^{*}\to \R$ for some Banach space $X$ and where $\Lambda$ is an operator from a space $H(\Omega; X)$ of $X$-valued functions  on $\Omega$ to a space $K(\Omega; X^{*})$ of $X^{*}$-valued functions  on $\Omega$. 
 
 Suppose now that for each $x\in \Omega$, the Lagrangian $L(x,\cdot, \cdot)$ is anti-selfdual on $X\times X^{*}$ and suppose that $\Lambda$ is skew-adjoint modulo the boundary $\Sigma$ in the following sense: \\
 There exists   a coercive map $J:X\to X^{*}$ (normally a duality map) such that for any $u,v$ in $H(\Omega; X)$ we have:
 \[
 \int_{\Omega}\langle (\Lambda u) (x), v(x)\rangle_{_{X, X^{*}}}dx=-\int_{\Omega}\langle (\Lambda v) (x), u(x)\rangle_{_{X, X^{*}}}dx+\int_{\partial \Omega}\langle Jv (x), u(x)\rangle_{_{X, X^{*}}}d\mu.
    \]
 On can then consider the Lagrangian 
\[
{\cal L}(u,p)=\int_{\Omega }L(x, u(x), \Lambda u (x) +p(x)) dx +\frac{1}{2}\int_{\partial \Omega}\langle Ju (x), u(x)\rangle_{_{X, X^{*}}}d\mu^{+}+\frac{1}{2}\int_{\partial \Omega}\langle Ju (x), u(x)\rangle_{_{X, X^{*}}}d\mu^{-}
\]
on $H(\Omega; X)\times K(\Omega; X^{*})$. Under the right conditions, the Lagrangian $\cal L$ has every chance to be anti-selfdual and therefore Theorem 4.1 can apply to yield that the infimum of the functional 
\[
I(u):={\cal L}(u,0)=\int_{\Omega }L(x, u(x), \Lambda u (x)) dx +\frac{1}{2}\int_{\partial \Omega}\langle Ju (x), u(x)\rangle_{_{X, X^{*}}}d\mu^{+}+\frac{1}{2}\int_{\partial \Omega}\langle Ju (x), u(x)\rangle_{_{X, X^{*}}}d\mu^{-}.
\]
is zero and that there exists $\bar u \in H(\Omega; X)$ such that 
\[
I(\bar u)=\inf\{ I (u);u\in H^{1}_{A}(\Omega)\}=0.
\]
It follows that
\begin{eqnarray*}
0&=&\int_{\Omega }L(x, u(x), \Lambda u (x)) dx +\frac{1}{2}\int_{\partial \Omega}\langle Ju (x), u(x)\rangle_{_{X, X^{*}}}d\mu^{+}+\frac{1}{2}\int_{\partial \Omega}\langle Ju (x), u(x)\rangle_{_{X, X^{*}}}d\mu^{-}\\
&=&\int_{\Omega }\left(L(x, u(x), \Lambda u (x)) +\langle u(x), (\Lambda u)(x)\rangle_{_{X, X^{*}}}\right) dx + \frac{1}{2}\int_{\partial \Omega}\langle Ju (x), u(x)\rangle_{_{X, X^{*}}}d\mu^{-}
\end{eqnarray*}
Since both integrands are non-negative and $J$ is coercive we deduce that $\bar u$ satisfies equation (\ref{last.eq.0}).

The key here is that ``pointwise'' the operator $\Lambda$ may not have any particular property but on the ``average'' it is skew-adjoint on the function space and this allows for the variational approach to apply. We illustrate the method on the following example. 

\subsubsection*{Example 14: An implicit transport equation}
Let $\Omega$ be a bounded domain in ${\R}^{n}$ and let ${\bf a}$ be a vector field as in Example (6) as well as its entrance set  $\Sigma_{-}$.  Consider the following ``toy'' equation for a function $u:\Omega \to \R$:
  \begin{eqnarray}
 \label{last.eq}
 \frac{1}{2}u^{2}+u+\frac{1}{2}(-{\bf a}\cdot \nabla u -1)^{2}-({\bf a}\cdot \nabla u)u&=&0\quad {\rm on}\quad \Omega \\ 
 u&=& 0 \quad {\rm on \quad \Sigma_{-}}.\nonumber
 \end{eqnarray}
Assuming that ${\rm div} {\bf a}=0$, we can solve the above equation by minimizing the following functional 
\[
I(u)=\int_{\Omega} \left( \frac{1}{2}u^{2}+u+\frac{1}{2}(-{\bf a}\cdot \nabla u -1)^{2}\right) dx+\frac{1}{2}\int_{\Sigma_{+}}|u(x)|^{2}{\bf n}(x)\cdot {\bf a}(x)\, d\sigma
-\frac{1}{2}\int_{\Sigma_{-}}|u(x)|^{2}{\bf n}(x)\cdot {\bf a}(x)\,  d\sigma 
\]
 on the space $H^{1}_{A}(\Omega)$. Indeed, we can rewrite
\begin{eqnarray}
\label{functional}
{I}(u)&=&\psi (u)+\psi^{*}(-\Lambda u)+\frac{1}{2}\int_{\Sigma_{+}}|u(x)|^{2}{\bf n}(x)\cdot {\bf a}(x)\, d\sigma
-\frac{1}{2}\int_{\Sigma_{-}}|u(x)|^{2}{\bf n}(x)\cdot {\bf a}(x)\,  d\sigma \nonumber  
\end{eqnarray}
 where $\psi$ is the convex functional  on $L^{2}(\Omega)$ defined by:
\[
\psi (u)=\frac{1}{2}\int_{\Omega}(|u(x)|^{2}+ u(x))dx
\]  
and $\psi^{*}$ is its Legendre conjugate
\[
\psi^{*}(v)=\frac{1}{2}\int_{\Omega}(v(x)-1)^{2} dx .
\]
Note that $\Lambda u ={\bf a}\cdot \nabla u$ is here skew-adjoint modulo the boundary. 
There exists then a solution $\bar u$ for (\ref{last.eq}) that is obtained as a minimizer of the problem:
\[
I(\bar u)=\inf\{ I (u);u\in H^{1}_{A}(\Omega)\}=0.
\]
Rewrite now 
\begin{eqnarray*}
0&=&I(\bar u)=\int_{\Omega} \left( \frac{1}{2}\bar u^{2}+\bar u+\frac{1}{2}(-{\bf a}\cdot \nabla \bar u -1)^{2}\right) dx+\frac{1}{2}\int_{\Sigma_{+}}|\bar u|^{2}{\bf n}\cdot {\bf a}\, d\sigma
-\frac{1}{2}\int_{\Sigma_{-}}|\bar u|^{2}{\bf n}\cdot {\bf a}\,  d\sigma\\ 
&=&\int_{\Omega}\left(\frac{1}{2}\bar u^{2}+\bar u+\frac{1}{2}(-{\bf a}\cdot \nabla \bar u -1)^{2}-({\bf a}\cdot \nabla \bar u)\bar u\right) dx +
\int_{\Sigma_{-}}|\bar u|^{2}|{\bf n}\cdot {\bf a}|\,  d\sigma.
\end{eqnarray*}
The equation follows because each integrand is non-negative.

\end{document}